\documentclass[a4paper]{amsart}
\usepackage[all]{xy}
\usepackage{amssymb}
\usepackage{amsopn}
\usepackage{amsmath}

\newcommand{\ie}{{\em i.e.}\ }
\newcommand{\cf}{{\em cf.}\ }
\newcommand{\eg}{{\em e.g.}\ }
\newcommand{\ko}{\: , \;}
\newcommand{\ul}[1]{\underline{#1}}

\numberwithin{equation}{subsection}
\newtheorem{theorem}{Theorem}
\numberwithin{theorem}{section}
\newtheorem{classification-theorem}[theorem]{Classification Theorem}
\newtheorem{decomposition-theorem}[theorem]{Decomposition Theorem}
\newtheorem{proposition-definition}[theorem]{Proposition-Definition}
\newtheorem{periodicity-conjecture}[theorem]{Periodicity Conjecture}
\newtheorem{lemma}[theorem]{Lemma}
\newtheorem{proposition}[theorem]{Proposition}
\newtheorem{corollary}[theorem]{Corollary}
\newtheorem{conjecture}[theorem]{Conjecture}
\newtheorem{question}[theorem]{Question}

\newcommand{\reminder}[1]{}

\newcommand{\opname}[1]{\operatorname{\mathsf{#1}}}

\renewcommand{\mod}{\opname{mod}\nolimits}
\newcommand{\coh}{\opname{coh}\nolimits}
\newcommand{\rep}{\opname{rep}\nolimits}

\newcommand{\Mod}{\opname{Mod}\nolimits}

\newcommand{\ind}{\opname{ind}\nolimits}
\newcommand{\per}{\opname{per}\nolimits}
\newcommand{\add}{\opname{add}\nolimits}

\newcommand{\Gr}{\opname{Gr}\nolimits}
\newcommand{\dimv}{\underline{\dim}\,}

\renewcommand{\rep}{\opname{rep}\nolimits}

\newcommand{\scr}{\mathcal}

\renewcommand{\ker}{\opname{ker}\nolimits}

\newcommand{\Z}{\mathbb{Z}}
\newcommand{\N}{\mathbb{N}}
\newcommand{\Q}{\mathbb{Q}}
\newcommand{\C}{\mathbb{C}}
\newcommand{\R}{\mathbb{R}}

\renewcommand{\P}{\mathbb{P}}

\newcommand{\la}{\leftarrow}
\newcommand{\iso}{\stackrel{_\sim}{\rightarrow}}
\newcommand{\id}{\mathbf{1}}

\newcommand{\Hom}{\opname{Hom}}

\newcommand{\RHom}{\opname{RHom}}

\newcommand{\Ext}{\opname{Ext}}

\newcommand{\End}{\opname{End}}
\newcommand{\rad}{\opname{rad}}
\newcommand{\irr}{\opname{irr}}

\newcommand{\ten}{\otimes}

\newcommand{\Tor}{\opname{Tor}}

\newcommand{\hh}{\opname{HH}}
\newcommand{\ca}{{\mathcal A}}

\newcommand{\cc}{{\mathcal C}}
\newcommand{\cd}{{\mathcal D}}

\newcommand{\cp}{{\mathcal P}}

\newcommand{\ct}{{\mathcal T}}
\newcommand{\cu}{{\mathcal U}}

\newcommand{\cx}{{\mathcal X}}

\newcommand{\eps}{\varepsilon}
\renewcommand{\phi}{\varphi}
\newcommand{\del}{\partial}

\renewcommand{\hat}[1]{\widehat{#1}}

\newcommand{\bt}{\bullet}

\newcommand{\sgn}{\mbox{sgn}}
\newcommand{\sign}{\mbox{sign}}
\newcommand{\Spec}{\mbox{Spec}}

\renewcommand{\tilde}[1]{\widetilde{#1}}

\newcommand{\arr}[1]{\stackrel{#1}{\rightarrow}}

\begin{document}

\date{July 2008, last modified on \today}

\title[Cluster algebras, representations,
triangulated categories]{Cluster algebras,
quiver representations and triangulated categories}
\author{Bernhard Keller}
\address{Universit\'e Paris Diderot -- Paris 7\\
    UFR de Math\'ematiques\\
   Institut de Math\'ematiques de Jussieu, UMR 7586 du CNRS \\
   Case 7012\\
   2, place Jussieu\\
   75251 Paris Cedex 05\\
   France }
\email{keller@math.jussieu.fr}

\begin{abstract}
  This is an introduction to some aspects of Fomin-Zele\-vinsky's
  cluster algebras and their links with the representation theory of
  quivers and with Calabi-Yau triangulated categories.  It is based on
  lectures given by the author at summer schools held in 2006
  (Bavaria) and 2008 (Jerusalem). In addition to by now classical
  material, we present the outline of a proof of the periodicity
  conjecture for pairs of Dynkin diagrams (details will appear
  elsewhere) and recent results on the interpretation of mutations as
  derived equivalences.
\end{abstract}

\maketitle
\tableofcontents

\section{Introduction}

\subsection{Context}
Cluster algebras were invented by S.~Fomin and A.~Zelevinsky
\cite{FominZelevinsky02} in the spring of the year 2000 in a project
whose aim it was to develop a combinatorial approach to the results obtained
by G.~Lusztig concerning total positivity in algebraic groups
\cite{Lusztig96}
on the one hand and canonical bases in quantum groups
\cite{Lusztig90} on the
other hand (let us stress that canonical bases were discovered
independently and simultaneously by M.~Kashiwara \cite{Kashiwara90}).
Despite great
progress during the last few years
\cite{FominZelevinsky03} \cite{BerensteinFominZelevinsky05}
\cite{FominZelevinsky07},
we are still relatively
far from these initial aims. Presently, the best results on the link between
cluster algebras and canonical bases are probably those of
C.~Geiss, B.~Leclerc and J.~Schr\"oer
\cite{GeissLeclercSchroeer05}
\cite{GeissLeclercSchroeer06}
\cite{GeissLeclercSchroeer08b}
\cite{GeissLeclercSchroeer07b}
\cite{GeissLeclercSchroeer08a}
but even they cannot construct canonical bases from cluster variables
for the moment. Despite these difficulties, the theory of
cluster algebras has witnessed spectacular growth thanks
notably to the many links that have been discovered with
a wide range of subjects including
\begin{itemize}
\item Poisson geometry
\cite{GekhtmanShapiroVainshtein03}
\cite{GekhtmanShapiroVainshtein05} \ldots ,
\item integrable systems \cite{FominZelevinsky03b} \ldots,
\item higher Teichm\"uller spaces
\cite{FockGoncharov09}
\cite{FockGoncharov07b}
\cite{FockGoncharov07a}
\cite{FockGoncharov06}
\ldots ,
\item combinatorics and the study of combinatorial
polyhedra like the Stasheff associahedra
\cite{ChapotonFominZelevinsky02}
\cite{Chapoton04}
\cite{Krattenthaler06}
\cite{FominReading05}
\cite{Musiker07}
\cite{FominShapiroThurston08}
\ldots ,
\item commutative and non commutative algebraic geometry,
in particular the study of stability
conditions in the sense of Bridgeland
\cite{Bridgeland06a} \cite{Bridgeland06}
\cite{Bridgeland07},
Calabi-Yau algebras
\cite{Ginzburg06} \cite{ChuangRouquier07},
Donaldson-Thomas invariants \cite{Szendroi08}
\cite{Kontsevich07a} \cite{Kontsevich07} \cite{KontsevichSoibelman08}
\ldots,
\item and last not least the representation theory of quivers and
finite-dimensional algebras, cf. for example the surveys
\cite{BuanMarsh06} \cite{Reiten06} \cite{Ringel07} .
\end{itemize}
We refer to the introductory
papers
\cite{FominZelevinsky03a}
\cite{Zelevinsky04}
\cite{Zelevinsky02}
\cite{Zelevinsky05}
\cite{Zelevinsky07a}
and to the cluster algebras portal \cite{Fomin07}
for more information on cluster algebras
and their links with other parts of mathematics.

The link between cluster algebras and quiver representations follows
the spirit of categorification: One tries to interpret cluster algebras
as combinatorial (perhaps $K$-theoretic) invariants associated with
categories of representations. Thanks to the rich structure of
these categories, one can then hope to prove results
on cluster algebras which seem beyond the scope of the purely combinatorial
methods. It turns out that the link becomes especially beautiful
if we use {\em triangulated categories} constructed from categories
of quiver representations.

\subsection{Contents}
We start with an informal presentation of
Fomin-Zelevinsky's classification theorem and of the cluster
algebras (without coefficients) associated with Dynkin diagrams.
Then we successively introduce quiver mutations, the cluster algebra
associated with a quiver, and the cluster algebra with coefficients
associated with an `ice quiver' (a quiver some of whose vertices
are frozen). We illustrate cluster algebras with coefficients on a number of examples
appearing as coordinate algebras of homogeneous varieties.

Sections~\ref{s:categorification-finite-case}, \ref{s:categorification-acyclic-case}
and \ref{s:categorification-2-Calabi-Yau} are devoted to the
(additive) categorification of cluster algebras. We start by
recalling basic notions from the representation theory of quivers.
Then we present a fundamental link between indecomposable representations
and cluster variables: the Caldero-Chapoton formula. After a brief reminder on
derived categories in general, we give the canonical presentation in
terms of generators and relations of the derived category
of a Dynkin quiver. This yields in particular a presentation
for the module category, which we use to sketch Caldero-Chapoton's
proof of their formula. Then we introduce the cluster category
and survey its many links to the cluster algebra in the
finite case. Most of these links are still valid, mutatis mutandis,
in the acyclic case, as we see in section~\ref{s:categorification-acyclic-case}.
Surprisingly enough, one can go even further and categorify
interesting classes of cluster algebras using generalizations
of the cluster category, which are still triangulated categories
and Calabi-Yau of dimension $2$. We present this relatively
recent theory in section~\ref{s:categorification-2-Calabi-Yau}.
In section~\ref{s:periodicity}, we apply it to sketch a proof of the
periodicity conjecture for pairs of Dynkin diagrams (details will
appear elsewhere \cite{Keller10a}). In the final
section~\ref{s:quiver-mutation-and-derived-equivalence}, we
give an interpretation of quiver mutation in terms of derived
equivalences. We use this framework to establish links
between various ways of lifting the mutation operation from
combinatorics to linear or homological algebra: mutation of
cluster-tilting objects, spherical collections and decorated
representations.

\subsection*{Acknowledgments}
These notes are based on lectures given at the IRTG-Summer\-school
2006 (Schloss Reisensburg, Bavaria) and at the Midrasha Mathematicae
2008 (Hebrew University, Jeru\-salem).  I thank the organizers of these
events for their generous invitations and for providing stimulating
working conditions. I am grateful to Thorsten Holm, Peter J\o rgensen
and Raphael Rouquier for their encouragment and for accepting to
include these notes in the proceedings of the `Workshop on
triangulated categories' they organized at Leeds in 2006. It is a
pleasure to thank to Carles Casacuberta, Andr\'e Joyal, Joachim Kock,
Amnon Neeman and Frank Neumann for an invitation to the Centre de
Recerca Matem\`atica, Barcelona, where most of this text was written
down. I thank Lingyan Guo, Sefi Ladkani and Dong Yang
for kindly pointing out misprints and inaccuracies.
I am indebted to Tom Bridgeland, Osamu Iyama, David
Kazhdan, Bernard Leclerc, Tomoki Nakanishi, Rapha\"el Rouquier
and Michel Van den Bergh for helpful conversations.

\section{An informal introduction to cluster-finite cluster algebras}
\label{s:informal-introduction}

\subsection{The classification theorem}\label{ss:classification-theorem}
Let us start with a remark on terminology: a {\em cluster} is a
group of similar things or people positioned or occurring closely
together \cite{Oxford04}, as in the combination `star cluster'.
In French, `star cluster' is translated as `amas d'\'etoiles', whence
the term `alg\`ebre amass\'ee' for cluster algebra.

We postpone the precise definition of a cluster algebra to
section~\ref{s:sym-cluster-alg-without-coeff}.
For the moment, the following description will suffice: A cluster algebra
is a commutative $\Q$-algebra endowed with a family of distinguished
generators (the {\em cluster variables}) grouped into overlapping
subsets (the {\em clusters}) of fixed finite cardinality, which are constructed
recursively using {\em mutations}.

The set of cluster variables in a cluster algebra may be finite
or infinite. The first important result of the theory is the
classification of those cluster algebras where it is finite:
the {\em cluster-finite} cluster algebras. This is the

\begin{classification-theorem}[Fomin-Zelevinsky \cite{FominZelevinsky03}]
The cluster-finite cluster algebras are para\-metrized by the
finite root systems (like semisimple complex Lie algebras).
\end{classification-theorem}

It follows that for each Dynkin diagram $\Delta$, there is a canonical
cluster algebra $\ca_\Delta$. It turns out that $\ca_\Delta$ occurs
naturally as a subalgebra of the field of rational functions $\Q(x_1,
\ldots, x_n)$, where $n$ is the number of vertices of $\Delta$.
Since $\ca_\Delta$ is generated by its cluster variables (like
any cluster algebra), it suffices to produce the (finite) list
of these variables in order to describe $\ca_\Delta$. Now
for the algebras $\ca_\Delta$, the recursive construction
via mutations mentioned above simplifies considerably. In fact,
it turns out that one can {\em directly construct the cluster
variables} without first constructing the clusters. This
is made possible by

\subsection{The knitting algorithm}\label{ss:knitting-algorithm}
The general algorithm will become clear from the following three
examples. We start with the simplest non trivial Dynkin diagram:
\[
\Delta=A_2 : \xymatrix{ \circ \ar@{-}[r] & \circ}\; .
\]
We first choose a numbering of its vertices and an orientation
of its edges:
\[
\vec{\Delta}=\vec{A}_2 : \xymatrix{ 1 \ar[r] & 2}\; .
\]
Now we draw the so-called {\em repetition} (or {\em Bratteli
diagram}) $\Z\vec{\Delta}$ associated with $\Delta$: We first draw
the product $\Z \times \vec{\Delta}$ made up of a countable number
of copies of $\Delta$ (drawn slanted upwards); then for each arrow
$\alpha: i \to j$ of $\Delta$, we add a new family of arrows
$(n,\alpha^*): (n,j) \to (n+1,i)$, $n\in \Z$ (drawn slanted
downwards). We refer to section~\ref{ss:der-cat-Dynkin-quiver}
for the formal definition.
Here is the result for $\vec{\Delta}=\vec{A}_2$:
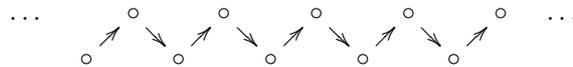
\[
\ldots\quad
\begin{xy} 0;<0.5pt,0pt>:<0pt,-0.5pt>::
(0,35) *+{\circ} ="0", (35,0) *+{\circ} ="1", (69,35) *+{\circ}
="2", (103,0) *+{\circ} ="3", (138,35) *+{\circ} ="4", (172,0)
*+{\circ} ="5", (207,35) *+{\circ} ="6", (241,0) *+{\circ} ="7",
(275,35) *+{\circ} ="8", (310,0) *+{\circ} ="9", "0", {\ar"1"}, "1",
{\ar"2"}, "2", {\ar"3"}, "3", {\ar"4"}, "4", {\ar"5"}, "5",
{\ar"6"}, "6", {\ar"7"}, "7", {\ar"8"}, "8", {\ar"9"},
\end{xy}
\quad
\ldots
\]
We will now assign a cluster variable to each vertex
of the repetition. We start by assigning $x_1$ and $x_2$ to the
vertices of the zeroth copy of $\vec{\Delta}$. Next, we construct
new variables $x'_1, x'_2, x''_1, \ldots$ by `knitting' from the
left to the right (an analogous procedure can be used to go from
the right to the left).
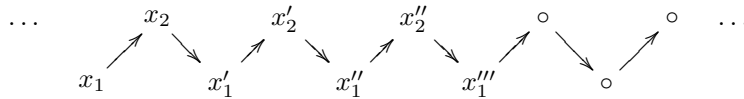
\[
\ldots\quad
\begin{xy} 0;<0.7pt,0pt>:<0pt,-0.7pt>::
(0,35) *+{x_1} ="0", (35,0) *+{x_2} ="1", (69,35) *+{x_1'}
="2", (103,0) *+{x_2'} ="3", (138,35) *+{x_1''} ="4", (172,0)
*+{x_2''} ="5", (207,35) *+{x_1'''} ="6", (241,0) *+{\circ} ="7",
(275,35) *+{\circ} ="8", (310,0) *+{\circ} ="9", "0", {\ar"1"}, "1",
{\ar"2"}, "2", {\ar"3"}, "3", {\ar"4"}, "4", {\ar"5"}, "5",
{\ar"6"}, "6", {\ar"7"}, "7", {\ar"8"}, "8", {\ar"9"},
\end{xy}
\quad
\ldots
\]
To compute $x_1'$, we consider its immediate predecessor $x_2$,
add $1$ to it and divide the result by the left translate of
$x_1'$, to wit the variable $x_1$. This yields
\[
x_1'=\frac{1+x_2}{x_1}.
\]
Similarly, we compute $x_2'$ by adding $1$ to its predecessor $x_1'$
and dividing the result by the left translate $x_2$:
\[
x_2'=\frac{1+x_1'}{x_2} = \frac{x_1+1+x_2}{x_1 x_2}.
\]
Using the same rule for $x_1''$ we obtain
\[
x_1''=\frac{1+x_2'}{x_1'} =
\left(\frac{x_1x_2+x_1+1+x_2}{x_1x_2}\right)/
\left(\frac{1+x_2}{x_1}\right) = \frac{1+x_1}{x_2}.
\]
Here something remarkable has happened: The numerator $x_1x_2 + x_1 + 1 + x_2$
is actually divisible by $1+x_2$ so that the denominator remains
a monomial (contrary to what one might expect). We continue with
\[
x_2''=\frac{1+x_1''}{x_2'} =
\left(\frac{x_2+1+x_1}{x_2}\right)/\left(\frac{x_1+1+x_2}{x_1x_2}\right)
= x_1 \ko
\]
a result which is perhaps even more surprising. Finally, we get
\[
x_1'''= \frac{1+x_2''}{x_1''} =(1+x_1)/\left(\frac{1+x_1}{x_2}\right) = x_2.
\]
Clearly, from here on, the pattern will repeat. We could have computed
`towards the left' and would have found the same repeating pattern.
In conclusion, there are the $5$ cluster variables $x_1, x_2, x_1', x_2'$
and $x_1''$ and the cluster algebra $\ca_{A_2}$ is the $\Q$-subalgebra
(not the subfield!) of $\Q(x_1, x_2)$ generated by these $5$ variables.

Before going on to a more complicated example, let us record the
remarkable phenomena we have observed:
\begin{itemize}
\item[(1)] All denominators of all cluster variables are monomials.
In other words, the cluster variables are Laurent polynomials.
This {\em Laurent phenomenon} holds for all cluster variables in
all cluster algebras, as shown by Fomin and Zelevinsky
\cite{FominZelevinsky02a}.
\item[(2)] The computation is {\em periodic} and thus only yields
finitely many cluster variables. Of course, this was to be expected
by the classification theorem above. In fact, the procedure generalizes
easily from Dynkin diagrams to arbitrary trees, and then periodicity
characterizes Dynkin diagrams among trees.
\item[(3)] {\em Numerology}: We have obtained $5$ cluster variables. Now
we have $5=2+3$ and this decomposition does correspond to a natural
partition of the set of cluster variables into the two {\em initial}
cluster variables $x_1$ and $x_2$ and the three non initial ones $x_1'$,
$x_2'$ and $x_1''$. The latter are in natural bijection with the positive roots
$\alpha_1$, $\alpha_1+\alpha_2$ and $\alpha_2$ of the root system of type $A_2$
with simple roots $\alpha_1$ and $\alpha_2$. To see this, it suffices to look
at the denominators of the three variables: The denominator $x_1^{d_1} x_2^{d_2}$
corresponds to the root $d_1 \alpha_1+d_2 \alpha_2$. It was proved by
Fomin-Zelevinsky \cite{FominZelevinsky03} that this generalizes to arbitrary
Dynkin diagrams. In particular, the number of cluster variables in
the cluster algebra $\ca_\Delta$ always equals the sum of the rank
and the number of positive roots of $\Delta$.
\end{itemize}
Let us now consider the example $A_3$: We choose the following linear
orientation:
\[
\xymatrix{1 \ar[r] & 2 \ar[r] & 3}.
\]
The associated repetition looks as follows:
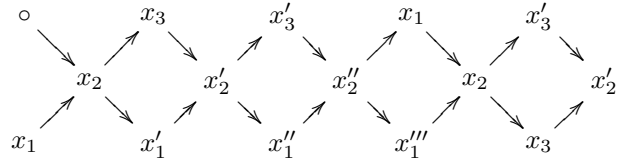
\[
\begin{xy} 0;<0.7pt,0pt>:<0pt,-0.7pt>::
(0,69) *+{x_1} ="0",
(0,0) *+{\circ} ="1",
(35,35) *+{x_2} ="2",
(69,69) *+{x_1'} ="3",
(69,0) *+{x_3} ="4",
(103,35) *+{x_2'} ="5",
(138,69) *+{x_1''} ="6",
(138,0) *+{x_3'} ="7",
(172,35) *+{x_2''} ="8",
(207,69) *+{x_1'''} ="9",
(207,0) *+{x_1} ="10",
(241,35) *+{x_2} ="11",
(275,69) *+{x_3} ="12",
(275,0) *+{x_3'} ="13",
(310,35) *+{x_2'} ="14",
"0", {\ar"2"},
"1", {\ar"2"},
"2", {\ar"3"},
"2", {\ar"4"},
"3", {\ar"5"},
"4", {\ar"5"},
"5", {\ar"6"},
"5", {\ar"7"},
"6", {\ar"8"},
"7", {\ar"8"},
"8", {\ar"9"},
"8", {\ar"10"},
"9", {\ar"11"},
"10", {\ar"11"},
"11", {\ar"12"},
"11", {\ar"13"},
"12", {\ar"14"},
"13", {\ar"14"},
\end{xy}
\]
The computation of $x_1'$ is as before:
\[
x_1'=\frac{1+x_2}{x_1} .
\]
However, to compute $x_2'$, we have to modify the rule, since $x_2'$ has
{\em two} immediate predecessors with associated variables $x_1'$ and $x_3$.
In the formula, we simply take the {\em product} over all immediate
predecessors:
\[
x_2'=\frac{1+x_1'x_3}{x_2} = \frac{x_1+x_3+x_2 x_3}{x_2 x_3}.
\]
Similarly, for the following variables $x_3'$, $x_1''$, \ldots.
We obtain the periodic pattern shown in the diagram above. In total,
we find $9=3+6$ cluster variables, namely
\begin{align*}
& x_1 \ko x_2 \ko x_3,
\frac{1+x_2}{x_1}\ko \frac{x_1+x_3+x_2x_3}{x_1x_2} \ko \frac{x_1+x_1 x_2 +x_3 + x_2 x_3}{x_1 x_2 x_3} \ko\\
& \frac{x_1+x_3}{x_2}\ko
\frac{x_1+x_1x_2+x_3}{x_2 x_3}\ko
\frac{1+x_2}{x_3}\;.
\end{align*}
The cluster algebra $\ca_{A_3}$ is the subalgebra of the
field $\Q(x_1, x_2, x_3)$ generated by these variables.
Again we observe that all denominators are monomials.
Notice also that $9=3+6$ and that $3$ is the rank of the root system associated
with $A_3$ and $6$ its number of positive roots. Moreover, if we look at
the denominators of the non initial cluster variables (those other than $x_1$, $x_2$, $x_3$), we
see a natural bijection with the positive roots
\[
\alpha_1, \alpha_1+\alpha_2, \alpha_1+\alpha_2+\alpha_3,
\alpha_2, \alpha_2+\alpha_3, \alpha_3
\]
of the root system of $A_3$, where $\alpha_1$, $\alpha_2$, $\alpha_3$ denote
the three simple roots.

Finally, let us consider the non simply laced Dynkin diagram $\Delta=G_2$:
\[
\xymatrix{ \circ \ar@{-}[r]^{(3,1)} & \circ}.
\]
The associated Cartan matrix is
\[
\left( \begin{array}{cc} 2 & -3 \\ -1 & 2 \end{array} \right)
\]
and the associated root system of rank $2$ looks as follows:
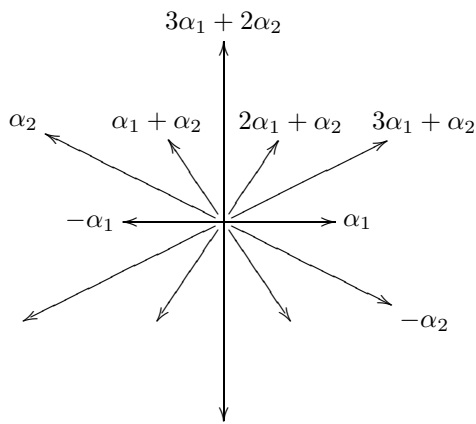
\[
\begin{xy} 0;<0.5pt,0pt>:<0pt,-0.5pt>::
(150,150) *+{} ="0",
(250,150) *+{\alpha_1} ="1",
(300,75) *+{3\alpha_1+\alpha_2} ="2",
(150,0) *+{3\alpha_1+2\alpha_2} ="3",
(200,75) *+{2\alpha_1+\alpha_2} ="4",
(100,75) *+{\alpha_1+\alpha_2} ="5",
(0,75) *+{\alpha_2} ="6",
(50,150) *+{-\alpha_1} ="7",
(0,225) *+{} ="8",
(100,225) *+{} ="9",
(200,225) *+{} ="10",
(150,300) *+{} ="11",
(300,225) *+{-\alpha_2} ="12",
"0", {\ar"1"},
"0", {\ar"2"},
"0", {\ar"3"},
"0", {\ar"4"},
"0", {\ar"5"},
"0", {\ar"6"},
"0", {\ar"7"},
"0", {\ar"8"},
"0", {\ar"9"},
"0", {\ar"10"},
"0", {\ar"11"},
"0", {\ar"12"},
\end{xy}
\]
We choose an orientation of the valued edge of $G_2$ to obtain
the following valued oriented graph:
\[
\vec{\Delta} : \xymatrix{ 1 \ar[r]^{(3,1)} & 2}.
\]
Now the repetition also becomes a valued oriented graph
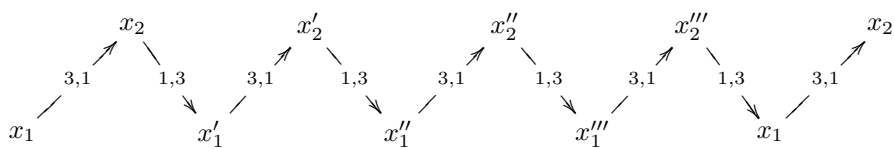
\[
\begin{xy} 0;<0.7pt,0pt>:<0pt,-0.7pt>::
(0,58) *+{x_1} ="0",
(59,0) *+{x_2} ="1",
(101,58) *+{x_1'} ="2",
(154,0) *+{x_2'} ="3",
(201,58) *+{x_1''} ="4",
(258,0) *+{x_2''} ="5",
(305,58) *+{x_1'''} ="6",
(358,0) *+{x_2'''} ="7",
(400,58) *+{x_1} ="8",
(459,0) *+{x_2} ="9",
"0", {\ar|*+{\scriptstyle 3,1}"1"},
"1", {\ar|*+{\scriptstyle 1,3}"2"},
"2", {\ar|*+{\scriptstyle 3,1}"3"},
"3", {\ar|*+{\scriptstyle 1,3}"4"},
"4", {\ar|*+{\scriptstyle 3,1}"5"},
"5", {\ar|*+{\scriptstyle 1,3}"6"},
"6", {\ar|*+{\scriptstyle 3,1}"7"},
"7", {\ar|*+{\scriptstyle 1,3}"8"},
"8", {\ar|*+{\scriptstyle 3,1}"9"},
\end{xy}
\]
The mutation rule is a variation on the one we are already
familiar with: In the recursion formula, each predecessor $p$ of
a cluster variable $x$ has to be raised to the power indicated
by the valuation `closest' to $p$. Thus, we have for example
\[
x_1' = \frac{1+x_2}{x_1} \ko
x_2'= \frac{1+(x_1')^3}{x_2} = \frac{1 + x_1^3 + 3x_2 + 3x_2^2 + x_2^3}{x_1^3x_2} \ko
x_1' = \frac{1+x_2'}{x_1'} = \frac{ \ldots} {x_1^2x_2} \quad \ko
\]
where we can read off the denominators from the decompositions of the
positive roots as linear combinations of simple roots given above. We find
$8=2+6$ cluster variables, which together generate the
cluster algebra $\ca_{G_2}$ as a subalgebra of $\Q(x_1, x_2)$.

\section{Symmetric cluster algebras without coefficients}
\label{s:sym-cluster-alg-without-coeff}

In this section, we will construct the cluster algebra associated
with an antisymmetric matrix with integer coefficients.
Instead of using matrices, we will use
quivers (without loops or $2$-cycles), since they are easy to visualize and well-suited
to our later purposes.

\subsection{Quivers}
Let us recall that a \emph{quiver} $Q$ is an oriented graph. Thus, it
is a quadruple given by a set $Q_0$ (the set of vertices), a set $Q_1$
(the set of arrows) and two maps $s:Q_1 \to Q_0$ and $t:Q_1\to Q_0$
which take an arrow to its source respectively its target. Our quivers
are `abstract graphs' but in practice we draw them as in this example:
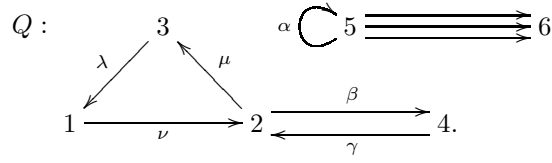
\[ Q:
\xymatrix{ & 3 \ar[ld]_\lambda & & 5 \ar@(dl,ul)[]^\alpha \ar@<1ex>[rr] \ar[rr] \ar@<-1ex>[rr] & & 6 \\
  1 \ar[rr]_\nu & & 2 \ar@<1ex>[rr]^\beta \ar[ul]_\mu & & 4.
  \ar@<1ex>[ll]^\gamma }
\]
A \emph{loop} in a quiver $Q$ is an arrow $\alpha$ whose source
coincides with its target; a \emph{$2$-cycle} is a pair of distinct
arrows $\beta\neq\gamma$ such that the source of $\beta$ equals the
target of $\gamma$ and vice versa. It is clear how to define
\emph{$3$-cycles}, \emph{connected components} \ldots . A quiver is
\emph{finite} if both, its set of vertices and its set of arrows, are
finite.

\subsection{Seeds and mutations} \label{ss:seeds}
Fix an integer $n\geq 1$.
A \emph{seed} is a pair $(R,u)$, where
\begin{itemize}
\item[a)] $R$ is a finite quiver without loops or $2$-cycles with vertex
set $\{1, \ldots, n\}$;
\item[b)] $u$ is a free generating set $\{u_1, \ldots, u_n\}$ of the field
$\Q(x_1, \ldots, x_n)$ of fractions of the polynomial ring $\Q[x_1, \ldots, x_n]$
in $n$ indeterminates.
\end{itemize}
Notice that in the quiver $R$ of a seed, all arrows between any two
given vertices point in the same direction (since $R$ does not have
$2$-cycles).  Let $(R,u)$ be a seed  and $k$ a vertex of $R$. The
\emph{mutation} $\mu_k(R,u)$ of $(R,u)$ at $k$ is the seed $(R',u')$,
where
\begin{itemize}
\item[a)] $R'$ is obtained from $R$ as follows:
\begin{itemize}
\item[1)] reverse all arrows incident with $k$;
\item[2)] for all vertices $i\neq j$ distinct from $k$, modify the number
of arrows between $i$ and $j$ as follows:
\[
\begin{array}{|c|c|}\hline
R & R'  \\ \hline
\xymatrix@=0.3cm{i \ar[rr]^{r} \ar[rd]_{s} & &  j \\
& k \ar[ru]_t & } &
\xymatrix@=0.3cm{i \ar[rr]^{r+st} & & j  \ar[ld]^t\\
& k \ar[lu]^{s}} \\\hline
\xymatrix@=0.3cm{i \ar[rr]^{r}  & &  j \ar[ld]^{t}\\
& k \ar[lu]^{s} & } &
\xymatrix@=0.3cm{i \ar[rr]^{r-st} \ar[rd]_{s} & &  j \\
& k \ar[ru]_t & } \\\hline
\end{array}
\]
where $r,s,t$ are non negative integers,
an arrow $\xymatrix@=0.3cm{i\ar[r]^l & j}$ with $l\geq 0$
means that $l$ arrows go from $i$ to $j$ and an arrow
$\xymatrix@=0.3cm{i\ar[r]^l & j}$ with $l\leq 0$
means that $-l$ arrows go from $j$ to $i$.
\end{itemize}
b) $u'$ is obtained from $u$ by replacing the element $u_k$ with
\begin{equation} \label{eq:exchange}
u_k'=\frac{1}{u_k} \left( \prod_{\mbox{\scriptsize arrows $i\to k$}} u_i + \prod_{\mbox{\scriptsize arrows $k\to j$}} u_j\right).
\end{equation}
\end{itemize}
In the {\em exchange relation}~(\ref{eq:exchange}),
if there are no arrows from $i$ with target $k$, the product is taken over
the empty set and equals $1$. It is not hard to see that $\mu_k(R,u)$ is
indeed a seed and that $\mu_k$ is an involution: we have $\mu_k(\mu_k(R,u))=(R,u)$.
Notice that the expression given in (\ref{eq:exchange}) for $u_k'$ is
{\em subtraction-free}.

To a quiver $R$ without loops or $2$-cycles with vertex set $\{1, \ldots, n\}$
there corresponds the $n\times n$ antisymmetric integer matrix $B$ whose
entry $b_{ij}$ is the number of arrows $i\to j$ minus the number of
arrows $j\to i$ in $R$ (notice that at least one of these numbers is
zero since $R$ does not have $2$-cycles). Clearly, this correspondence
yields a bijection. Under this bijection, the matrix $B'$ corresponding
to the mutation $\mu_k(R)$ has the entries
\[
b'_{ij} = \left\{ \begin{array}{ll}
-b_{ij} & \mbox{if $i=k$ or $j=k$;} \\
b_{ij}+\sgn(b_{ik})[b_{ik} b_{kj}]_+ & \mbox{else,}
\end{array} \right.
\]
where $[x]_+= \max(x,0)$. This is matrix mutation as it was defined
by Fomin-Zelevinsky in their seminal paper \cite{FominZelevinsky02},
\cf also \cite{FominZelevinsky07}.

\subsection{Examples of seed and quiver mutations}
Let $R$ be the cyclic quiver
\begin{equation} \label{quiver1}
\begin{xy} 0;<0.3pt,0pt>:<0pt,-0.3pt>::
(94,0) *+{1} ="0",
(0,156) *+{2} ="1",
(188,156) *+{3} ="2",
"1", {\ar"0"},
"0", {\ar"2"},
"2", {\ar"1"},
\end{xy}
\end{equation}
and $u=\{x_1, x_2, x_3\}$. If we mutate at $k=1$, we obtain the quiver
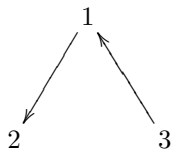
\[
\begin{xy} 0;<0.3pt,0pt>:<0pt,-0.3pt>::
(92,0) *+{1} ="0",
(0,155) *+{2} ="1",
(188,155) *+{3} ="2",
"0", {\ar"1"},
"2", {\ar"0"},
\end{xy}
\]
and the set of fractions given by $u'_1=(x_2+x_3)/x_1$, $u'_2=u_2=x_2$ and $u'_3=u_3=x_3$.
Now, if we mutate again at $1$, we obtain the original seed. This is a
general fact: Mutation at $k$ is an involution. If, on the other hand, we
mutate $(R', u')$ at $2$, we obtain the quiver
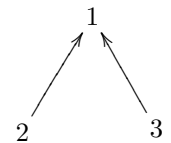
\[
\begin{xy} 0;<0.3pt,0pt>:<0pt,-0.3pt>::
(87,0) *+{1} ="0",
(0,145) *+{2} ="1",
(167,141) *+{3} ="2",
"1", {\ar"0"},
"2", {\ar"0"},
\end{xy}
\]
and the set $u''$ given by $u''_1=u'_1=(x_2+x_3)/x_1$,
$u'_2=\frac{x_1 + x_2 + x_3}{x_1 x_2}$ and
$u''_3=u'_3=x_3$.

An important special case of quiver mutation is the mutation at
a source (a vertex without incoming arrows) or a sink (a vertex without
outgoing arrows). In this case, the mutation only reverses the arrows
incident with the mutating vertex. It is easy to see that all orientations
of a tree are mutation equivalent and that only sink and source mutations
are needed to pass from one orientation to any other.

Let us consider the following, more complicated quiver
glued together from four $3$-cycles:
\begin{equation} \label{quiver2}
\begin{xy} 0;<0.4pt,0pt>:<0pt,-0.4pt>::
(74,0) *+{1} ="0",
(38,62) *+{2} ="1",
(110,62) *+{3} ="2",
(0,123) *+{4} ="3",
(74,104) *+{5} ="4",
(148,123) *+{6.} ="5",
"1", {\ar"0"},
"0", {\ar"2"},
"2", {\ar"1"},
"3", {\ar"1"},
"1", {\ar"4"},
"4", {\ar"2"},
"2", {\ar"5"},
"4", {\ar"3"},
"5", {\ar"4"},
\end{xy}
\end{equation}
If we successively perform mutations at the vertices $5$, $3$, $1$ and $6$,
we obtain the sequence of quivers (we use \cite{KellerQuiverMutationApplet})
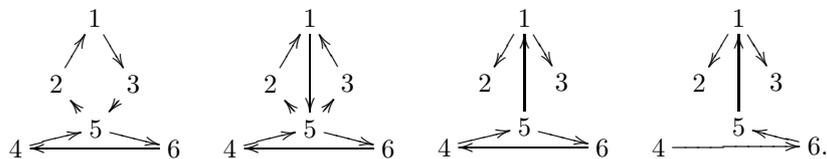
\[
\quad
\begin{xy} 0;<0.4pt,0pt>:<0pt,-0.4pt>::
(74,0) *+{1} ="0",
(38,62) *+{2} ="1",
(110,62) *+{3} ="2",
(0,123) *+{4} ="3",
(75,104) *+{5} ="4",
(148,123) *+{6} ="5",
"1", {\ar"0"},
"0", {\ar"2"},
"4", {\ar"1"},
"2", {\ar"4"},
"3", {\ar"4"},
"5", {\ar"3"},
"4", {\ar"5"},
\end{xy}
\quad
\begin{xy} 0;<0.4pt,0pt>:<0pt,-0.4pt>::
(75,0) *+{1} ="0",
(38,62) *+{2} ="1",
(110,61) *+{3} ="2",
(0,123) *+{4} ="3",
(75,104) *+{5} ="4",
(148,123) *+{6} ="5",
"1", {\ar"0"},
"2", {\ar"0"},
"0", {\ar"4"},
"4", {\ar"1"},
"4", {\ar"2"},
"3", {\ar"4"},
"5", {\ar"3"},
"4", {\ar"5"},
\end{xy}
\quad
\begin{xy} 0;<0.4pt,0pt>:<0pt,-0.4pt>::
(75,0) *+{1} ="0",
(38,61) *+{2} ="1",
(110,60) *+{3} ="2",
(0,122) *+{4} ="3",
(75,103) *+{5} ="4",
(148,122) *+{6} ="5",
"0", {\ar"1"},
"0", {\ar"2"},
"4", {\ar"0"},
"3", {\ar"4"},
"5", {\ar"3"},
"4", {\ar"5"},
\end{xy}
\quad
\begin{xy} 0;<0.4pt,0pt>:<0pt,-0.4pt>::
(75,0) *+{1} ="0",
(38,61) *+{2} ="1",
(110,60) *+{3} ="2",
(0,122) *+{4} ="3",
(75,103) *+{5} ="4",
(149,121) *+{6.} ="5",
"0", {\ar"1"},
"0", {\ar"2"},
"4", {\ar"0"},
"3", {\ar"5"},
"5", {\ar"4"},
\end{xy}
\quad
\]
Notice that the last quiver no longer has any oriented cycles
and is in fact an orientation of the Dynkin diagram of type $D_6$.
The sequence of new fractions appearing in these steps is
\begin{eqnarray*}
u'_5 & = &\frac{x_3 x_4 + x_2 x_6}{x_5} \ko \quad
u'_3  = \frac{x_3 x_4 + x_1 x_5 + x_2 x_6}{x_3 x_5}\ko \\
u'_1 & = & \frac{x_2 x_3 x_4 + x_3^2 x_4 + x_1 x_2 x_5 + x_2^2 x_6 + x_2 x_3x_6}{x_1 x_3 x_5} \ko\quad
u'_6=\frac{x_3 x_4 + x_4 x_5 + x_2 x_6}{x_5 x_6}\;.
\end{eqnarray*}
It is remarkable that all the denominators appearing here are monomials
and that all the coefficients in the numerators are positive.

Finally,
let us consider the quiver
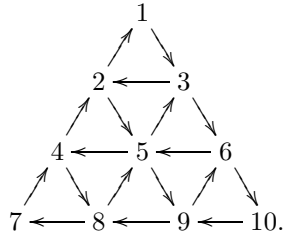
\begin{equation} \label{quiver3}
\begin{xy} 0;<0.6pt,0pt>:<0pt,-0.6pt>::
(79,0) *+{1} ="0",
(52,44) *+{2} ="1",
(105,44) *+{3} ="2",
(26,88) *+{4} ="3",
(79,88) *+{5} ="4",
(131,88) *+{6} ="5",
(0,132) *+{7} ="6",
(52,132) *+{8} ="7",
(105,132) *+{9} ="8",
(157,132) *+{10.} ="9",
"1", {\ar"0"},
"0", {\ar"2"},
"2", {\ar"1"},
"3", {\ar"1"},
"1", {\ar"4"},
"4", {\ar"2"},
"2", {\ar"5"},
"4", {\ar"3"},
"6", {\ar"3"},
"3", {\ar"7"},
"5", {\ar"4"},
"7", {\ar"4"},
"4", {\ar"8"},
"8", {\ar"5"},
"5", {\ar"9"},
"7", {\ar"6"},
"8", {\ar"7"},
"9", {\ar"8"},
\end{xy}
\end{equation}
One can show \cite{KellerReiten06} that it is impossible to transform it
into a quiver without
oriented cycles by a finite sequence of mutations. However, its mutation
class (the set of all quivers obtained from it by iterated mutations)
contains many quivers with just one oriented cycle, for example
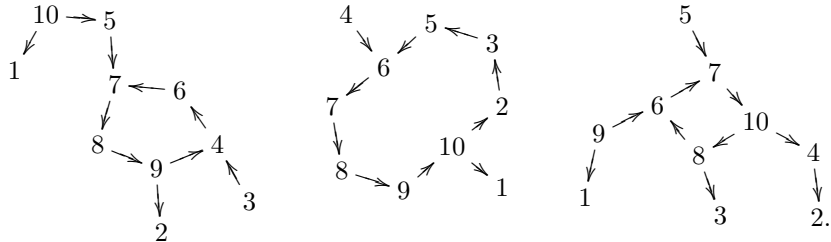
\[
\begin{xy} 0;<0.3pt,0pt>:<0pt,-0.3pt>::
(0,70) *+{1} ="0",
(183,274) *+{2} ="1",
(293,235) *+{3} ="2",
(253,164) *+{4} ="3",
(119,8) *+{5} ="4",
(206,96) *+{6} ="5",
(125,88) *+{7} ="6",
(104,164) *+{8} ="7",
(177,194) *+{9} ="8",
(39,0) *+{10} ="9",
"9", {\ar"0"},
"8", {\ar"1"},
"2", {\ar"3"},
"3", {\ar"5"},
"8", {\ar"3"},
"4", {\ar"6"},
"9", {\ar"4"},
"5", {\ar"6"},
"6", {\ar"7"},
"7", {\ar"8"},
\end{xy}
\quad\quad
\begin{xy} 0;<0.3pt,0pt>:<0pt,-0.3pt>::
(212,217) *+{1} ="0",
(212,116) *+{2} ="1",
(200,36) *+{3} ="2",
(17,0) *+{4} ="3",
(123,11) *+{5} ="4",
(64,66) *+{6} ="5",
(0,116) *+{7} ="6",
(12,196) *+{8} ="7",
(89,221) *+{9} ="8",
(149,166) *+{10} ="9",
"9", {\ar"0"},
"1", {\ar"2"},
"9", {\ar"1"},
"2", {\ar"4"},
"3", {\ar"5"},
"4", {\ar"5"},
"5", {\ar"6"},
"6", {\ar"7"},
"7", {\ar"8"},
"8", {\ar"9"},
\end{xy}
\quad\quad
\begin{xy} 0;<0.3pt,0pt>:<0pt,-0.3pt>::
(0,230) *+{1} ="0",
(294,255) *+{2.} ="1",
(169,253) *+{3} ="2",
(285,174) *+{4} ="3",
(125,0) *+{5} ="4",
(90,114) *+{6} ="5",
(161,73) *+{7} ="6",
(142,177) *+{8} ="7",
(17,150) *+{9} ="8",
(213,135) *+{10} ="9",
"8", {\ar"0"},
"3", {\ar"1"},
"7", {\ar"2"},
"9", {\ar"3"},
"4", {\ar"6"},
"5", {\ar"6"},
"7", {\ar"5"},
"8", {\ar"5"},
"6", {\ar"9"},
"9", {\ar"7"},
\end{xy}
\]
In fact, in this example, the mutation class is finite and it can be
completely computed
using, for example, \cite{KellerQuiverMutationApplet}:
It consists of $5739$
quivers up to isomorphism. The above quivers are members of the mutation
class containing relatively few arrows. The initial quiver is the unique
member of its mutation class with the largest number of arrows. Here
are some other quivers in the mutation class with a relatively
large number of arrows:
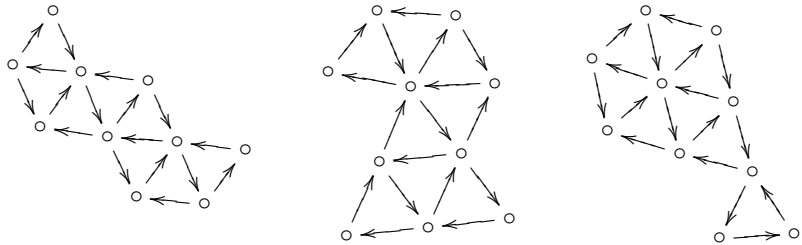
\[
\begin{xy} 0;<0.3pt,0pt>:<0pt,-0.3pt>::
(290,176) *+{\circ} ="0",
(154,235) *+{\circ} ="1",
(34,147) *+{\circ} ="2",
(50,0) *+{\circ} ="3",
(239,244) *+{\circ} ="4",
(0,69) *+{\circ} ="5",
(169,89) *+{\circ} ="6",
(205,165) *+{\circ} ="7",
(85,78) *+{\circ} ="8",
(118,159) *+{\circ} ="9",
"4", {\ar"0"},
"0", {\ar"7"},
"4", {\ar"1"},
"1", {\ar"7"},
"9", {\ar"1"},
"5", {\ar"2"},
"2", {\ar"8"},
"9", {\ar"2"},
"5", {\ar"3"},
"3", {\ar"8"},
"7", {\ar"4"},
"8", {\ar"5"},
"6", {\ar"7"},
"6", {\ar"8"},
"9", {\ar"6"},
"7", {\ar"9"},
"8", {\ar"9"},
\end{xy}
\quad
\quad
\begin{xy} 0;<0.3pt,0pt>:<0pt,-0.3pt>::
(0,78) *+{\circ} ="0",
(226,262) *+{\circ} ="1",
(23,284) *+{\circ} ="2",
(61,0) *+{\circ} ="3",
(208,92) *+{\circ} ="4",
(159,7) *+{\circ} ="5",
(125,273) *+{\circ} ="6",
(64,191) *+{\circ} ="7",
(166,180) *+{\circ} ="8",
(103,96) *+{\circ} ="9",
"0", {\ar"3"},
"9", {\ar"0"},
"1", {\ar"6"},
"8", {\ar"1"},
"6", {\ar"2"},
"2", {\ar"7"},
"5", {\ar"3"},
"3", {\ar"9"},
"5", {\ar"4"},
"8", {\ar"4"},
"4", {\ar"9"},
"9", {\ar"5"},
"7", {\ar"6"},
"6", {\ar"8"},
"8", {\ar"7"},
"7", {\ar"9"},
"9", {\ar"8"},
\end{xy}
\quad
\quad
\begin{xy} 0;<0.3pt,0pt>:<0pt,-0.3pt>::
(159,287) *+{\circ} ="0",
(252,281) *+{\circ} ="1",
(19,152) *+{\circ} ="2",
(67,0) *+{\circ} ="3",
(0,61) *+{\circ} ="4",
(200,203) *+{\circ} ="5",
(109,180) *+{\circ} ="6",
(155,26) *+{\circ} ="7",
(176,115) *+{\circ} ="8",
(87,92) *+{\circ} ="9",
"0", {\ar"1"},
"5", {\ar"0"},
"1", {\ar"5"},
"4", {\ar"2"},
"6", {\ar"2"},
"2", {\ar"9"},
"4", {\ar"3"},
"7", {\ar"3"},
"3", {\ar"9"},
"9", {\ar"4"},
"5", {\ar"6"},
"8", {\ar"5"},
"6", {\ar"8"},
"9", {\ar"6"},
"7", {\ar"8"},
"9", {\ar"7"},
"8", {\ar"9"},
\end{xy}
\]
Only $84$ among the $5739$ quivers in the mutation class contain
double arrows (and none contain arrows of multiplicity $\geq 3$).
Here is a typical example
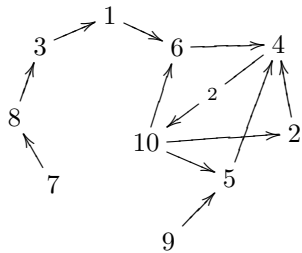
\[
\begin{xy} 0;<0.4pt,0pt>:<0pt,-0.4pt>::
(89,0) *+{1} ="0", (262,111) *+{2} ="1", (24,29) *+{3} ="2",
(247,27) *+{4} ="3", (201,153) *+{5} ="4", (152,30) *+{6} ="5",
(36,159) *+{7} ="6", (0,96) *+{8} ="7", (144,213) *+{9} ="8",
(123,120) *+{10} ="9", "2", {\ar"0"}, "0", {\ar"5"}, "1", {\ar"3"},
"9", {\ar"1"}, "7", {\ar"2"}, "4", {\ar"3"}, "5", {\ar"3"}, "3",
{\ar|*+{\scriptstyle 2}"9"}, "8", {\ar"4"}, "9", {\ar"4"}, "9",
{\ar"5"}, "6", {\ar"7"},
\end{xy}
\]
The classification of the quivers with a finite mutation class
has recently been settled \cite{FeliksonShapiroTumarkin08}:
All possible examples occur in \cite{FominShapiroThurston08}
and \cite{DerksenOwen08}.

The quivers (\ref{quiver1}), (\ref{quiver2}) and (\ref{quiver3}) are part
of a family which appears in the study of the cluster
algebra structure on the coordinate algebra of the
subgroup of upper unitriangular matrices in
$SL(n,\C)$, \cf section~\ref{ss:max-unipotent-sln}.
The quiver (\ref{quiver3}) is associated with the
elliptic root system $E_8^{(1,1)}$ in the notations of Saito \cite{Saito85},
\cf Remark~19.4 in \cite{GeissLeclercSchroeer05}.
The study of coordinate algebras on varieties associated
with reductive algebraic groups
(in particular, double Bruhat cells) has provided a major impetus
for the development of cluster algebras, \cf
\cite{BerensteinFominZelevinsky05}.

\subsection{Definition of cluster algebras}
Let $Q$ be a finite quiver without loops or $2$-cycles with
vertex set $\{1, \ldots, n\}$. Consider the {\em initial seed $(Q,x)$}
consisting of $Q$ and the set $x$ formed by the variables
$x_1, \ldots, x_n$. Following \cite{FominZelevinsky02} we
define
\begin{itemize}
\item the \emph{clusters with respect to $Q$} to be the
sets $u$ appearing in seeds $(R,u)$ obtained from $(Q,x)$ by
iterated mutation,
\item the \emph{cluster variables} for $Q$ to be the elements of all clusters,
\item the \emph{cluster algebra $\ca_Q$} to be the $\Q$-subalgebra of
the field $\Q(x_1, \ldots, x_n)$ generated by all the cluster variables.
\end{itemize}
Thus, the cluster algebra consists of all $\Q$-linear combinations
of monomials in the cluster variables. It is useful to define
yet another combinatorial object associated with this recursive construction:
The \emph{exchange graph} associated with $Q$ is the graph whose
vertices are the seeds modulo simultaneous renumbering of the vertices and
the associated cluster variables and whose edges correspond to mutations.

A remarkable theorem due to Gekhtman-Shapiro-Vainshtein states that
each cluster $u$ occurs in a unique seed $(R,u)$, \cf
\cite{GekhtmanShapiroVainshtein07}.

Notice that the knitting algorithm only produced the cluster variables
whereas this definition yields additional structure: the clusters.

\subsection{The example $A_2$} Here the computation of the exchange
graph is essentially equivalent to performing the knitting algorithm.
If we denote the cluster variables by $x_1$, $x_2$, $x_1'$, $x_2'$ and
$x_1''$ as in section~\ref{ss:knitting-algorithm}, then the exchange
graph is the pentagon
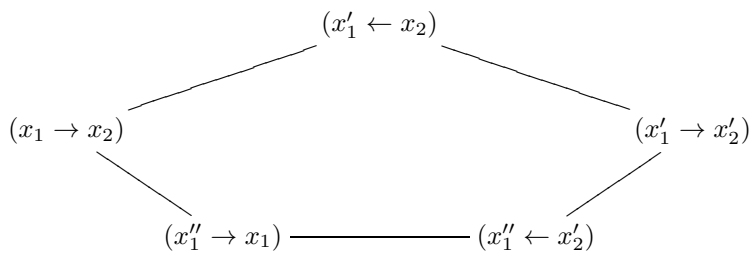
\[
\xymatrix@C=0.3cm{
 &  & (x_1' \la x_2) \ar@{-}[lld] \ar@{-}[rrd] & & \\
 (x_1 \to x_2) \ar@{-}[rd] & & & & (x_1' \to x_2') \ar@{-}[ld] \\
 & (x_1''\to x_1) \ar@{-}[rr] & & (x_1'' \la x_2') &
}
\]
where we have written $x_1 \to x_2$ for the seed $(1\to 2, \{x_1, x_2\})$.
Notice that it is not possible to find a consistent labeling of the
edges by $1$'s and $2$'s. The reason for this is that the vertices
of the exchange graph are not the seeds but the seeds up to renumbering
of vertices and variables. Here the clusters are precisely the pairs of
consecutive variables in the cyclic ordering of $x_1$, \ldots, $x_1''$.

\subsection{The example $A_3$}
Let us consider the quiver
\[
Q: \xymatrix{1 \ar[r] & 2 \ar[r] & 3}
\]
obtained by endowing the Dynkin diagram $A_3$ with a linear orientation.
By applying the recursive construction to the initial seed $(Q,x)$
one finds exactly fourteen seeds (modulo simultaneous renumbering of
vertices and cluster variables). These are the vertices of the exchange
graph, which is isomorphic to the third Stasheff associahedron \cite{Stasheff63}
\cite{ChapotonFominZelevinsky02}:
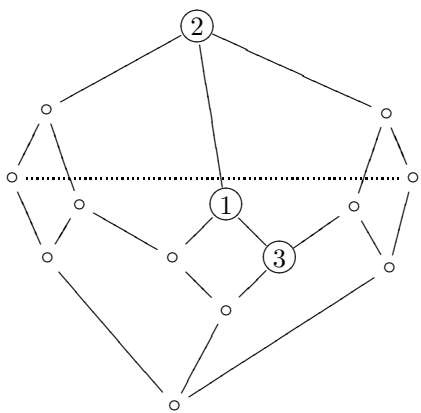
\[
\begin{xy} 0;<0.4pt,0pt>:<0pt,-0.4pt>::
(173,0) *+<8pt>[o][F]{2} ="0",
(0,143) *+{\circ} ="1",
(63,168) *+{\circ} ="2",
(150,218) *+{\circ} ="3",
(250,218) *+<8pt>[o][F]{3} ="4",
(375,143) *+{\circ} ="5",
(350,82) *+{\circ} ="6",
(152,358) *+{\circ} ="7",
(200,168) *+<8pt>[o][F]{1} ="8",
(200,268) *+{\circ} ="9",
(32,79) *+{\circ} ="10",
(33,218) *+{\circ} ="11",
(320,170) *+{\circ} ="12",
(353,228) *+{\circ} ="13",
"0", {\ar@{-}"6"},
"0", {\ar@{-}"8"},
"0", {\ar@{-}"10"},
"1", {\ar@{.}"5"},
"1", {\ar@{-}"10"},
"11", {\ar@{-}"1"},
"2", {\ar@{-}"3"},
"10", {\ar@{-}"2"},
"2", {\ar@{-}"11"},
"3", {\ar@{-}"8"},
"9", {\ar@{-}"3"},
"8", {\ar@{-}"4"},
"4", {\ar@{-}"9"},
"4", {\ar@{-}"12"},
"6", {\ar@{-}"5"},
"5", {\ar@{-}"13"},
"12", {\ar@{-}"6"},
"9", {\ar@{-}"7"},
"11", {\ar@{-}"7"},
"13", {\ar@{-}"7"},
"13", {\ar@{-}"12"},
\end{xy}
\]
The vertex labeled $1$ corresponds to $(Q,x)$, the vertex $2$ to
$\mu_2(Q,x)$, which is given by
\[
\xymatrix{1 \ar@/^1pc/[rr] & 2 \ar[l] & 3 \ar[l]} \ko \{ x_1, \frac{x_1+x_3}{x_2}, x_3\} \ko
\]
and the vertex $3$ to $\mu_1(Q,x)$, which is given by
\[
\xymatrix{1 & 2 \ar[l] \ar[r] & 3} \ko \{\frac{1+x_2}{x_1}, x_2, x_3\}.
\]
As expected (section~\ref{ss:knitting-algorithm}),
we find a total of $3+6=9$ cluster variables, which correspond bijectively
to the faces of the exchange graph.
The clusters ${x_1, x_2, x_3}$ and ${x_1', x_2, x_3}$ also appear
naturally as slices of the repetition, where by a {\em slice}, we
mean a full connected subquiver containing a representative of
each orbit under the horizontal translation (a subquiver is {\em full}
if, with any two vertices, it contains all the arrows between them).
\[
\begin{xy} 0;<0.7pt,0pt>:<0pt,-0.7pt>::
(0,69) *+{x_1} ="0",
(0,0) *+{\circ} ="1",
(35,35) *+{x_2} ="2",
(69,69) *+{x_1'} ="3",
(69,0) *+{x_3} ="4",
(103,35) *+{x_2'} ="5",
(138,69) *+{x_1''} ="6",
(138,0) *+{x_3'} ="7",
(172,35) *+{x_2''} ="8",
(207,69) *+{x_1'''} ="9",
(207,0) *+{x_1} ="10",
(241,35) *+{x_2} ="11",
(275,69) *+{x_3} ="12",
(275,0) *+{x_3'} ="13",
(310,35) *+{x_2'} ="14",
"0", {\ar"2"},
"1", {\ar"2"},
"2", {\ar"3"},
"2", {\ar"4"},
"3", {\ar"5"},
"4", {\ar"5"},
"5", {\ar"6"},
"5", {\ar"7"},
"6", {\ar"8"},
"7", {\ar"8"},
"8", {\ar"9"},
"8", {\ar"10"},
"9", {\ar"11"},
"10", {\ar"11"},
"11", {\ar"12"},
"11", {\ar"13"},
"12", {\ar"14"},
"13", {\ar"14"},
\end{xy}
\]
In fact, as it is
easy to check, each slice yields a cluster. However, some clusters
do not come from slices, for example the cluster ${x_1, x_3, x_1''}$
associated with the seed $\mu_2(Q,x)$.

\subsection{Cluster algebras with finitely many cluster variables}
The phenomena observed in the above examples are explained by the following
key theorem:

\begin{theorem}[Fomin-Zelevinsky \cite{FominZelevinsky03}]
\label{thm:cluster-finite-classification}
Let $Q$ be
  a finite connected quiver without loops or $2$-cycles with vertex
  set $\{1, \ldots, n\}$. Let $\ca_Q$ be the associated cluster
  algebra.
\begin{itemize}
\item[a)] All cluster variables are Laurent polynomials, i.e. their
  denominators are monomials.
\item[b)] The number of cluster variables is finite iff $Q$ is
  mutation equivalent to an orientation of a simply
  laced Dynkin diagram $\Delta$. In this case, $\Delta$ is unique and
  the non initial cluster variables are in bijection with the positive
  roots of $\Delta$; namely, if we denote the simple roots by
  $\alpha_1, \ldots, \alpha_n$, then for each positive root $\sum d_i
  \alpha_i$, there is a unique non initial cluster variable whose
  denominator is $\prod x_i^{d_i}$.
\item[c)] The knitting algorithm yields all cluster variables iff
the quiver $Q$ has two vertices or is an orientation of a simply
laced Dynkin diagram $\Delta$.
\end{itemize}
\end{theorem}

The theorem can be extended to the non simply laced case if
we work with valued quivers as in the example of $G_2$ in
section~\ref{ss:knitting-algorithm}.

It is not hard to check that the knitting algorithm yields exactly
the cluster variables obtained by iterated mutations at sinks and
sources. Remarkably, in the Dynkin case, all cluster variables can be
obtained in this way.

The construction of the
cluster algebra shows that if the quiver $Q$ is mutation-equivalent
to $Q'$, then we have an isomorphism
\[
\ca_{Q'} \iso \ca_Q
\]
preserving clusters and cluster variables. Thus, to prove that
the condition in b) is sufficient, it suffices to show that
$\ca_Q$ is cluster-finite if the underlying graph of $Q$ is
a Dynkin diagram.

No normal form for mutation-equivalence is known in
general and it is unkown how to decide whether two given
quivers are mutation-equivalent.
However, for certain restricted classes, the answer to
this problem is known: Trivially, two quivers with two vertices are
mutation-equivalent iff they are isomorphic. But it is
already a non-trivial problem to decide when a quiver
\[
\xymatrix@=0.3cm{1 \ar[rr]^{r}  & &  2 \ar[ld]^{t}\\
& 3 \ar[lu]^{s} & } ,
\]
where $r$, $s$ and $t$ are non negative integers,
is mutation-equivalent to a quiver without a $3$-cycle:
As shown in \cite{BeinekeBruestleHille06}, this is the
case iff the `Markoff inequality'
\[
r^2+s^2+t^2-rst >4
\]
holds or one among $r$, $s$ and $t$ is $< 2$.

For a general quiver $Q$,
a criterion for $\ca_Q$ to be cluster-finite in terms of
quadratic forms was given in \cite{BarotGeissZelevinsky06}.
In practice, the quickest way to
decide whether a concretely given quiver is cluster-finite and to determine its
cluster-type is to compute its mutation-class using
\cite{KellerQuiverMutationApplet}.
For example, the reader can easily check that for $3\leq n\leq 8$,
the following quiver glued together from $n-2$ triangles
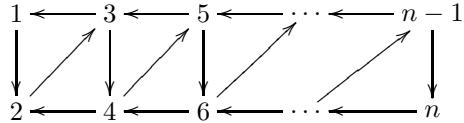
\[
\xymatrix{ 1 \ar[d] & 3 \ar[l] \ar[d] & 5 \ar[l] \ar[d] & \ldots \ar[l] &
n-1 \ar[l]\ar[d] \\
2 \ar[ur] & 4 \ar[l] \ar[ur] & 6 \ar[ur]\ar[l] & \ldots \ar[l]\ar[ur]
& n \ar[l]
}
\]
is cluster-finite of respective cluster-type $A_3$, $D_4$, $D_5$, $E_6$, $E_7$
and $E_8$ and that it is not cluster-finite if $n>8$.

\section{Cluster algebras with coefficients}
\label{s:cluster-alg-with-coeff}

In their combinatorial properties, cluster algebras with coefficients
are very similar to those without coefficients which we have
considered up to now. The great virtue of cluster algebras with
coefficients is that they proliferate in nature as algebras of
coordinates on homogeneous varieties. We will define cluster algebras
with coefficients and illustrate their ubiquity on several examples.

\subsection{Definition} \label{ss:def-cluster-alg-coeff}
Let $1\leq n \leq m$ be integers. An {\em ice quiver of type $(n,m)$} is
a quiver $\tilde{Q}$ with vertex set
\[
\{1, \ldots, m\} = \{1,\ldots, n\} \cup \{n+1, \ldots, m\}
\]
such that there are no arrows between any vertices $i,j$ which are strictly
greater than $n$. The {\em principal part} of $\tilde{Q}$ is the full
subquiver $Q$ of $\tilde{Q}$ whose vertex set is
$\{1, \ldots, n\}$ (a subquiver is {\em full} if, with any two vertices,
it contains all the arrows between them).
The vertices $n+1$, \ldots, $m$ are often called
{\em frozen vertices}. The cluster algebra
\[
\ca_{\tilde{Q}} \subset \Q(x_1, \ldots, x_m)
\]
is defined as before but
\begin{itemize}
\item[-] only mutations with respect to vertices in the principal
part are allowed and no arrows are drawn between the vertices
$>n$,
\item[-] in a cluster
\[
u=\{u_1, \ldots, u_n, x_{n+1}, \ldots, x_m\}
\]
only $u_1$, \ldots, $u_n$ are called cluster variables; the elements
$x_{n+1}$, \ldots, $x_m$ are called {\em coefficients}; to make things
clear, the set $u$ is often called an {\em extended cluster};
\item[-] the {\em cluster type of $\tilde{Q}$} is that of $Q$ if it
is defined.
\end{itemize}
Often, one also considers localizations of $\ca_{\tilde{Q}}$ obtained
by inverting certain coefficients.
Notice that the datum of $\tilde{Q}$ corresponds to that of the
integer $m\times n$-matrix $\tilde{B}$ whose top $n\times n$-submatrix
$B$ is antisymmetric and whose entry $b_{ij}$ equals the number
of arrows $i\to j$ or the opposite of the number of arrows
$j\to i$. The matrix $B$ is called the {\em principal part} of
$\tilde{B}$. One can also consider valued ice quivers, which
will correspond to $m\times n$-matrices whose principal part
is antisymmetrizable.

\subsection{Example: $SL(2,\C)$}
Let us consider the algebra of regular functions
on the algebraic group $SL(2,\C)$, \ie the algebra
\[
\C[a,b,c,d]/(ad-bc-1).
\]
We claim that this algebra {\em has a cluster algebra structure},
namely that it is isomorphic to the complexification of the
cluster algebra with coefficients associated with the following
ice quiver
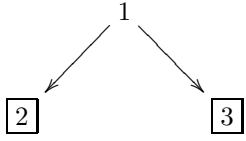
\[
\xymatrix{ & 1 \ar[dl] \ar[dr] & \\
\framebox{2} & & \framebox{3}
}
\]
where we have framed the frozen vertices. Indeed, here the
principal part $Q$ only consists of the vertex $1$ and we can
only perform one mutation, whose associated exchange relation reads
\[
x_1 x_1' = 1 + x_2 x_3 \mbox{ or } x_1 x_1' - x_2 x_3 = 1.
\]
We obtain an isomorphism as required by sending $x_1$ to $a$,
$x_1'$ to $d$, $x_2$ to $b$ and $x_3$ to $c$. We describe this
situation by saying that the quiver
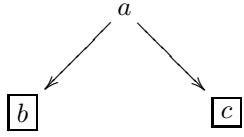
\[
\xymatrix{ & a \ar[dl] \ar[dr] & \\
\framebox{$b$} & & \framebox{$c$}
}
\]
whose vertices are labeled by the images of the corresponding
variables, is an initial seed for a cluster structure
on the algebra $A$. Notice that this cluster structure is not unique.

\subsection{Example: Planes in affine space}
As a second example, let us consider the algebra $A$ of polynomial functions
on the cone over the Grassmannian of planes in $\C^{n+3}$. This
algebra is the $\C$-algebra generated by the Pl\"ucker coordinates
$x_{ij}$, $1\leq i<j\leq n+3$, subject to the Pl\"ucker relations:
for each quadruple of integers $i<j<k<l$, we have
\[
x_{ik} x_{jl} = x_{ij} x_{kl} + x_{jk} x_{il}.
\]
Each plane $P$ in $\C^{n+3}$ gives rise to a straight line in this cone,
namely the one generated by the $2\times 2$-minors $x_{ij}$ of any
$(n+3)\times 2$-matrix whose columns generate $P$. Notice that
the monomials in the Pl\"ucker relation are naturally associated
with the sides and the diagonals of the square
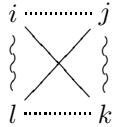
\[
\xymatrix{ i \ar@{.}[r] \ar@{~}[d] \ar@{-}[dr] & j \ar@{-}[dl] \ar@{~}[d] \\
l \ar@{.}[r] & k}
\]
The relation expresses the product of the variables associated with
the diagonals as the sum of the monomials associated with the two
pairs of opposite sides.

Now the idea is that the Pl\"ucker relations are
exactly the exchange relations for a suitable structure of cluster
algebra with coefficients on the coordinate ring.  To formulate this
more precisely, let us consider a regular $(n+3)$-gon in the plane
with vertices numbered $1$, \ldots, $n+2$, and consider the variable
$x_{ij}$ as associated with the segment $[ij]$ joining the vertices
$i$ and $j$.

\begin{proposition}[\protect{\cite[Example 12.6]{FominZelevinsky03}}]
The algebra $A$ has a structure of cluster algebra
with coefficients such that
\begin{itemize}
\item[-] the coefficients are the variables $x_{ij}$ associated with
the sides of the $(n+3)$-gon;
\item[-] the cluster variables are the variables $x_{ij}$ associated
with the diagonals of the $(n+3)$-gon;
\item[-] the clusters are the $n$-tuples of variables whose associated
  diagonals form a triangulation of the $(n+3)$-gon.
\end{itemize}
Moreover, the exchange relations are exactly the Pl\"ucker relations and
the cluster type is $A_n$.
\end{proposition}

Thus, a triangulation of the $(n+3)$-gon determines an initial seed for
the cluster algebra and hence an ice quiver $\tilde{Q}$ whose frozen
vertices correspond to the sides of the $(n+3)$-gon and whose non frozen
variables to the diagonals in the triangulation. The arrows of the
quiver are determined by the exchange relations which appear when we
wish to replace one diagonal $[ik]$ of the triangulation by its {\em flip},
i.e. the unique diagonal $[jl]$ different from $[ik]$ which does not
cross any other diagonal of the triangulation. It is not hard to see
that this means that the underlying graph of $\tilde{Q}$ is the
graph dual to the triangulation and that the orientation of the edges
of this graph is induced by the choice of an orientation of the plane.
Here is an example of a triangulation and the associated ice quiver:
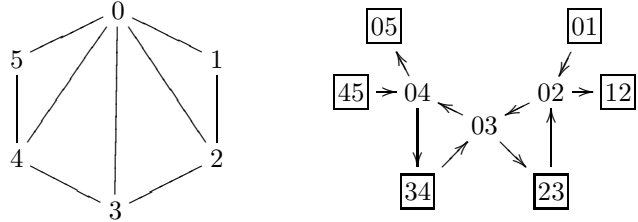
\[
\begin{xy} 0;<0.5pt,0pt>:<0pt,-0.5pt>::
(76,0) *+{0} ="0",
(150,37) *+{1} ="1",
(150,112) *+{2} ="2",
(74,151) *+{3} ="3",
(0,112) *+{4} ="4",
(0,37) *+{5} ="5",
(275,12) *+{\framebox[3ex]{05}} ="6",
(300,62) *+{04} ="7",
(350,87) *+{03} ="8",
(400,62) *+{02} ="9",
(425,12) *+{\framebox[3ex]{01}} ="10",
(250,62) *+{\framebox[3ex]{45}} ="11",
(300,137) *+{\framebox[3ex]{34}} ="12",
(400,137) *+{\framebox[3ex]{23}} ="13",
(450,62) *+{\framebox[3ex]{12}} ="14",
"0", {\ar@{-}"1"},
"0", {\ar@{-}"2"},
"0", {\ar@{-}"3"},
"0", {\ar@{-}"4"},
"5", {\ar@{-}"0"},
"1", {\ar@{-}"2"},
"2", {\ar@{-}"3"},
"3", {\ar@{-}"4"},
"4", {\ar@{-}"5"},
"7", {\ar"6"},
"8", {\ar"7"},
"11", {\ar"7"},
"7", {\ar"12"},
"9", {\ar"8"},
"12", {\ar"8"},
"8", {\ar"13"},
"10", {\ar"9"},
"13", {\ar"9"},
"9", {\ar"14"},
\end{xy}
\]

\subsection{Example: The big cell of the Grassmannian}
We consider the cone over the big
cell in the Grassmannian of $k$-dimensional subspaces of the space
of rows $\C^n$, where $1\leq k\leq n$ are fixed integers such that
$l=n-k$ is greater or equal to $2$. In more detail,
let $G$ be the group $SL(n,\C)$ and $P$ the subgroup of $G$ formed
by the block lower triangular matrices with diagonal blocks
of sizes $k\times k$ and $l\times l$. The quotient
$P\setminus G$ identifies with our Grassmannian. The big cell is the
image under $\pi: G \to P\setminus G$ of the space of block upper triangular
matrices whose diagonal is the identity matrix and whose upper
right block is an arbitrary $k\times l$-matrix $Y$. The
projection $\pi$ induces an isomorphism between the space
$M_{k\times l}(\C)$ of these matrices and the big cell.
In particular, the algebra $A$ of regular functions on the
big cell is the algebra of polynomials in the coefficients
$y_{ij}$, $1\leq i\leq k$, $1\leq j\leq l$, of $Y$.
Now for $1\leq i\leq k$ and $1 \leq j\leq l$, let $F_{ij}$
be the largest square submatrix of $Y$ whose lower left
corner is $(i,j)$ and let $s(i,j)$ be its size. Put
\[
f_{ij} = (-1)^{(k-i)(s(i,j)-1)} \det(F_{ij}).
\]

\begin{theorem}[\cite{GekhtmanShapiroVainshtein03}]
The algebra $A$ has the structure of a cluster
algebra with coefficients whose initial seed is given by
\end{theorem}

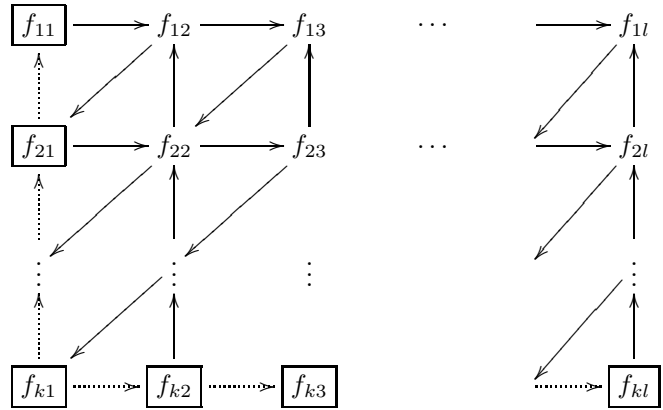
\[
\xymatrix{
\framebox{$f_{11}$} \ar[r]
        & f_{12} \ar[r] \ar[dl]
                 & f_{13} \ar[dl]
                          & \ldots & \ar[r] & f_{1l} \ar[dl] \\
\framebox{$f_{21}$} \ar[r] \ar@{.>}[u]
        & f_{22} \ar[r] \ar[u]\ar[dl]
                 &  f_{23} \ar[u] \ar[dl]
                          & \ldots & \ar[r] & f_{2l} \ar[u] \ar[dl]   \\
\vdots \ar@{.>}[u]
        & \vdots \ar[dl] \ar[u]
                 & \vdots &        &        & \vdots \ar[dl]\ar[u]    \\
\framebox{$f_{k1}$} \ar@{.>}[r] \ar@{.>}[u]
        & \framebox{$f_{k2}$} \ar[u] \ar@{.>}[r]
                 & \framebox{$f_{k3}$}
                          &        &  \ar@{.>}[r]
                                            & \framebox{$f_{kl}$} \ar[u]
}
\]

The following table indicates when these algebras are
cluster-finite and what their cluster-type is:
\[
\begin{array}{c|cccccc}
k \setminus n & 2 & 3 & 4 & 5 & 6 & 7 \\ \hline
2 & A_1 & A_2 & A_3 & A_4 & A_5 & A_6 \\
3 & A_2 & D_4 & E_6 & E_8 &       &       \\
4 & A_3 & E_6 &       &       &       &       \\
5 & A_4 & E_8   &       &       &       &
\end{array}
\]
The homogeneous coordinate ring of the Grassmannian
$Gr(k,n)$ itself also has a cluster algebra structure
\cite{Scott06} and so have partial flag varieties,
double Bruhat cells, Schubert varieties \ldots,
\cf \cite{GeissLeclercSchroeer08b} \cite{BerensteinFominZelevinsky05}.

\subsection{Compatible Poisson structures}
Recall that the group $SL(n,\C)$ has a canonical Poisson
structure given by the {\em Sklyanin bracket}, which is
defined by
\[
\omega(x_{ij}, x_{\alpha \beta}) = (\sign(\alpha-i) - \sign(\beta-j))
x_{i\beta} x_{\alpha j}
\]
where the $x_{ij}$ are the coordinate functions
on $SL(n,\C)$. This bracket makes $G=SL(n,\C)$ into a
Poisson-Lie group and $P\setminus G$ into a Poisson $G$-variety
for each subgroup $P$ of $G$ containing the subgroup $B$ of lower
triangular matrices. In particular, the big cell of the
Grassmannian considered above inherits a Poisson bracket.

\begin{theorem}[\cite{GekhtmanShapiroVainshtein03}]
This bracket is compatible with the cluster algebra structure
in the sense that each extended cluster is a {\em log-canonical} coordinate
system, i.e. we have
\[
\omega(u_i, u_j) = \omega^{(u)}_{ij} u_i u_j
\]
for certain (integer) constants $\omega^{(u)}_{ij}$ depending on the
extended cluster $u$. Moreover, the coefficients are central for
$\omega$.
\end{theorem}

This theorem admits the following generalization: Let $\tilde{Q}$ be an
ice quiver. Define the {\em cluster variety} $\cx(\tilde{Q})$ to
be obtained by glueing the complex tori indexed by the clusters $u$
\[
T^{(u)}=(\C^*)^n = \Spec(\C[u_1, u_1^{-1}, \ldots, u_m, u_m^{-1}])
\]
using the exchange relations as glueing maps, where $m$ is
the number of vertices of $\tilde{Q}$.

\begin{theorem}[\cite{GekhtmanShapiroVainshtein03}] Suppose that the
  principal part $Q$ of $\tilde{Q}$ is connected and that the matrix
  $\tilde{B}$ associated with $\tilde{Q}$ is of maximal rank. Then the
  vector space of Poisson structures on $\cx(\tilde{Q})$ compatible
  with the cluster algebra structure is of dimension
\[
1 + {{m-n}\choose{2}}.
\]
\end{theorem}

Notice that in general, the cluster variety $\cx(\tilde{Q})$ is
an open subset of the spectrum of the (complexified) cluster
algebra. For example, for the cluster algebra associated
with $SL(2,\C)$ which we have considered above, the cluster
variety is the union of the elements
\[
\left( \begin{array}{cc} a & b \\ c & d \end{array} \right)
\]
of $SL(2,\C)$ such that we have $abc\neq 0$ or $bcd\neq 0$.
The cluster variety is always regular, but the spectrum
of the cluster algebra may be singular. For example,
the spectrum of the cluster algebra associated with
the ice quiver
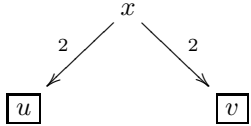
\[
\xymatrix{ & x \ar[dl]_2 \ar[dr]^2 & \\
\framebox{$u$} & & \framebox{$v$}
}
\]
is the hypersurface in $\C^4$ defined by the equation $x x'=u^2 +v^2$,
which is singular at the origin. The corresponding cluster variety
is obtained by removing the points with $x=x'=u^2+v^2=0$ and is
regular.

In the above theorem, the assumption that $\tilde{B}$ be of
full rank is essential. Otherwise, there may not exist any
Poisson bracket compatible with the cluster algebra structure.
However, as shown in \cite{GekhtmanShapiroVainshtein05},
for any cluster algebra with coefficients, there are
`dual Poisson structures', namely certain $2$-forms,
which are compatible with the cluster algebra structure.

\subsection{Example: The maximal unipotent subgroup of $SL(n+1,\C)$}
\label{ss:max-unipotent-sln}
Let $n$ be a non negative integer and $N$ the subgroup of $SL(n+1,\C)$
formed by the upper triangular matrices with all diagonal coefficients
equal to $1$. For $1\leq i,j \leq n+1$ and $g\in N$, let $F_{ij}(g)$ be the maximal
square submatrix of $g$ whose lower left corner is $(i,j)$. Let
$f_{ij}(g)$ be the determinant of $F_{ij}(g)$. We consider the
functions $f_{ij}$ for $1 \leq i \leq n$ and $i+j\leq n+1$.

\begin{theorem}[\cite{BerensteinFominZelevinsky05}] The coordinate algebra $\C[N]$
has an upper cluster algebra structure whose initial seed is given by
\end{theorem}
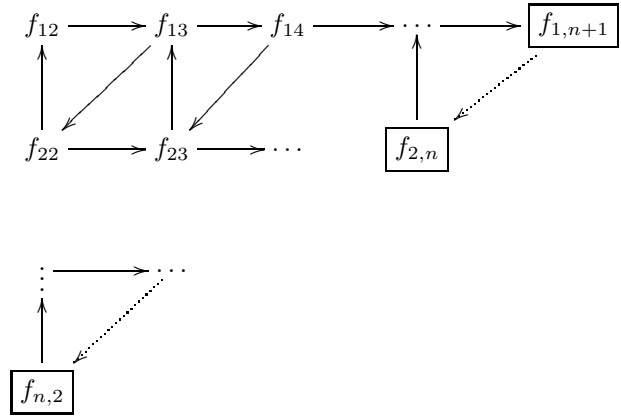
\[
\xymatrix{
f_{12} \ar[r] & f_{13} \ar[r] \ar[dl] & f_{14} \ar[dl] \ar[r] & \ldots \ar[r] & \framebox{$f_{1,n+1}$} \ar@{.>}[dl] \\
f_{22} \ar[u] \ar[r]
             & f_{23} \ar[u] \ar[r]  & \ldots         & \framebox{$f_{2,n}$} \ar[u]     &   \\
\vdots \ar[r]       &  \ldots \ar@{.>}[dl]                     &          &               &   \\
\framebox{$f_{n,2}$} \ar[u] &                    &                &               &
}
\]
We refer to \cite{BerensteinFominZelevinsky05} for the
notion of `upper' cluster algebra structure.
It is not hard to check that this structure is of cluster type $A_3$ for
$n=3$, $D_6$ for $n=4$ and cluster-infinite for $n\geq 5$. For
$n=5$, this cluster algebra is related to the
elliptic root system $E_8^{(1,1)}$ in the notations of Saito \cite{Saito85},
\cf \cite{GeissLeclercSchroeer05}.

A theorem of Fekete \cite{Fekete1912} generalized in \cite{BerensteinFominZelevinsky96}
claims that a square matrix of order $n+1$ is {\em totally positive} (\ie all its minors are $>0$)
if and only if the following $(n+1)^2$ minors of g are positive:
all minors occupying several initial rows and several consecutive columns,
and all minors occupying several initial columns and several consecutive
rows. It follows that an element $g$ of $N$ is totally positive
if $f_{ij}(g)>0$ for the $f_{ij}$
belonging to the initial seed above. The same holds for the
$u_1, \ldots, u_m$ in place of these $f_{ij}$ for any
cluster $u$ of this cluster algebra because each exchange
relation expresses the new variable {\em subtraction-free}
in the old variables.

Geiss-Leclerc-Schr\"oer have shown \cite{GeissLeclercSchroeer06}
that each monomial in the variables of an arbitrary cluster
belongs to Lusztig's dual semicanonical basis of $\C[N]$
\cite{Lusztig00}. They also show that the dual semicanonical
basis of $\C[N]$ is different from the dual canonical basis
of Lusztig and Kashiwara except in types $A_2$, $A_3$ and $A_4$
\cite{GeissLeclercSchroeer05}.

\section{Categorification via cluster categories: the finite case}
\label{s:categorification-finite-case}

\subsection{Quiver representations and Gabriel's theorem}
We refer to the books \cite{Ringel84} \cite{GabrielRoiter92}
\cite{AuslanderReitenSmaloe95} and \cite{AssemSimsonSkowronski06} for
a wealth of information on the representation theory of quivers and
finite-dimensional algebras. Here, we will only need very basic notions.

Let $Q$ be a finite quiver without oriented cycles. For example, $Q$ can
be an orientation of a simply laced Dynkin diagram or
the quiver
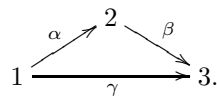
\[
\xymatrix@R=10pt{ & 2 \ar[rd]^\beta & \\
1 \ar[rr]_{\gamma} \ar[ru]^{\alpha} & & 3.
}
\]
Let $k$ be an algebraically closed field. A {\em
representation of $Q$} is a diagram of finite-dimensional vector
spaces of the shape given by $Q$. More formally, a representation of
$Q$ is the datum $V$ of
\begin{itemize}
\item a finite-dimensional vector space $V_i$ for each vertex $i$ or $Q$,
\item a linear map $V_\alpha: V_i \to V_j$ for each arrow $\alpha: i\to j$
of $Q$.
\end{itemize}
Thus, in the above example, a representation of $Q$ is a (not
necessarily commutative) diagram
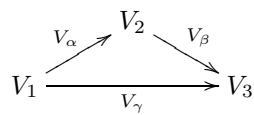
\[
\xymatrix@R=10pt{ & V_2 \ar[rd]^{V_\beta} & \\
V_1 \ar[rr]_{V_{\gamma}} \ar[ru]^{V_\alpha} & & V_3
}
\]
formed by three finite-dimensional vector spaces and three linear maps.
A \emph{morphism of representations} is a morphism of diagrams. More formally,
a morphism of representations $f:V \to W$ is the datum of a linear map
$f_i: V_i \to W_i$ for each vertex $i$ of $Q$ such that the square
\[
\xymatrix{
V_i \ar[d]_{f_i} \ar[r]^{V_\alpha} & V_j \ar[d]^{f_j} \\
W_i \ar[r]_{W_\alpha} & W_j
}
\]
commutes for all arrows $\alpha:i \to j$ of $Q$. The {\em composition}
of morphisms is defined in the natural way.
We thus obtain the \emph{category of representations} $\rep(Q)$.
A morphism $f:V \to W$ of this category is an isomorphism iff its
components $f_i$ are invertible for all vertices $i$ of $Q_0$.

For example, let $Q$ be the quiver
\[
\xymatrix{1 \ar[r] & 2} \ko
\]
and
\[
V: \xymatrix{V_1 \ar[r]^{V_\alpha} & V_2}
\]
a representation of $Q$. By choosing basis in the spaces $V_1$ and $V_2$
we find an isomorphism of representations
\[
\xymatrix{V_1 \ar[r]^{V_\alpha} & V_2 \\
k^n \ar[u] \ar[r]_{A}  & k^p \ko \ar[u]}
\]
where, by abuse of notation, we denote by $A$ the multiplication by
a $p\times n$-matrix $A$. We know that we have
\[
P A Q = \left[ \begin{array}{cc} I_r & 0 \\ 0 & 0 \end{array} \right]
\]
for invertible matrices $P$ and $Q$, where $r$ is the rank of $A$.
Let us denote the right hand side by $I_r \oplus 0$. Then we have an
isomorphism of representations
\[
\xymatrix{
k^n \ar[r]_{A} & k^p  \\
k^n \ar[u]^Q \ar[r]_{I_r\oplus 0} & k^p \ar[u]_{P^{-1}}
}
\]
We thus obtain a normal form for the representations of this quiver.

Now the category $\rep_k(Q)$ is in fact an {\em abelian category}:
Its direct sums, kernels and cokernels are computed componentwise.
Thus, if $V$ and $W$ are two representations, then the
{\em direct sum} $V\oplus W$
is the representation given by
\[
(V\oplus W)_i = V_i \oplus W_i \mbox{ and }
(V\oplus W)_\alpha= V_\alpha \oplus W_\alpha \ko
\]
for all vertices $i$ and all arrows $\alpha$ of $Q$. For
example, the above representation in normal form is isomorphic to
the direct sum
\[
(\xymatrix{k \ar[r]^{1}  & k})^r \oplus (\xymatrix{k \ar[r] & 0})^{n-r}
\oplus (\xymatrix{0 \ar[r] & k})^{p-r}.
\]
The {\em kernel} of a morphism of representations
$f: V \to W$ is given by
\[
\ker(f)_i = \ker(f_i: V_i \to W_i)
\]
endowed with the maps induced by the $V_\alpha$ and
similarly for the cokernel. A {\em subrepresentation} $V'$
of a representation $V$ is given by a family of subspaces
$V'_i \subset V_i$, $i\in Q_0$, such that the image of
$V'_i$ under $V_\alpha$ is contained in $V'_j$ for each
arrow $\alpha: i \to j$ of $Q$. A sequence
\[
\xymatrix{ 0 \ar[r] & U \ar[r] & V \ar[r] & W \ar[r] & 0}
\]
of representations is a {\em short exact sequence} if the sequence
\[
\xymatrix{ 0 \ar[r] & U_i \ar[r] & V_i \ar[r] & W_i \ar[r] & 0}
\]
is exact for each vertex $i$ of $Q$.

A representation $V$ is {\em simple} if it is non zero and
if for each subrepresentation $V'$ of $V$ we have $V'=0$ or
$V/V'=0$. Equivalently, a representation is simple if it
has exactly two subrepresentations.
A representation $V$ is {\em indecomposable} if it
is non zero and in each decomposition $V=V'\oplus V''$, we
have $V'=0$ or $V''=0$. Equivalently, a representation
is indecomposable if it has exactly two direct factors.

In the above example, the representations
\[
\xymatrix{ k \ar[r] & 0} \mbox{ and } \xymatrix{ 0 \ar[r] & k}
\]
are simple. The representation
\[
V=(\xymatrix{k \ar[r]^1 & k})
\]
is not simple: It has the non trivial subrepresentation
$\xymatrix{ 0 \ar[r] & k}$. However, it is indecomposable.
Indeed, each endomorphism $f: V \to V$ is given by two
equal components $f_1=f_2$ so that the endomorphism algebra
of $V$ is one-dimensional. If $V$ was a direct sum
$V'\oplus V''$ for two non-zero subspaces, the endomorphism
algebra of $V$ would contain the product of the endomorphism
algebras of $V'$ and $V''$ and thus would have to be at least
of dimension $2$. Since $V$ is indecomposable, the exact
sequence
\[
0 \to (\xymatrix{0 \ar[r] & k}) \to (\xymatrix{k\ar[r]^1 & k}) \to (\xymatrix{k \ar[r] & 0}) \to 0
\]
is not a split exact sequence.

If $Q$ is an arbitrary quiver, for each vertex $i$, we define
the representation $S_i$ by
\[
(S_i)_j = \left\{ \begin{array}{ll} k & i=j \\ 0 & \mbox{else.} \end{array} \right.
\]
Then clearly the representations $S_i$ are simple and
pairwise non isomorphic. As an exercise, the reader may show that
if $Q$ does not have oriented cycles, then each representation
admits a finite filtration whose subquotients are among the
$S_i$. Thus, in this case, each simple representation is
isomorphic to one of the representations $S_i$.

Recall that a (possibly non commutative) ring is {\em local}
if its non invertible elements form an ideal.

\begin{decomposition-theorem}[Azumaya-Fitting-Krull-Remak-Schmidt]
\label{thm:decomposition}
\quad
\begin{itemize}
\item[a)] A representation is indecomposable iff its endomorphism
algebra is local.
\item[b)] Each representation decomposes into a finite sum of
indecomposable representations, unique up to isomorphism and
permutation.
\end{itemize}
\end{decomposition-theorem}

As we have seen above, for quivers without oriented cycles,
the classification of the simple representations is trivial.
On the other hand, the problem of classifying the indecomposable
representations is non trivial. Let us examine this problem
in a few examples:
For the quiver $1 \to 2$, we have checked the existence in
part b) directly. The uniqueness in b) then implies that each
indecomposable representation is isomorphic to exactly one
of the representations $S_1$, $S_2$ and
\[
\xymatrix{k \ar[r]^1 & k}.
\]
Similarly, using elementary linear algebra it is not hard to
check that each indecomposable representation of the quiver
\[
\vec{A}_n : \xymatrix{1 \ar[r] & 2 \ar[r] & \ldots \ar[r] & n}
\]
is isomorphic to a representation $I[p,q]$, $1 \leq p < q \leq n$,
which takes the vertices $i$ in the interval $[p,q]$ to $k$,
the arrows linking them to the identity and all other vertices
to zero. In particular, the number of isomorphism classes
of indecomposable representations of $\vec{A}_n$ is $n(n+1)/2$.

The representations of the quiver
\[
\xymatrix{ 1 \ar@(ur,dr)[]^\alpha }
\]
are the pairs $(V_1, V_\alpha)$ consisting of a finite-dimensional vector
space and an endomorphism and the morphisms of representations
are the `intertwining operators'. It follows from the existence and
uniqueness of the Jordan normal form that a system of representatives
of the isomorphism classes of indecomposable representations is
formed by the representations $(k^n, J_{n,\lambda}$), where $n\geq 1$
is an integer, $\lambda$ a scalar and $J_{n,\lambda}$ the
Jordan block of size $n$ with eigenvalue $\lambda$.

The {\em Kronecker quiver}
\[
\xymatrix{1 \ar@<1ex>[r] \ar@<-1ex>[r] & 2}
\]
admits the following infinite family of pairwise non isomorphic
representations:
\[
\xymatrix{k \ar@<1ex>[r]^\lambda \ar@<-1ex>[r]_\mu & k} \ko
\]
where $(\lambda:\mu)$ runs through the projective line.

\begin{question} For which quivers are there only finitely many
isomorphism classes of indecomposable representations?
\end{question}

To answer this question, we define the {\em dimension vector}
of a representation $V$ to be the sequence $\dimv V$ of
the dimensions $\dim V_i$, $i\in Q_0$. For example, the
dimension vectors of the indecomposable representations
of $\vec{A}_2$ are the pairs
\[
\dimv S_1 = [1 0]\ko \dimv S_2 = [0 1] \ko \dimv (k\to k) = [1 1].
\]
We define the {\em Tits form}
\[
q_Q: \Z^{Q_0} \to \Z
\]
by
\[
q_Q(v) = \sum_{i\in Q_0} v_i^2 - \sum_{\alpha\in Q_1} v_{s(\alpha)} v_{t(\alpha)}.
\]
Notice that the Tits form does not depend on the orientation
of the arrows of $Q$ but only on its underlying graph. We say
that the quiver $Q$ is {\em representation-finite} if, up to
isomorphism, it has only finitely many indecomposable representations.
We say that a vector $v\in \Z^{Q_0}$ is a {\em root} of $q_Q$
if $q_Q(v)=1$ and that it is {\em positive} if its components are $\geq 0$.

\begin{theorem}[Gabriel \protect{\cite{Gabriel72}}] Let $Q$ be a connected quiver
and assume that $k$ is algebraically closed. The following are equivalent:
\begin{itemize}
\item[(i)] $Q$ is representation-finite;
\item[(ii)] $q_Q$ is positive definite;
\item[(iii)] The underlying graph of $Q$ is a simply laced Dynkin diagram $\Delta$.
\end{itemize}
Moreover, in this case, the map taking a representation to its dimension
vector yields a bijection from the set of isomorphism classes of indecomposable
representations to the set of positive roots of the Tits form $q_Q$.
\end{theorem}

It is not hard to check that if the conditions hold, the positive
roots of $q_Q$ are in turn in bijection with the positive roots
of the root system $\Phi$ associated with $\Delta$, via the map taking
a positive root $v$ of $q_Q$ to the element
\[
\sum_{i\in Q_0} v_i \alpha_i
\]
of the root lattice of $\Phi$.

Let us consider the example of the
quiver $Q=\vec{A}_2$. In this case, the Tits form is given by
\[
q_Q(v)= v_1^2+v_2^2-v_1 v_2.
\]
It is positive definite and its positive roots are indeed
precisely the dimension vectors
\[
[0 1] \ko [1 0] \ko [1 1]
\]
of the indecomposable representations.

Gabriel's theorem has been generalized to non algebraically closed
ground fields by Dlab and Ringel \cite{DlabRingel76}. Let us illustrate
the main idea on one simple example: Consider the category of
diagrams
\[
V: \xymatrix{V_1 \ar[r]^f & V_2}
\]
where $V_1$ is a finite-dimensional real vector space, $V_2$ a finite-dimensional
complex vector space and $f$ an $\R$-linear map. Morphisms are given in
the natural way. Then we have the following complete list of representatives
of the isomorphism classes of indecomposables:
\[
\R \to 0 \ko \R^2 \to \C \ko \R \to \C \ko 0 \to \C.
\]
The corresponding dimension vectors are
\[
[1 0] \ko [2 1] \ko [1 1] \ko [0 1].
\]
They correspond bijectively to the positive roots of the root
system $B_2$.

\subsection{Tame and wild quivers}
The quivers with infinitely many isomorphism classes of indecomposables
can be further subdivided into two important classes: A quiver is
{\em tame} if it has infinitely many isomorphism classes of
indecomposables but these occur in `families of at most one
parameter' (we refer to \cite{Ringel84} \cite{AssemSimsonSkowronski06}
for the precise definition). The Kronecker quiver is a
typical example. A quiver is {\em wild} if there are `families
of indecomposables of $\geq 2$ parameters'. One can show that
in this case, there are families of an arbitrary number of
parameters and that the classification of the indecomposables
over any fixed wild algebra would entail the classification of
the the indecomposables over all finite-dimensional algebras.
The following
three quivers are representation-finite, tame and wild
respectively:
\[
\xymatrix{  & 1 \ar[dl] \ar[d] \ar[dr] & \\ 2 & 3 & 4 } \quad\quad
\xymatrix@C=0.2cm{ & & 1 \ar[dll] \ar[dl] \ar[dr] \ar[drr] & & \\ 2 & 3 & & 4 & 5 } \quad\quad
\xymatrix@C=0.2cm{ & & 1 \ar[dll] \ar[dl] \ar[d] \ar[dr] \ar[drr] & & \\ 2 & 3 & 4 & 5 & 6. }
\]

\begin{theorem}[\protect{Donovan-Freislich \cite{DonovanFreislich73}, Nazarova \cite{Nazarova73}}]
Let $Q$ be a connected quiver and assume that $k$ is algebraically closed. Then
$Q$ is tame iff the underlying graph of $Q$ is a simply laced extended Dynkin
diagram.
\end{theorem}
Let us recall the list of simply laced extended Dynkin quivers. In each
case, the number of vertices of the diagram $\tilde{D}_n$ equals $n+1$.
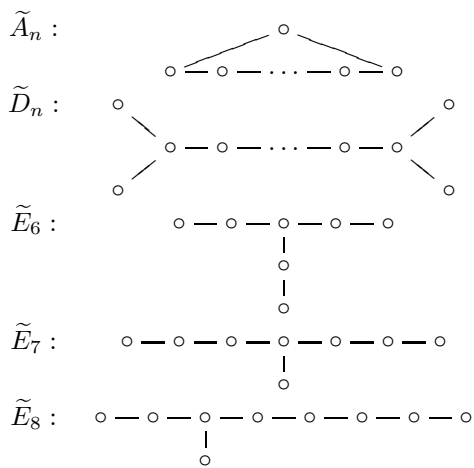
\[
\begin{array}{cc} \tilde{A}_n : & \xymatrix@R=0.2cm@C=0.3cm{ & & \circ \ar@{-}[drr] & & \\
\circ \ar@{-}[urr] \ar@{-}[r] & \circ \ar@{-}[r] & \ldots \ar@{-}[r] & \circ \ar@{-}[r] & \circ
} \\
\tilde{D}_n : & \xymatrix@R=0.2cm@C=0.3cm{ \circ \ar@{-}[dr] & & & & & & \circ \\
 & \circ \ar@{-}[r] & \circ \ar@{-}[r] & \ldots \ar@{-}[r] &
 \circ \ar@{-}[r] & \circ \ar@{-}[ru] \ar@{-}[rd] \\
 \circ \ar@{-}[ru] & & & & & & \circ} \\
\tilde{E}_6 : & \xymatrix@R=0.2cm@C=0.3cm{\circ \ar@{-}[r] & \circ \ar@{-}[r] & \circ \ar@{-}[r] \ar@{-}[d] & \circ \ar@{-}[r] & \circ \\
 & & \circ \ar@{-}[d] & & \\
 & & \circ & & } \\
\tilde{E}_7 : & \xymatrix@R=0.2cm@C=0.3cm{
\circ \ar@{-}[r] & \circ \ar@{-}[r] & \circ \ar@{-}[r] & \circ \ar@{-}[r] \ar@{-}[d] &
\circ \ar@{-}[r] & \circ \ar@{-}[r] & \circ \\
 & & & \circ & & &} \\
\tilde{E}_8 : & \xymatrix@R=0.2cm@C=0.3cm{
\circ \ar@{-}[r] & \circ \ar@{-}[r] & \circ \ar@{-}[r] \ar@{-}[d] &
\circ \ar@{-}[r] & \circ \ar@{-}[r] & \circ \ar@{-}[r] & \circ \ar@{-}[r] & \circ \\
 & & \circ }
\end{array}
\]

The following theorem is a first illustration of the close connection
between cluster algebras and the representation theory of quivers.
Let $Q$ be a finite quiver without oriented cycles and let $\nu(Q)$
be the supremum of the multiplicities of the arrows occurring in
all quivers mutation-equivalent to $Q$.

\begin{theorem}
\begin{itemize}
\item[a)] $Q$ is representation-finite iff $\nu(Q)$ equals $1$.
\item[b)] $Q$ is tame iff $\nu(Q)$ equals $2$.
\item[c)] $Q$ is wild iff $\nu(Q)\geq 3$ iff $\nu(Q)=\infty$.
\item[d)] The mutation class of $Q$ is finite iff $Q$ has two vertices,
is representation-finite or tame.
\end{itemize}
\end{theorem}

Here, part a) follows from Gabriel's theorem and part (iii) of
Theorem~1.8 in \cite{FominZelevinsky03}.
Part b) follows from parts a) and c) by exclusion of the third.
For part c), let us first assume that $Q$ is wild. Then
it is proved at the end of the proof of theorem~3.1 in \cite{BuanReiten06} that
$\nu(Q)=\infty$. Conversely, let us assume that $\nu(Q)\geq 3$. Then using
Theorem~5 of \cite{CalderoKeller06} we obtain that $Q$ is wild.
Part d) is proved in \cite{BuanReiten06}.

\subsection{The Caldero-Chapoton formula} \label{ss:Caldero-Chapoton-formula}
Let $\Delta$ be a simply laced Dynkin diagram and $Q$ a quiver with
underlying graph $\Delta$. Suppose that the set of vertices of
$\Delta$ and $Q$ is the set of the natural numbers $1$, $2$, \ldots, $n$.
We already know from part b) of theorem~\ref{thm:cluster-finite-classification}
that for each positive root
\[
\alpha=\sum_{i=1}^n d_i \alpha_i
\]
of the corresponding root system, there is a unique non initial
cluster variable $X_\alpha$ with denominator
\[
x_1^{d_1} \ldots x_n^{d_n}.
\]
By combining this with Gabriel's theorem, we get the

\begin{corollary} The map taking an indecomposable representation $V$
with dimension vector $(d_i)$ of $Q$ to the unique non initial cluster
variable $X_V$ whose denominator is $x_1^{d_1} \ldots x_n^{d_n}$
induces a bijection from the set of ismorphism classes of indecomposable
representations to the set of non initial cluster variables.
\end{corollary}

Let us consider this bijection for $Q=\vec{A}_2$:
\begin{align*}
S_2 &=(0 \to k) & P_1 &=(k\to k) & S_1 &=(k \to 0) \\
X_{S_2} &=\frac{1+x_1}{x_2} & X_{P_1} &=\frac{x_1+1+x_2}{x_1 x_2} & S_{S_1} &=\frac{1+x_2}{x_1}
\end{align*}
We observe that for the two simple representations, the
numerator contains exactly two terms: the number of subrepresentations
of the simple representation! Moreover,
the representation $P_1$ has exactly three subrepresentations and
the numerator of $X_{P_1}$ contains three terms. In fact, it turns out that this
phenomenon is general in type $A$. But now let us consider the
following quiver with underlying graph $D_4$
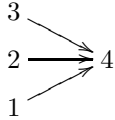
\[
\xymatrix@R=0.2cm{ 3 \ar[rd] & \\ 2 \ar[r] & 4 \\ 1 \ar[ur] & }
\]
and the dimension vector $d$ with $d_1=d_2=d_3=1$ and $d_4=2$. The
unique (up to isomorphism) indecomposable representation $V$ with
dimension vector $d$ consists of a plane $V_4$ together with
three lines in general position $V_i\subset V_4$, $i=1,2,3$.
The corresponding cluster variable is
\[
X_4 = \frac{1}{x_1x_2x_3x_4^2}\,(1+3x_4+3x_4^2+x_4^3+2x_1x_2x_3+3x_1x_2x_3x_4+x_1^2x_2^2x_3^2).
\]
Its numerator contains a total of $14$ monomials. On the other
hand, it is easy to see that $V_4$ has only $13$ types of submodules:
twelve submodules are determined by their dimension vectors but
for the dimension vector $e=(0,0,0,1)$, we have a family of
submodules: Each submodule of this dimension vector corresponds to the choice
of a line in $V_4$. Thus for this dimension vector $e$,
the family of submodules is parametrized by a projective
line. Notice that the Euler characteristic of the projective
line is $2$ (since it is a sphere: the Riemann sphere).
So if we attribute weight $1$ to the submodules
determined by their dimension vector and weight $2$ to
this $\mathbb{P}^1$-family, we find a `total submodule weight' equal to
the number of monomials in the numerator. These considerations
led Caldero-Chapoton \cite{CalderoChapoton06} to the following definition,
whose ingredients we describe below: Let $Q$ be a finite quiver with vertices $1$, \ldots, $n$,
and $V$ a finite-dimensional representation of $Q$. Let $d$ be the
dimension vector of $V$. Define
\[
CC(V)=\frac{1}{x_1^{d_1} x_2^{d_2} \ldots x_n^{d_n}}\,(\sum_{0\leq e \leq d}
\chi(\mbox{Gr}_e(V)) \prod_{i=1}^n x_i^{\sum_{j\to i} e_j + \sum_{i\to j}(d_j-e_j)} ).
\]
Here the sum is taken over all
vectors $e\in \N^n$ such that $0 \leq e_i \leq d_i$ for all $i$.
For each such vector $e$,
the {\em quiver Grassmannian} $\mbox{Gr}_e(V)$ is the variety
of $n$-tuples of subspaces $U_i \subset V_i$ such that $\dim U_i=e_i$
and the $U_i$ form a subrepresentation of $V$. By taking such a
subrepresentation to the family of the $U_i$, we obtain a
map
\[
\mbox{Gr}_e(V) \to \prod_{i=1}^n \mbox{Gr}_{e_i}(V_i) \ko
\]
where $\mbox{Gr}_{e_i}(V_i)$ denotes the ordinary Grassmannian
of $e_i$-dimensional subspaces of $V_i$. Recall that the Grassmannian
carries a canonical structure of projective variety. It is not
hard to see that for a family of subspaces $(U_i)$ the condition
of being a subrepresentation is a closed condition so that
the quiver Grassmannian identifies with a projective subvariety
of the product of ordinary Grassmannians. If $k$ is the field
of complex numbers, the Euler characteristic $\chi$ is taken
with respect to singular cohomology with coefficients in $\Q$
(or any other field). If $k$ is an arbitrary algebraically
closed field, we use \'etale cohomology to define $\chi$.
The most important properties of $\chi$ are (\cf \eg section~7.4 in
\cite{GeissLeclercSchroeer07a})
\begin{itemize}
\item[(1)] $\chi$ is additive with respect to disjoint unions;
\item[(2)] if $p: E \to X$ is a morphism of algebraic varieties
such that the Euler characteristic of the fiber over a point
$x\in X$ does not depend on $x$, then $\chi(E)$ is the product
of $\chi(X)$ by the Euler characteristic of the fiber over any
point $x\in X$.
\end{itemize}

\begin{theorem}[Caldero-Chapoton \protect{\cite{CalderoChapoton06}}]
Let $Q$ be a Dynkin quiver and $V$ an indecomposable representation.
Then we have $CC(V)=X_V$, the cluster variable obtained from $V$
by composing Fomin-Zelevinsky's bijection with Gabriel's.
\end{theorem}

Caldero-Chapoton's proof of the theorem was by induction.
One of the aims of the following sections is to explain
`on what' they did the induction.

\subsection{The derived category} \label{ss:derived-category}

Let $k$ be an algebraically closed field and $Q$ a (possibly
infinite) quiver without oriented cycles (we will impose more
restrictive conditions on $Q$ later).
For example, $Q$ could be the quiver
\[
\xymatrix{1 \ar[d]_\gamma \ar[r] & 3 \ar[d]^\alpha\\
2 \ar[ur]_\beta \ar[r] & 4}
\]
A {\em path} of $Q$ is a formal composition of $\geq 0$
arrows. For example, the sequence $(4|\alpha|\beta|\gamma|1)$
is a path of length $3$ in the above example (notice that we
include the source and target vertices of the path in the
notation). For each vertex $i$ of $Q$, we have the {\em lazy path}
$e_i=(i|i)$, the unique path of length $0$ which starts at $i$
and stops at $i$ and does nothing in between. The {\em path category}
has set of objects $Q_0$ (the set of vertices of $Q$) and,
for any vertices $i$, $j$, the morphism space from $i$ to $j$
is the vector space whose basis consists of all paths from $i$
to $j$. Composition is induced by composition of paths and
the unit morphisms are the lazy paths. If $Q$ is finite,
we define the {\em path algebra} to be the matrix algebra
\[
kQ=\bigoplus_{i,j\in Q_0} \Hom(i,j)
\]
where multiplication is matrix multiplication. Equivalently,
the path algebra has as a basis all paths and its product
is given by concatenating composable paths and equating
the product of non composable paths to zero. The path
algebra has the sum of the lazy paths as its unit element
\[
1 = \sum_{i\in Q_0} e_i.
\]
The idempotent $e_i$ yields the {\em projective right module}
\[
P_i = e_i kQ.
\]
The modules $P_i$ generate the {\em category of $k$-finite-dimensional
right modules $\mod kQ$}. Each arrow $\alpha$ from $i$ to
$j$ yields a map $P_i \to P_j$ given by left multiplication
by $\alpha$. (If we were to consider -- heaven forbid -- left
modules, the analogous map would be given by right
multiplication by $\alpha$ and it would go in the
direction opposite to that of $\alpha$. Whence our
preference for right modules).

Notice that we have an equivalence of categories
\[
\rep_k(Q^{op}) \to \mod kQ
\]
sending a representation $V$ of the opposite quiver $Q^{op}$
to the sum
\[
\bigoplus_{i\in Q_0} V_i
\]
endowed with the natural right action of the path algebra.
Conversely, a $kQ$-module $M$ gives rise to the representation
$V$ with $V_i=M e_i$ for each vertex $i$ of $Q$ and $V_\alpha$ given by
right multiplication by $\alpha$ for each arrow $\alpha$ of $Q$.
The category $\mod kQ$ is {\em abelian}, \ie it is additive,
has kernels and cokernels and for each morphism $f$ the
cokernel of its kernel is canonically isomorphic to the
kernel of its cokernel.

The category $\mod kQ$ is {\em hereditary}. Recall
from \cite{CartanEilenberg56} that this means that
submodules of projective modules are projective; equivalently,
that all extension groups in degrees $i\geq 2$ vanish:
\[
\Ext^i_{kQ}(L,M) =0;
\]
equivalently, that $kQ$ is of global dimension $\leq 1$; \ldots
Thus, in the spirit of noncommutative algebraic geometry approached
via abelian categories, we should think of $\mod kQ$ as a
`non commutative curve'.

We define $\cd_Q$ to be the {\em bounded derived category} $\cd^b(\mod kQ)$
of the abelian category $\mod kQ$. Thus, the objects of $\cd_Q$ are the
bounded complexes of (right) $kQ$-modules
\[
\ldots \to 0 \to \ldots \to M^p \arr{d^p} M^{p+1} \to \ldots \to 0 \to \ldots .
\]
Its morphisms are obtained from morphisms of complexes by
formally inverting all quasi-isomorphisms. We refer to
\cite{Verdier96} \cite{KashiwaraSchapira94} \ldots for
in depth treatments of the fundamentals of this construction.
Below, we will give a complete and elementary description of
the category $\cd_Q$ if $Q$ is a Dynkin quiver. We have the following
general facts: The functor
\[
\mod kQ \to \cd_Q
\]
taking a module $M$ to the complex concentrated in
degree $0$
\[
\ldots \to 0 \to M \to 0 \to \ldots
\]
is a fully faithful embedding. From now on, we will identify
modules with complexes concentrated in degree $0$. If $L$
and $M$ are two modules, then we have a canonical
isomorphism
\[
\Ext^i_{kQ}(L,M) \iso \Hom_{\cd_Q}(L, M[i])
\]
for all $i\in\Z$, where $M[i]$ denotes the complex $M$
shifted by $i$ degrees to the left: $M[i]^p=M^{p+i}$, $p\in\Z$,
and endowed with the differential $d_{M[i]}=(-1)^i d_M$.
The category $\cd_Q$ has all finite direct sums (and they
are given by direct sums of complexes) and the decomposition
theorem~\ref{thm:decomposition} holds. Moreover, each
object is isomorphic to a direct sum of shifted copies
of modules (this holds more generally in the derived
category of any hereditary abelian category, for example
the derived category of coherent sheaves on an algebraic
curve). The category $\cd_Q$ is abelian if and only if the
quiver $Q$ does not have any arrows. However, it is
always {\em triangulated}. This means that it is
$k$-linear (it is additive, and the morphism sets
are endowed with $k$-vector space structures so that
the composition is bilinear) and endowed with
the following extra structure:
\begin{itemize}
\item[a)] a {\em suspension (or shift) functor} $\Sigma: \cd_Q \to \cd_Q$,
namely the functor taking a complex $M$ to $M[1]$;
\item[b)] a class of {\em triangles} (sometimes called `distinguished
triangles'), namely the sequences
\[
L \to M \to N \to \Sigma L
\]
which are `induced' by short exact sequences of complexes.
\end{itemize}
The class of triangles satisfies certain axioms, \cf \eg \cite{Verdier96}.
The most important consequence of these axioms is that the triangles
induce long exact sequences in the functors $\Hom(X,?)$ and
$\Hom(?,X)$, \ie for each object $X$ of $\cd_Q$, the sequences
\[
\ldots (X, \Sigma^{-1} N) \to (X,L) \to (X,M) \to (X,N) \to (X,\Sigma L) \to \ldots
\]
and
\[
\ldots (\Sigma^{-1} N, X) \la (L,X) \la (M,X) \la (N,X) \la (\Sigma L,X) \la \ldots
\]
are exact.

\subsection{Presentation of the derived category of a Dynkin quiver}
\label{ss:der-cat-Dynkin-quiver}
From now on, we assume that $Q$ is a Dynkin quiver.
Let $\Z Q$ be its repetition (\cf section~\ref{ss:knitting-algorithm}).
So the vertices of $\Z Q$ are the pairs $(p,i)$, where $p$ is an integer
and $i$ a vertex of $Q$ and the arrows of $\Z Q$ are obtained as follows:
each arrow $\alpha: i\to j$ of $Q$ yields the arrows
\[
(p,\alpha) : (p,i) \to (p,j) \ko p\in \Z\ko
\]
and the arrows
\[
\sigma(p, \alpha) : (p-1, j) \to (p,i) \ko p\in \Z.
\]
We extend $\sigma$ to a map defined on all arrows of $\Z Q$ by
defining
\[
\sigma(\sigma(p,\alpha)) = (p-1, \alpha).
\]
We endow $\Z Q$ with the map $\sigma$ and with the automorphism
$\tau: \Z Q \to \Z Q$ taking $(p,i)$ to $(p-1, i)$ and $(p,\alpha)$
to $(p-1,\alpha)$ for all vertices $i$ of $Q$, all arrows $\alpha$
of $Q$ and all integers $p$.

For a vertex $v$ of $\Z Q$ the {\em mesh ending at $v$} is the
full subquiver
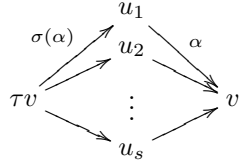
\begin{equation}\label{eq:mesh}
\xymatrix@R=0.1cm{ & u_1 \ar[rdd]^\alpha & \\
& u_2 \ar[rd] & \\
\tau v \ar[ruu]^{\sigma(\alpha)} \ar[ru] \ar[rd] & \vdots & v \\
 & u_s \ar[ru] &
}
\end{equation}
formed by $v$, $\tau(v)$ and all sources $u$ of arrows
$\alpha: u\to v$ of $\Z Q$ ending in $v$.
We define the {\em mesh ideal} to be the (two-sided) ideal
of the path category of $\Z Q$ which is generated by
all {\em mesh relators}
\[
r_v = \sum_{\mbox{\tiny arrows }\alpha: u\to v} \alpha \sigma(\alpha) \ko
\]
where $v$ runs through the vertices of $\Z Q$. The {\em mesh category}
is the quotient of the path category of $\Z Q$ by the mesh ideal.

\begin{theorem}[Happel \protect{\cite{Happel87}}]
\begin{itemize}
\item[a)] There is a canonical bijection $v\mapsto M_v$ from the set
of vertices of $\Z Q$ to the set of isomorphism classes of indecomposables of
$\cd_Q$ which takes the vertex $(1,i)$ to the indecomposable projective $P_i$.
\item[b)] Let $\ind \cd_Q$ be the full subcategory of indecomposables
of $\cd Q$. The bijection of a) lifts to an equivalence of
categories from the mesh category of $\Z Q$ to the category $\ind \cd_Q$.
\end{itemize}
\end{theorem}

\begin{figure}[p]
\[
\begin{xy}<1.2cm,0cm>:
(9,2.5)*{\ind \cd_{\vec{A}_5}},
(0,0)*\xybox{<1.2cm,0cm>:<0.6cm,0.6cm>::
@={(0,0),(1,0),(2,0),(3,0),(4,0),(5,0),(6,0),(7,0),(8,0),(9,0)}
@@{*\xybox{<1.2cm,0cm>:<0.6cm,0.6cm>::
0*\cir<2pt>{}="A"*\cir<2pt>{}; "A"+(0,1)*\cir<2pt>{} **\dir{-}
?>*\dir{>}, "A"+(-1,2)
*\cir<2pt>{}="B"; "B" +(1,-1)*\cir<2pt>{} **\dir{-} ?>*\dir{>},
"B"+ (0,1)*\cir<2pt>{} **\dir{-} ?>*\dir{>},
"A"+(-2,4)*\cir<2pt>{}="B"; "B" + (1,-1)*\cir<2pt>{} **\dir{-}
?>*\dir{>}, "A"+(0,1)
*\cir<2pt>{}="B"; "B" +(1,-1)*\cir<2pt>{} **\dir{-} ?>*\dir{>},
"B"+ (0,1)*\cir<2pt>{} **\dir{-} ?>*\dir{>}, "A"+(-1,3)
*\cir<2pt>{}="B"; "B" +(1,-1)*\cir<2pt>{} **\dir{-} ?>*\dir{>},
"B"+ (0,1)*\cir<2pt>{} **\dir{-} ?>*\dir{>},}},
@={(1,0),(1,1),(1,2),(1,3),(1,4),(2,0),(2,1),(2,2),(2,3),(3,0),(3,1),(3,2),(4,0),(4,1),(5,0)},
@@{*{\bt}}, @={(7,0), (7,1),(7,2),(7,3), (7,4), (8,0), (8,1), (8,2),
(8,3), (9,0), (9,1), (9,2)}, @@{*{\bt}},
@={(1,2),(1,3),(1,4),(2,1),(2,2),(2,3),(2,4),(3,0),(3,1),(3,2),(3,3),(3,4),(4,0),(4,1),(5,0)},
(1,1)*!L(2){{\scriptstyle P_2}},
(1,2)*!L(2){{\scriptstyle P_3}},
(3,3)*!L(1.5){\scriptstyle{\Sigma P_2}},
(7,1)*!L(1.5){\scriptstyle{\Sigma^2 P_2}},
(1,0)*!R(1.6){{\scriptstyle (1,1)}}, (1,1)*!R(1.6){{\scriptstyle (1,2)}},
(1,2)*!R(1.6){{\scriptstyle (1,3)}}, (1,3)*!R(1.6){{\scriptstyle
(1,4)}}, (1,4)*!R(1.6){{\scriptstyle (1,5)}},
(2,0)*!R(1.6){{\scriptstyle (2,1)}}, }
\end{xy}
\]
\caption{The repetition of type $A_n$} \label{fig:an}
\end{figure}
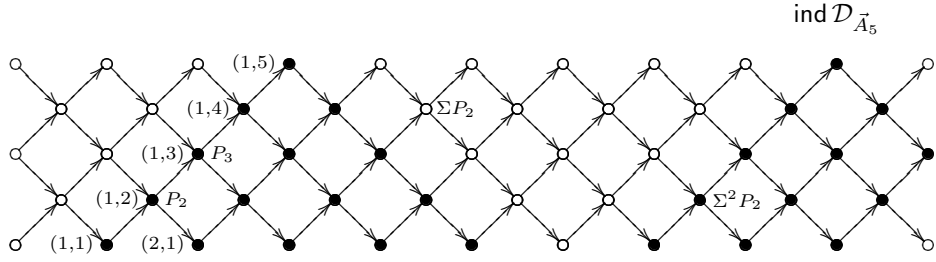

In figure~\ref{fig:an}, we see the repetition for $Q=\vec{A}_5$
and the map taking its vertices to the indecomposable objects
of the derived category. The vertices marked $\bt$ belonging
to the left triangle are mapped to indecomposable modules.
The vertex $(1,i)$ corresponds to the indecomposable projective
$P_i$. The arrow $(1,i)\to (1,i+1)$, $1\leq i\leq 5$, is
mapped to the left multiplication by the arrow $i \to i+1$.
The functor takes a mesh~(\ref{eq:mesh}) to a triangle
\begin{equation} \label{eq:ARtriangle}
\xymatrix{ M_{\tau v} \ar[r] & \bigoplus_{i=1}^s M_{u_i} \ar[r] & M_v \ar[r] & \Sigma M_{\tau v}}
\end{equation}
called an {\em Auslander-Reiten triangle} or {\em almost split
triangle}, \cf \cite{Happel88}. If $M_v$ and $M_{\tau v}$ are modules, then so
is the middle term and the triangle
comes from an exact sequence of modules
\[
\xymatrix{0 \ar[r] & M_{\tau v} \ar[r] & \bigoplus_{i=1}^s M_{u_i} \ar[r] & M_v \ar[r] & 0}
\]
called an {\em Auslander-Reiten sequence} or {\em almost split sequence},
\cf \cite{AuslanderReitenSmaloe95}.
These almost split triangles respectively sequences can be
characterized intrinsically in $\cd_Q$ respectively $\mod kQ$.

\begin{figure}
\[
\begin{xy} 0;<0.5pt,0pt>:<0pt,-0.5pt>::
(0,193) *+{10000} ="0", (49,147) *+{11000} ="1", (96,98) *+{11100} ="2",
(145,49) *+{11110} ="3", (194,0) *+{11111} ="4", (96,193) *+{01000} ="5",
(145,147) *+{01100} ="6", (194,98) *+{01110} ="7", (243,49) *+{01111} ="8",
(292,0) *+{-10000} ="9", (194,193) *+{00100} ="10", (243,147) *+{00110} ="11",
(292,98) *+{00111} ="12", (341,49) *+{-11000} ="13", (389,0) *+{-01000} ="14",
(292,193) *+{00010} ="15", (341,147) *+{00011} ="16", (389,98) *+{-11100}
="17", (438,49) *+{-01100} ="18", (487,0) *+{-00100} ="19", (389,193) *+{00001}
="20", (438,147) *+{-11110} ="21", (487,98) *+{-01110} ="22", (536,49)
*+{-00110} ="23", (585,0) *+{-00010} ="24", (487,193) *+{-11111} ="25",
(536,147) *+{-01111} ="26", (585,98) *+{-00111} ="27", (632,49) *+{-00011}
="28", (681,0) *+{-00001} ="29", "0", {\ar"1"}, "1", {\ar"2"}, "1",
{\ar"5"}, "2", {\ar"3"}, "2", {\ar"6"}, "3", {\ar"4"}, "3",
{\ar"7"}, "4", {\ar"8"}, "5", {\ar"6"}, "6", {\ar"7"}, "6",
{\ar"10"}, "7", {\ar"8"}, "7", {\ar"11"}, "8", {\ar"9"}, "8",
{\ar"12"}, "9", {\ar"13"}, "10", {\ar"11"}, "11", {\ar"12"}, "11",
{\ar"15"}, "12", {\ar"13"}, "12", {\ar"16"}, "13", {\ar"14"}, "13",
{\ar"17"}, "14", {\ar"18"}, "15", {\ar"16"}, "16", {\ar"17"}, "16",
{\ar"20"}, "17", {\ar"18"}, "17", {\ar"21"}, "18", {\ar"19"}, "18",
{\ar"22"}, "19", {\ar"23"}, "20", {\ar"21"}, "21", {\ar"22"}, "21",
{\ar"25"}, "22", {\ar"23"}, "22", {\ar"26"}, "23", {\ar"24"}, "23",
{\ar"27"}, "24", {\ar"28"}, "25", {\ar"26"}, "26", {\ar"27"}, "27",
{\ar"28"}, "28", {\ar"29"},
\end{xy}
\]
\caption{Some dimension vectors of indecomposables in $\cd_{\vec{A}_5}$}\
\label{fig:dimvectA5}
\end{figure}
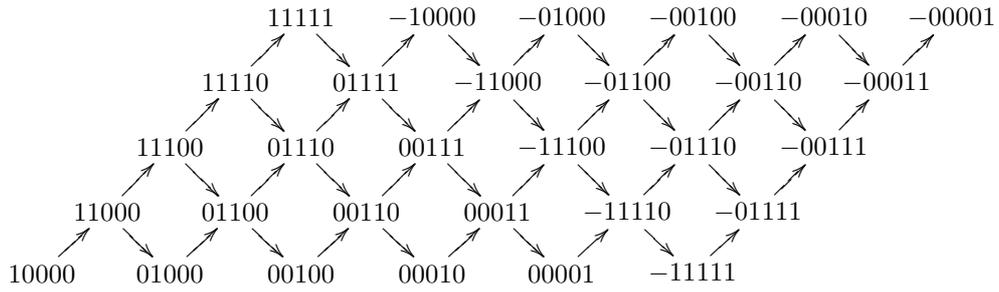

Recall that the Grothendieck group $K_0(\ct)$ of a triangulated
category is the quotient of the free abelian group on the
isomorphism classes $[X]$ of objects $X$ of $\ct$ by the
subgroup generated by all elements
\[
[X] -[Y] + [Z]
\]
arising from triangles $(X,Y,Z)$ of $\ct$. In the case of
$\cd_Q$, the natural map
\[
K_0(\mod kQ) \to K_0(\cd_Q)
\]
is an isomorphism (its inverse sends a complex to the alternating
sum of the classes of its homologies). Since $K_0(\mod kQ)$ is
free on the classes $[S_i]$ associated with the simple modules,
the same holds for $K_0(\cd_Q)$ so that its elements are given
by $n$-tuples of integers. We write $\dimv M$ for the image
in $K_0(\cd_Q)$ of an object $M$ of $K_0(\cd_Q)$ and call
$\dimv M$ the {\em dimension vector of $M$}. Then each
triangle~(\ref{eq:ARtriangle}) yields an equality
\[
\dimv M_v = \sum_{i=1}^s \dimv M_{u_i} - \dimv M_{\tau v}.
\]
Using these equalities, we can easily determine $\dimv M$
for each indecomposable $M$ starting from the known
dimension vectors $\dimv P_i$, $1\leq i\leq n$.
In the above example, we find the dimension vectors listed
in figure~(\ref{fig:dimvectA5}).

Thanks to the theorem, the automorphism $\tau$ of the
repetition yields a $k$-linear automorphism, still
denoted by $\tau$, of the derived category $\cd_Q$.
This automorphism has several intrinsic descriptions:

1) As shown in \cite{Gabriel80}, it is the right derived functor of the left
exact Coxeter functor $\rep(Q^{op}) \to \rep(Q^{op})$
introduced by Bernstein-Gelfand-Ponomarev
\cite{BernsteinGelfandPonomarev73}
in their proof of Gabriel's theorem. If we identify
$K_0(\cd_Q)$ with the root lattice via Gabriel's theorem,
then the automorphism induced by $\tau^{-1}$ equals the
the Coxeter transformation $c$. As shown by Gabriel \cite{Gabriel80},
the identity $c^h=\id$, where $h$ is the Coxeter number,
lifts to an isomorphism of functors
\begin{equation} \label{eq:tau-h-equals-Sigma2}
\tau^{-h} \iso \Sigma^2.
\end{equation}

2) It can be expressed in terms of the Serre functor
of $\cd_Q$: Recall that for a $k$-linear triangulated
category $\ct$ with finite-dimensional morphism spaces,
a Serre functor is an autoequivalence $S: \ct \to \ct$
such that the Serre duality formula holds: We have
bifunctorial isomorphisms
\[
D\Hom(X,Y) \iso \Hom(Y, SX) \ko X,Y\in \ct\ko
\]
where $D$ is the duality $\Hom_k(?,k)$ over the ground field.
Notice that this determines the functor $S$ uniquely
up to isomorphism.
In the case of $\cd_Q=\cd^b(\mod kQ)$, it is not hard to
prove that a Serre functor exists (it is given by the left
derived functor of the tensor product by the bimodule
$D(kQ)$). Now the autoequivalence $\tau$, the
suspension functor $\Sigma$ and the Serre functor $S$
are linked by the fundamental isomorphism
\begin{equation} \label{eq:tau-Sigma-equals-S}
\tau \Sigma \iso S.
\end{equation}

\subsection{Caldero-Chapoton's proof}
The above description of the derived category yields
in particular a description of the module category, which
is a full subcategory of the derived category. This description
was used by Caldero-Chapoton \cite{CalderoChapoton06} to
prove their formula. Let us sketch the main steps in their
proof: Recall that we have defined a surjective map
$v \mapsto X_v$ from the set of vertices of the repetition
to the set of cluster variables such that
\begin{itemize}
\item[a)] we have $X_{(0,i)} = x_i$ for $1\leq i\leq n$ and
\item[b)] we have
\[
X_{\tau v} X_v = 1 + \prod_{\mbox{\tiny arrows } w \to v} X_w
\]
for all vertices $v$ of the repetition.
\end{itemize}

We wish to show that
we have
\[
X_{v} = CC(M_v)
\]
for all vertices $v$ such that $M_v$ is an indecomposable module.
This is done by induction on the distance of $v$ from the
vertices $(1,i)$ in the quiver $\Z Q$. More precisely, one
shows the following
\begin{itemize}
\item[a)] We have $CC(P_i)=X_{(1,i)}$ for each indecomposable
projective $P_i$. Here we use the fact that submodules of projectives
are projective in order to explicitly compute $CC(P_i)$.
\item[b1)] For each split exact sequence
\[
0 \to L \to E \to M \to 0 \ko
\]
we have
\[
CC(L) CC(M) = CC(E).
\]
Thus, if $E = E_1 \oplus \ldots E_s$ is a decomposition into
indecomposables, then
\[
CC(E) = \prod_{i=1}^s CC(E_i).
\]
\item[b2)] If
\[
0 \to L \to E \to M \to 0
\]
is an almost split exact sequence, then we have
\[
CC(E) +1 = CC(L) CC(M).
\]
\end{itemize}
It is now clear how to prove the equality $X_{v}=CC(M_v)$
by induction by proceeding from the projective indecomposables
to the right.

\subsection{The cluster category}
The {\em cluster category}
\[
\cc_Q = \cd_Q/(\tau^{-1}\Sigma)^{\Z} = \cd_Q/(S^{-1}\Sigma^2)^\Z
\]
is the orbit category of the derived category under the action of
the cyclic group generated by the autoequivalence
$\tau^{-1}\Sigma = S^{-1} \Sigma^2$. This means that the objects
of $\cc_Q$ are the same as those of the derived category $\cd_Q$ and
that for two objects $X$ and $Y$, the morphism space from $X$ to
$Y$ in $\cc_Q$ is
\[
\cc_Q(X,Y)=\bigoplus_{p\in \Z} \cd_Q(X, (S^{-1} \Sigma^2)^p Y).
\]
Morphisms are composed in the natural way.
This definition is due to
Buan-Marsh-Reineke-Reiten-Todorov \cite{BuanMarshReinekeReitenTodorov06},
who were trying to obtain a better understanding of the `decorated
quiver representations' introduced by Reineke-Marsh-Zelevinsky
\cite{MarshReinekeZelevinsky03}.
For quivers of type $A$, an equivalent category was defined independently
by Caldero-Chapoton-Schiffler \cite{CalderoChapotonSchiffler06}
using an entirely different description. Clearly the
category $\cc_Q$ is $k$-linear. It is not
hard to check that its morphism spaces are finite-dimensional.

One can show \cite{Keller05} that $\cc_Q$ admits a canonical structure of
triangulated category such that the projection functor
$\pi: \cd_Q \to \cc_Q$ becomes a triangle functor (in
general, orbit categories of triangulated categories are
no longer triangulated). The Serre functor $S$ of $\cd_Q$
clearly induces a Serre functor in $\cc_Q$, which we still
denote by $S$. Now, by the definition of $\cc_Q$ (and its
triangulated structure), we have an isomorphism of triangle
functors
\[
S \iso \Sigma^2.
\]
This means that $\cc_Q$ is {\em $2$-Calabi-Yau}. Indeed,
for an integer $d\in\Z$, a triangulated category $\ct$ with finite-dimensional
morphism spaces is $d$-Calabi-Yau if it admits a Serre
functor isomorphic as a triangle functor to the $d$th
power of its suspension functor.

\subsection{From cluster categories to cluster algebras}
\label{ss:from-cluster-categories-to-cluster-algebras}
We keep the notations and hypotheses of the previous section.
The suspension functor $\Sigma$  and the Serre functor $S$
induce automorphisms of the repetition $\Z Q$ which we still
denote by $\Sigma$ and $S$ respectively. The orbit quiver
$\Z Q/(\tau^{-1} \Sigma)^\Z$ inherits the automorphism
$\tau$ and the map $\sigma$ (defined on arrows only) and
thus has a well-defined mesh category. Recall that we
write $\Ext^1(X,Y)$ for $\Hom(X,\Sigma Y)$ in any
triangulated category.

\begin{theorem}[\protect{\cite{BuanMarshReinekeReitenTodorov06} \cite{BuanMarshReiten08}}]
\label{thm:cluster-tilting}
\begin{itemize}
\item[a)] The decomposition theorem holds for the cluster
category and the mesh category of $\Z Q/(\tau^{-1}\Sigma)^\Z$
is canonically equivalent to the full subcategory $\ind \cc_Q$
of the indecomposables of $\cc_Q$. Thus, we have an induced
bijection $L \mapsto X_L$ from the set of isomorphism classes
of indecomposables of $\cc_Q$ to the set of all cluster
variables of $\ca_Q$ which takes the shifted projective
$\Sigma P_i$ to the initial variable $x_i$, $1 \leq i\leq n$.
\item[b)] Under this bijection, the clusters correspond
to the {\em cluster-tilting sets}, \ie the sets
of pairwise non isomorphic indecomposables $T_1, \ldots, T_n$
such that we have
\[
Ext^1(T_i, T_j)=0
\]
for all $i,j$.
\item[c)] If $T_1, \ldots, T_n$ is cluster-tilting,
then the quiver (\cf below) of the endomorphism algebra
of the sum $T=\bigoplus_{i=1}^n T_i$ does not have loops
nor $2$-cycles and the associated antisymmetric matrix
is the exchange matrix of the unique seed containing
the cluster $X_{T_1}, \ldots, X_{T_n}$.
\end{itemize}
\end{theorem}

In part b), the condition implies in particular that
$\Ext^1(T_i, T_i)$ vanishes. However, for a Dynkin quiver
$Q$, we have $\Ext^1(L,L)=0$ for each indecomposable $L$
of $\cc_Q$. A {\em cluster-tilting object of $\cc_Q$}
is the direct sum of the objects $T_1, \ldots, T_n$ of a
cluster-tilting set. Since these are pairwise non-isomorphic
indecomposables, the datum of $T$ is equivalent to that of the $T_i$.
A {\em cluster-tilted algebra of type $Q$} is the endomorphism algebra of
a cluster-tilting object of $\cc_Q$.
In part c), the most subtle point is that
the quiver does not have loops or $2$-cycles
\cite{BuanMarshReiten08}.
Let us recall what one means by the {\em quiver of
a finite-dimensional algebra over an algebraically closed
field}:

\begin{proposition-definition}[Gabriel] \label{prop:quiver}
Let $B$ be a finite-dimensional
algebra over the algebraically closed ground field $k$.
\begin{itemize}
\item[a)] There exists a quiver $Q_B$ , unique up to isomorphism,
such that $B$ is Morita equivalent to the algebra $kQ_B/I$, where
$I$ is an ideal of $kQ_B$ contained in the square of the ideal generated
by the arrows of $Q_B$.
\item[b)] The ideal $I$ is not unique in general but we have
$I=0$ iff $B$ is hereditary.
\item[c)] There is a bijection $i \mapsto S_i$ between the vertices
of $Q_B$ and the isomorphism classes of simple $B$-modules. The
number of arrows from a vertex $i$ to a vertex $j$ equals the
dimension of $\Ext^1_B(S_j, S_i)$.
\end{itemize}
\end{proposition-definition}
In our case, the algebra $B$ is the endomorphism algebra of the
sum $T$ of the cluster-tilting set $T_1$, \ldots, $T_n$
in $\cc_Q$. In this case, the Morita equivalence of a) even becomes
an isomorphism (because the $T_i$ are pairwise non isomorphic).
For a suitable choice of this isomorphism, the idempotent $e_i$
associated with the vertex $i$ is sent to the identity of $T_i$ and
the images of the arrows from $i$ to $j$ yield a basis of the space
of {\em irreducible morphisms}
\[
\irr_T(T_i,T_j)=\rad_T(T_i,T_j)/\rad_T^2(T_i,T_j) \ko
\]
where $\rad_T(T_i,T_j)$ denotes the vector space of non isomorphisms
from $T_i$ to $T_j$ (thanks to the locality of the endomorphism rings,
this set is indeed closed under addition) and $\rad_T^2$ the subspace
of non isomorphisms admitting a non trivial factorization:
\[
\rad_T^2(T_i, T_j) = \sum_{r=1}^n \rad_T(T_r, T_j) \rad_T(T_i, T_r).
\]

\begin{figure}[p]
\[
\begin{xy} 0;<1pt,0pt>:<0pt,-1pt>::
(0,100) *+{0} ="0",
(25,75) *+{1} ="1",
(50,50) *+{2} ="2",
(75,25) *+{3} ="3",
(100,0) *+{4} ="4",
(50,100) *+{T_1} ="5",
(75,75) *+{6} ="6",
(100,50) *+{T_2} ="7",
(125,25) *+{8} ="8",
(150,0) *+{T_3} ="9",
(100,100) *+{10} ="10",
(125,75) *+{11} ="11",
(150,50) *+{12} ="12",
(175,25) *+{13} ="13",
(150,100) *+{T_4} ="14",
(175,75) *+{15} ="15",
(200,50) *+{16} ="16",
(200,100) *+{17} ="17",
(225,75) *+{18} ="18",
(250,100) *+{T_5} ="19",
(200,0) *+{20} ="20",
(225,25) *+{21} ="21",
(250,50) *+{22} ="22",
(275,75) *+{23} ="23",
(300,100) *+{24} ="24",
"0", {\ar"1"},
"1", {\ar"2"},
"1", {\ar"5"},
"2", {\ar"3"},
"2", {\ar"6"},
"3", {\ar"4"},
"3", {\ar"7"},
"4", {\ar"8"},
"5", {\ar"6"},
"6", {\ar"7"},
"6", {\ar"10"},
"7", {\ar"8"},
"7", {\ar"11"},
"8", {\ar"9"},
"8", {\ar"12"},
"9", {\ar"13"},
"10", {\ar"11"},
"11", {\ar"12"},
"11", {\ar"14"},
"12", {\ar"13"},
"12", {\ar"15"},
"13", {\ar"16"},
"13", {\ar"20"},
"14", {\ar"15"},
"15", {\ar"16"},
"15", {\ar"17"},
"16", {\ar"18"},
"16", {\ar"21"},
"17", {\ar"18"},
"18", {\ar"19"},
"18", {\ar"22"},
"19", {\ar"23"},
"20", {\ar"21"},
"21", {\ar"22"},
"22", {\ar"23"},
"23", {\ar"24"},
\end{xy}
\]
\caption{A cluster-tilting set in $A_5$} \label{fig:cluster-tilting}
\end{figure}
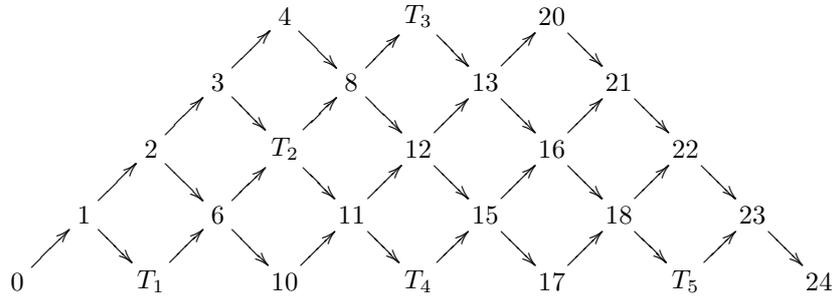

As an illustration of theorem~\ref{thm:cluster-tilting}, we consider
the cluster-tilting set $T_1, \ldots, T_5$ in $\cc_{\vec{A_5}}$
depicted in figure~\ref{fig:cluster-tilting}. Here the vertices
labeled $0$, $1$, \ldots, $4$ have to be identified with the
vertices labeled $20$, $21$, \ldots, $24$ (in this order) to
obtain the orbit quiver $\Z Q/(\tau^{-1} \Sigma)^\Z$. In the
orbit category, we have $\tau \iso \Sigma$ so that $\Sigma T_1$
is the indecomposable associated to vertex $0$, for example.
Using this and the description of the morphisms in the mesh
category, it is easy to check that we do have
\[
\Ext^1(T_i, T_j)=0
\]
for all $i,j$. It is also easy to determine the spaces
of morphisms
\[
\Hom_{\cc_Q}(T_i,T_j)
\]
and the compositions of morphisms. Determining these is
equivalent to determining the endomorphism
algebra
\[
\End(T)=\Hom(T,T)=\bigoplus_{i,j} \Hom(T_i, T_j).
\]
This algebra is easily seen to be {\em isomorphic} to the algebra given
by the following quiver $Q'$
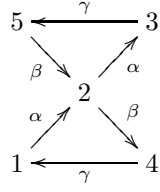
\[
\xymatrix@R=0.5cm@C=0.5cm{
5 \ar[rd]_\beta & & 3 \ar[ll]_\gamma \\
 & 2 \ar[ru]_\alpha \ar[rd]^\beta & \\
 1 \ar[ru]^\alpha & & 4 \ar[ll]^\gamma
 }
 \]
 with the relations
 \[
\alpha\beta=0 \ko \beta\gamma=0 \ko \gamma\alpha=0.
\]
Thus the quiver of $\End(T)$ is $Q'$. It encodes the
exchange matrix of the associated cluster
\begin{align*}
X_{T_1} &=\frac{1+x_2}{x_1} \\
X_{T_2} &=\frac{x_1 x_2+x_1x_4+x_3 x_4 + x_2 x_3 x_4}{x_1 x_2 x_3} \\
X_{T_3} &=\frac{x_1 x_2 x_3 + x_1 x_2 x_3 x_4 + x_1 x_2 x_5 +x_1 x_4 x_5 +
x_3 x_4 x_5 + x_2 x_3 x_4 x_5}{x_1 x_2 x_3 x_4 x_5} \\
X_{T_4} &=\frac{x_2+x_4}{x_3} \\
X_{T_5} &=\frac{1+x_4}{x_5}.
\end{align*}

\subsection{A $K$-theoretic interpretation of the exchange matrix}
Keep the notations and hypotheses of the preceding section. Let
$T_1, \ldots, T_n$ be a cluster-tilting set, $T$ the sum of the
$T_i$ and $B$ its endomorphism algebra. For two finite-dimensional
right $B$-modules $L$ and $M$ put
\[
\langle L, M \rangle_a = \dim \Hom(L,M) -\dim \Ext^1(L,M)
                        -\dim \Hom(M,L) +\dim \Ext^1(M,L).
\]
This is the antisymmetrization of a truncated Euler form. A priori
it is defined on the split Grothendieck group of the category
$\mod B$ (\ie the quotient of the free abelian group on the
isomorphism classes divided by the subgroup generated by
all relations obtained from direct sums in $\mod B$).

\begin{proposition}[Palu] The form $\langle, \rangle_a$ descends
to an antisymmetric form on $K_0(\mod B)$. Its matrix in the basis
of the simples is the exchange matrix associated with the cluster
corresponding to $T_1, \ldots, T_n$.
\end{proposition}

\subsection{Mutation of cluster-tilting sets}
\label{ss:mutation-cluster-tilting-sets} Let us recall two axioms of
triangulated categories:
\begin{itemize}
\item[TR1] For each morphism $u: X \to Y$, there exists a
triangle
\[
X \arr{u} Y \to Z \to \Sigma X.
\]
\item[TR2] A sequence
\[
X \arr{u} Y \arr{v} Z \arr{w} \Sigma X
\]
is a triangle if and only if the sequence
\[
Y \arr{v} Z \arr{w} \Sigma X \arr{-u} \Sigma Y
\]
is a triangle.
\end{itemize}
One can show that in TR1, the triangle is unique up to
(non unique) isomorphism. In particular, up to isomorphism,
the object $Z$ is uniquely determined by $u$. Notice the sign
in TR2. It follows from TR1 and TR2 that a given morphism
also occurs as the second (respectively third) morphism in
a triangle.

Now, with the notations and hypotheses of the preceding
section, suppose that $T_1, \ldots, T_n$ is a cluster-tilting
set and $Q'$ the quiver of the endomorphism algebra $B$ of
the sum of the $T_i$. As explained after
proposition-definition~\ref{prop:quiver},
we have a surjective
algebra morphism
\[
kQ' \to \bigoplus_{i,j} \Hom(T_i, T_j)
\]
which takes the idempotent $e_i$ to the identity of $T_i$
and the arrows $i \to j$ to irreducible morphisms
$T_i \to T_j$, for all vertices $i$, $j$ of $Q'$ (\cf
the above example computation of $B$ and $Q'=Q_B$).

Now let $k$ be a vertex of $Q'$ (the mutating vertex).
We choose triangles
\[
T_k \arr{u} \bigoplus_{\stackrel{\mbox{\tiny arrows}}{k\to i}} T_i
\to T_k^* \to \Sigma T_k
\]
and
\[
\mbox{ }^* \! T_k \to \bigoplus_{\stackrel{\mbox{\tiny arrows}}{j\to k}} T_j
\arr{v} T_k \to \Sigma ^*\! T_k \ko
\]
where the component of $u$ (respectively $v$) corresponding
to an arrow $\alpha: k \to i$ (respectively $j\to k$) is the
corresponding morphism $T_k \to T_i$ (respectively $T_j \to T_k$).
These triangles are unique up to isomorphism and
called the {\em exchange triangles} associated with
$k$ and $T_1, \ldots, T_n$.

\begin{theorem}[\protect{\cite{BuanMarshReinekeReitenTodorov06}}]
\label{thm:mutation}
\begin{itemize}
\item[a)] The objects $T_k^*$ and $\mbox{}^*\! T_k$ are isomorphic.
\item[b)] The set obtained from $T_1, \ldots, T_n$ by replacing
$T_k$ with $T_k^*$ is cluster-tilting and its associated
cluster is the mutation at $k$ of the cluster associated
with $T_1, \ldots, T_n$.
\item[c)] Two indecomposables $L$ and $M$ appear as the
the pair $(T_k, T_k^*)$ associated with an exchange if
and only if the space $\Ext^1(L,M)$ is one-dimensional.
In this case, the exchange triangles are the unique
(up to isomorphism) non split triangles
\[
L \to E \to M \to \Sigma L \mbox{ and }
M \to E' \to L \to \Sigma M.
\]
\end{itemize}
\end{theorem}

Let us extend the map $L \mapsto X_L$ from
indecomposable to decomposable objects of $\cc_Q$
by requiring that we have
\[
X_N = X_{N_1} X_{N_2}
\]
whenever $N=N_1 \oplus N_2$ (this is compatible
with the muliplicativity of the Caldero-Chapoton map).
We know that if $u_1, \ldots, u_n$ is a cluster and
$B=(b_{ij})$ the associated exchange matrix, then
the mutation at $k$ yields the variable $u_k'$ such
that
\[
u_k u'_k = \prod_{\stackrel{\mbox{\tiny arrows}}{k\to i}} u_i +
\prod_{\stackrel{\mbox{\tiny arrows}}{j\to k}} u_j .
\]
By combining this with the exchange triangles,
we see that in the situation of c), we have
\[
X_L X_M = X_E + X_{E'}.
\]
We would like to generalize this identity to the
case where the space $\Ext^1(L,M)$ is of higher
dimension. For three objects $L$, $M$ and $N$ of
$\cc_Q$, let $Ext^1(L,M)_N$ be the subset of
$\Ext^1(L,M)$ formed by those morphisms $\eps: L\to \Sigma M$
such that in the triangle
\[
M \to E \to L \arr{\eps} \Sigma M \ko
\]
the object $E$ is isomorphic to $N$ (we do not
fix an isomorphism). Notice that this subset is a
cone (\ie stable under multiplication by non zero
scalars) in the vector space $\Ext^1(L,M)$.

\begin{proposition}[\protect{\cite{CalderoKeller08}}]
The subset $\Ext^1(L,M)_N$ is
constructible in $\Ext^1(L,M)$. In particular, it
is a union of algebraic subvarieties. It is
empty for all but finitely isomorphism classes
of objects $N$.
\end{proposition}

If $k$ is the field of complex numbers, we denote
by $\chi$ the Euler characteristic with respect to
singular cohomology with coefficients in a field.
If $k$ is an arbitrary algebraically closed field,
we denote by $\chi$ the Euler characteristic with
respect to \'etale cohomology with proper support.

\begin{theorem}[\protect{\cite{CalderoKeller08}}]
Suppose that $L$ and $M$ are objects of
$\cc_Q$ such that $\Ext^1(L,M)\neq 0$. Then we have
\[
X_L X_M = \sum_{N} \frac{\chi(\mathbb{P}\Ext^1(L,M)_N) +
\chi(\mathbb{P}\Ext^1(M,L)_N)}{\chi(\mathbb{P}\Ext^1(L,M))} \, X_N \ko
\]
where the sum is taken over all isomorphism classes of
objects $N$ of $\cc_Q$.
\end{theorem}

Notice that in the theorem, the objects $L$ and $M$
may be decomposable so that $X_L$ and $X_M$ will not
be cluster variables in general and the $X_N$ do
not form a linearly independent set in the cluster
algebra. Thus, the formula should be considered as
a relation rather than as an alternative definition
for the multiplication of the cluster algebra.
Notice that it nevertheless bears a close resemblance
to the product formula in a dual Hall algebra:
For two objects $L$ and $M$ in a
finitary abelian category of finite global dimension,
we have
\[
[L] * [M] =  \sum_{[N]} \frac{|\Ext^1(L,M)_N|}{|\Ext^1(L,M)|} \, [N] \ko
\]
where the brackets denote isomorphism classes and the vertical bars
the cardinalities of the underlying sets,
\cf Proposition~1.5 of \cite{Schiffmann06}.

\section{Categorification via cluster categories: the acyclic case}
\label{s:categorification-acyclic-case}

\subsection{Categorification}
\label{ss:acyclic-categorificaton}
Let $Q$ be a connected finite quiver without oriented
cycles with vertex set $\{1, \ldots, n\}$. Let $k$
be an algebraically closed field. We have seen in
section~\ref{ss:derived-category} how to define the
bounded derived category $\cd_Q$. We still have
a fully faithful functor from the mesh category
of $\Z Q$ to the category of indecomposables
of $\cd_Q$ but this functor is very far from
being essentially surjective. In fact, its image
does not even contain the injective indecomposable
$kQ$-modules. The methods of the preceding section
therefore do not generalize but most of the results
continue to hold. The derived category $\cd_Q$ still
has a Serre functor (the total left derived functor
of the tensor product functor $?\ten_B D(kQ)$). We can
form the cluster category
\[
\cc_Q = \cd_Q /(S^{-1}\Sigma^2)^\Z
\]
as before and it is still a triangulated category in a
canonical way such that the projection $\pi: \cd_Q \to \cc_Q$ becomes
a triangle functor \cite{Keller05}. Moreover, the decomposition theorem
\ref{thm:decomposition} holds for $\cc_Q$ and each object
$L$ of $\cc_Q$ decomposes into a direct sum
\[
L = \pi(M) \oplus \bigoplus_{i=1}^n \pi(\Sigma P_i)^{m_i}
\]
for some module $M$ and certain multiplicities $m_i$,
$1 \leq i\leq n$, \cf \cite{BuanMarshReinekeReitenTodorov06}. We put
\[
X_L = CC(M) \, \prod_{i=1}^n x_i^{m_i} \ko
\]
where $CC(M)$ is defined as in section~\ref{ss:Caldero-Chapoton-formula}
Notice that in general, $X_L$ can only be expected to be
an element of the fraction field $\Q(x_1, \ldots, x_n)$,
not of the cluster algebra $\ca_Q$ inside this field.
(The exponents in the formula for $X_L$ are perhaps
more transparent in equation~\ref{eq:cluster-char} below).

\begin{theorem} \label{thm:acyclic-main-thm}
Let $Q$ be a finite quiver without oriented cycles with
vertex set $\{1, \ldots, n\}$.
\begin{itemize}
\item[a)] The map $L \mapsto X_L$ induces a bijection from the set
of isomorphism classes of rigid indecomposables of the cluster category $\cc_Q$
onto the set of cluster variables of the cluster algebra $\ca_Q$.
\item[b)] Under this bijection, the clusters correspond exactly to
the cluster-tilting sets, i.e. the sets $T_1, \ldots, T_n$
of rigid indecomposables such that
\[
Ext^1(T_i,T_j)=0
\]
for all $i,j$.
\item[c)] For a cluster-tilting set $T_1$, \ldots, $T_n$, the quiver of the
endomorphism algebra of the sum of the $T_i$ does
not have loops nor $2$-cycles and encodes the exchange
matrix of the \cite{GekhtmanShapiroVainshtein07} seed
containing the corresponding cluster.
\item[d)] If $L$ and $M$ are rigid indecomposables such that the
space $\Ext^1(L,M)$ is one-dimensional, then we have the generalized
exchange relation
\begin{equation} \label{eq:gen-exchange}
X_L X_M = X_B + X_{B'}
\end{equation}
where $B$ and $B'$ are the middle terms of `the' non split triangles
\[
\xymatrix{L \ar[r] & B \ar[r] & M \ar[r] & \Sigma L} \mbox{ and }
\xymatrix{M \ar[r] & B' \ar[r] & L \ar[r] & \Sigma M}.
\]
\end{itemize}
\end{theorem}

Parts a), b) and d) of the theorem are proved in \cite{CalderoKeller06}
and part c) in \cite{BuanMarshReiten08}.
The proofs build on work by many authors
notably
Buan-Marsh-Reiten-Todorov \cite{BuanMarshReitenTodorov07}
Buan-Marsh-Reiten \cite{BuanMarshReiten08},
Buan-Marsh-Reineke-Reiten-Todorov \cite{BuanMarshReinekeReitenTodorov06},
Marsh-Reineke-Zelevinsky \cite{MarshReinekeZelevinsky03},
 \ldots\  and especially on Caldero-Chapoton's explicit
formula for $X_L$ proved in \cite{CalderoChapoton06} for orientations
of simply laced Dynkin diagrams. Another
crucial ingredient of the proof is the Calabi-Yau property
of the cluster category.
An alternative proof of part c) was given by A.~Hubery \cite{Hubery06}
for quivers whose underlying graph is an extended simply
laced Dynkin diagram.

We describe the main steps of the proof of a). The mutation of cluster-tilting
sets is defined using the construction of
section~\ref{ss:mutation-cluster-tilting-sets}.
\begin{itemize}
\item[1)] If $T$ is a cluster-tilting object, then the quiver $Q_T$
of its endomorphism algebra does not have loops or $2$-cycles.
If $T'$ is obtained from $T$ by mutation at the summand $T_1$,
then the quiver $Q_{T'}$ of the endomorphism algebra of $T'$ is the
mutation at the vertex $1$ of the quiver $Q_T$, \cf \cite{BuanMarshReiten08}.
\item[2)] Each rigid indecomposable is contained in a cluster-tilting set.
Any two cluster-tilting sets are linked by a finite sequence
of mutations. This is deduced in \cite{BuanMarshReinekeReitenTodorov06}
from the work of Happel-Unger \cite{HappelUnger05}.
\item[3)] If $(T_1, T_1^*)$ is an exchange pair and
\[
T_1^* \to E \to T_1 \to \Sigma T_1^* \mbox{ and }
T_1 \to E' \to T_1^* \to \Sigma T_1
\]
are the exchange triangles, then we have
\[
X_{T_1} X_{T_1^*} = X_E + X_{E'}.
\]
This is shown in \cite{CalderoKeller06}.
\end{itemize}
It follows from 1)-3) that the map $L \to X_L$ does take rigid
indecomposables to cluster variables and that each cluster variable
is obtained in this way. It remains to be shown that a rigid indecomposable
$L$ is determined up to isomorphism by $X_L$. This follows from
\begin{itemize}
\item[4)] If $M$ is a rigid indecomposable module, the denominator
of $X_M$ is $x_1^{d_1} \ldots x_n^{d_n}$, \cf \cite{CalderoKeller06}.
\end{itemize}
Indeed,  a rigid indecomposable module $M$ is determined,
up to isomorphism, by its dimension vector.

We sum up the relations between the cluster algebra and
the cluster category in the following table
\begin{center}
\begin{tabular}{|c|c|} \hline
cluster algebra & cluster category \\ \hline
multiplication & direct sum \\
addition & ? \\
cluster variables & rigid indecomposables \\
clusters & cluster-tilting sets \\
mutation & mutation \\
exchange relation & exchange triangles \\
$x x^* = m+ m'$ & $T_k   \to M  \to T_k^* \to \Sigma T_k$ \\
                & $T_k^* \to M' \to T_k   \to \Sigma T_k^*$ \\ \hline
\end{tabular}
\end{center}

\subsection{Two applications}
Theorem~\ref{thm:acyclic-main-thm} does shed new light on cluster
algebras. In particular, thanks to the theorem,
Caldero and Reineke \cite{CalderoReineke08}
have made significant progress towards the

\begin{conjecture}
Suppose that $Q$ does not have
oriented cycles. Then all cluster variables of $\ca_Q$ belong to
$\N[x_1^{\pm}, \ldots, x_n^{\pm}]$.
\end{conjecture}

This conjecture is a consequence of a general conjecture of
Fomin-Zelevinsky \cite{FominZelevinsky02}, which here is specialized
to the case of cluster algebras associated with acyclic quivers, for
cluster expansions in the initial cluster. Caldero-Reineke's
work in \cite{CalderoReineke08} is based on
Lusztig's \cite{Lusztig98} and in this sense it does not quite live up
to the hopes that cluster theory ought to explain Lusztig's results.
Notice that in \cite{CalderoReineke08}, the above conjecture
is stated as a theorem. However, a gap in the proof was found
by Nakajima \cite{Nakajima09}: the authors incorrectly identify
their parameter $q$ with Lusztig's parameter $v$, whereas the
correct identification is $v=-\sqrt{q}$.

Here are two applications to the exchange graph of the cluster
algebra associated with an acyclic quiver $Q$:

\begin{corollary}[\cite{CalderoKeller06}]
\begin{itemize}
\item[a)]
For any cluster variable $x$, the set of seeds whose clusters contain $x$
form a connected subgraph of the exchange graph.
\item[b)] The set of seeds whose quiver does not have oriented cycles form a connected subgraph
(possibly empty) of the exchange graph.
\end{itemize}
\end{corollary}

For acyclic cluster algebras, parts a) and b) confirm conjecture~4.14
parts (3) and (4) by Fomin-Zelevinsky in \cite{FominZelevinsky03a}.
By b), the cluster algebra associated with a quiver without oriented
cycles has a well-defined cluster-type.

\subsection{Cluster categories and singularities} The construction
of cluster categories may seem a bit artificial. Nevertheless, cluster categories
do occur `in nature'. In particular, certain triangulated categories
associated with singularities are equivalent to cluster categories.
We illustrate this on the following example: Let the cyclic
group $G$ of order $3$ act on a three-dimensional complex vector space
$V$ by scalar multiplication with a primitive third root of unity.
Let $S$ be the completion at the origin of the coordinate
algebra of $V$ and let $R=S^G$
the fixed point algebra, corresponding to the completion of
the singularity at the origin of the quotient $V/\!/G$.
The algebra $R$ is a Gorenstein ring,
\cf \eg \cite{Watanabe74}, and an
isolated singularity of dimension $3$, \cf \eg
Corollary~8.2 of \cite{IyamaYoshino08}.
The category $\mbox{CM}(R)$ of maximal Cohen-Macaulay modules is
an exact Frobenius category and its stable category
$\underline{\mbox{CM}}(R)$ is a triangulated category. By Auslander's results
\cite{Auslander76}, \cf Lemma 3.10 of \cite{Yoshino90}, it
is $2$-Calabi Yau. One can show that it is equivalent
to the cluster category $\cc_Q$ for the quiver
\[
Q : \xymatrix{ 1 \ar[r] \ar@<1ex>[r] \ar@<-1ex>[r] & 2}
\]
by an equivalence which takes the cluster-tilting
object $T=kQ$ to $S$ considered as an $R$-module.
This example can be found in \cite{KellerReiten06}, where
it is deduced from an abstract characterization of cluster categories.
A number of similar examples can be found in \cite{BurbanIyamaKellerReiten08}
and \cite{KellerVandenBergh08}.

\section{Categorification via 2-Calabi-Yau categories}
\label{s:categorification-2-Calabi-Yau}

The extension of the results of the preceding sections to quivers
containing oriented cycles is the subject of ongoing research,
\cf for example \cite{DerksenWeymanZelevinsky08} \cite{GeissLeclercSchroeer07b}
\cite{BuanIyamaReitenSmith08}.
Here we present an approach based on the fact that many arguments
developed for cluster categories apply more generally to suitable
triangulated categories, whose most important property it is to be
Calabi-Yau of dimension $2$. As an application, we will sketch a proof
of the periodicity conjecture in section~\ref{s:periodicity} (details
will appear elsewhere \cite{Keller10a}).

\subsection{Definition and main examples}
\label{ss:CY-definition} Let $k$ be an
algebraically closed field and let $\cc$ a triangulated
category with suspension functor $\Sigma$ where all idempotents
split (\ie each idempotent endomorphism $e$ of an object
$M$ is the projection onto $M_1$ along $M_2$ in a decomposition $M=M_1 \oplus M_2$).
We assume that
\begin{itemize}
\item[1)] $\cc$ is $\Hom$-finite (\ie we have $\dim \cc(L,M)<\infty$ for all
$L$, $M$ in $\cc$) and the decomposition theorem~\ref{thm:decomposition} holds for $\cc$;
\item[2)] $\cc$ is $2$-Calabi-Yau, \ie we are given bifunctorial isomorphisms
\[
D\cc(L,M) \iso \cc(M, \Sigma^2 L) \ko L,M \in \cc ;
\]
\item[3)] $\cc$ admits a cluster-tilting object $T$, \ie
\begin{itemize}
\item[a)] $T$ is the sum of pairwise non-isomorphic indecomposables,
\item[b)] $T$ is rigid and
\item[c)] for each object $L$ of $\cc$, if $\Ext^1(T,L)$ vanishes, then
$L$ belongs to the category $\add(T)$ of direct factors of finite direct
sums of copies of $T$.
\end{itemize}
\end{itemize}

If all these assumptions hold, we say that {\em $(\cc,T)$ is
a $2$-Calabi-Yau category with cluster tilting object}.
If $Q$ is a finite quiver, we say that a $2$-Calabi-Yau category
$\cc$ with cluster-tilting object $T$
is a {\em $2$-Calabi-Yau realization of $Q$}
if $Q$ is the quiver of the endomorphism algebra of $T$.

For example, if $\cc$ is the cluster category of a finite quiver
$Q$ without oriented cycles, conditions 1) and 2) hold and an
object $T$ is cluster-tilting in the above sense iff it is the
direct sum of a cluster-tilting set $T_1$, \ldots, $T_n$, where
$n$ is the number of vertices of $Q$, \cf \cite{BuanMarshReinekeReitenTodorov06}.
The `initial' cluster-tilting
object in this case is $T=kQ$ (the image in $\cc_Q$ of the free
module of rank one) and $(\cc, kQ)$ is the {\em canonical
$2$-Calabi-Yau realization} of the
quiver $Q$ (without oriented cycles).

The second main class of examples comes from the work
of Geiss-Leclerc-Schr\"oer: Let $\Delta$ be a simply laced
Dynkin diagram, $\vec{\Delta}$ a quiver with underlying graph
$\Delta$ and $\overline{\Delta}$ the {\em doubled quiver} obtained from
$\vec{\Delta}$ by adjoining an arrow $\alpha^* : j\to i$
for each arrow $\alpha: i\to j$. The {\em preprojective algebra}
$\Lambda=\Lambda(\Delta)$ is the quotient of the path algebra
of $\overline{\Delta}$ by the ideal generated by the relator
\[
\sum_{\alpha \in Q_1} \alpha \alpha^* - \alpha^* \alpha.
\]
For example, if $\vec{\Delta}$ is the quiver
\[
\xymatrix{ 1 \ar[r]^\alpha & 2 \ar[r]^\beta & 3 } \ko
\]
then $\overline{\Delta}$ is the quiver
\[
\xymatrix{ 1 \ar@<0.5ex>[r]^\alpha & 2 \ar@<0.5ex>[r]^\beta
\ar@<0.5ex>[l]^{\alpha^*} & 3 \ar@<0.5ex>[l]^{\beta^*}}
\]
and the ideal generated by the above sum of commutators
is also generated by the elements
\[
\alpha^* \alpha, \alpha\alpha^* - \beta^* \beta, \beta\beta^*.
\]

It is classical, \cf \eg \cite{Ringel98}, that the algebra $\Lambda=\Lambda(\Delta)$ is
a finite-dimensional (!) selfinjective algebra (\ie $\Lambda$ is
injective as a right $\Lambda$-module over itself). Let $\mod \Lambda$
denote the category of $k$-finite-dimensional right $\Lambda$-modules.
The {\em stable module category $\ul{\mod} \Lambda$} is the quotient of $\mod\Lambda$ by
the ideal of all morphisms factoring through a projective
module. This category carries a canonical triangulated
structure (like any stable module category of a self-injective
algebra): The suspension is constructed by choosing exact
sequences of modules
\[
0 \to L \to IL \to \Sigma L \to 0
\]
where $IL$ is injective (but not necessarily functorial
in $L$; the object $\Sigma L$ becomes functorial in
$L$ when we pass to the stable category). The triangles
are by definition isomorphic to standard triangles obtained
from exact sequences of modules as follows: Let
\[
0 \to L \arr{i} M \arr{p} N \to 0
\]
be a short exact sequence of $\mod\Lambda$. Choose
a commutative diagram
\[
\xymatrix{
0 \ar[r] & L \ar[d]_{\id} \ar[r]^i & M \ar[r]^p \ar@{.>}[d] & N \ar[r] \ar@{.>}[d]^e & 0 \\
0 \ar[r] & L \ar[r] & IL \ar[r] & \Sigma L \ar[r] & 0
}
\]
Then the image of $(i,p,e)$ is a standard triangle in the
stable module category. As shown in \cite{CrawleyBoevey00},
the stable module category $\cc=\ul{\mod}\Lambda$ is
$2$-Calabi-Yau and it is easy to check that assumption 1)
holds.

\begin{theorem}[Geiss-Leclerc-Schr\"oer] The category
$\cc=\ul{\mod} \Lambda$ admits a cluster-tilting object
$T$ such that the quiver of $\End(T)$ is obtained
from that of the category of indecomposable
$k\Delta$-modules by deleting the injective vertices
and adding an arrow $v \to \tau v$ for each non
projective vertex $v$.
\end{theorem}

Here, by the quiver of the category of indecomposable
$k\Delta$-modules, we mean the full subquiver of the
repetition $\Z Q$ which is formed by the vertices corresponding
to modules (complexes concentrated in degree $0$). Thus
for $\Delta=\vec{A}_5$, this quiver is as follows:
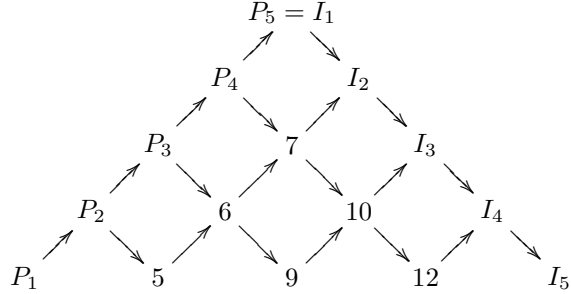
\[
\begin{xy} 0;<1pt,0pt>:<0pt,-1pt>::
(0,100) *+{P_1} ="0",
(25,75) *+{P_2} ="1",
(50,50) *+{P_3} ="2",
(75,25) *+{P_4} ="3",
(100,0) *+{P_5=I_1} ="4",
(50,100) *+{5} ="5",
(75,75) *+{6} ="6",
(100,50) *+{7} ="7",
(125,25) *+{I_2} ="8",
(100,100) *+{9} ="9",
(125,75) *+{10} ="10",
(150,50) *+{I_3} ="11",
(150,100) *+{12} ="12",
(175,75) *+{I_4} ="13",
(200,100) *+{I_5} ="14",
"0", {\ar"1"},
"1", {\ar"2"},
"1", {\ar"5"},
"2", {\ar"3"},
"2", {\ar"6"},
"3", {\ar"4"},
"3", {\ar"7"},
"4", {\ar"8"},
"5", {\ar"6"},
"6", {\ar"7"},
"6", {\ar"9"},
"7", {\ar"8"},
"7", {\ar"10"},
"8", {\ar"11"},
"9", {\ar"10"},
"10", {\ar"11"},
"10", {\ar"12"},
"11", {\ar"13"},
"12", {\ar"13"},
"13", {\ar"14"},
\end{xy}
\]
where we have marked the indecomposable projectives $P_i$ and the indecomposable injectives
$I_j$. If we remove the vertices corresponding to the indecomposable injectives
(and all the arrows incident with them) and add an arrow $v\to \tau v$ for each
vertex not corresponding to an indecomposable projective, we obtain the
following quiver
\[
\begin{xy} 0;<1pt,0pt>:<0pt,-1pt>::
(75,0) *+{0} ="0",
(50,25) *+{1} ="1",
(100,25) *+{2} ="2",
(25,50) *+{3} ="3",
(75,50) *+{4} ="4",
(125,50) *+{5} ="5",
(0,75) *+{6} ="6",
(50,75) *+{7} ="7",
(100,75) *+{8} ="8",
(150,75) *+{9} ="9",
"1", {\ar"0"},
"0", {\ar"2"},
"2", {\ar"1"},
"3", {\ar"1"},
"1", {\ar"4"},
"4", {\ar"2"},
"2", {\ar"5"},
"4", {\ar"3"},
"6", {\ar"3"},
"3", {\ar"7"},
"5", {\ar"4"},
"7", {\ar"4"},
"4", {\ar"8"},
"8", {\ar"5"},
"5", {\ar"9"},
"7", {\ar"6"},
"8", {\ar"7"},
"9", {\ar"8"},
\end{xy}
\]
In a series of
papers \cite{GeissLeclercSchroeer05} \cite{GeissLeclercSchroeer07d} \cite{GeissLeclercSchroeer06}
\cite{GeissLeclercSchroeer08b} \cite{GeissLeclercSchroeer07b}
\cite{GeissLeclercSchroeer07a}
\cite{GeissLeclercSchroeer08a}, Geiss-Leclerc-Schr\"oer have obtained
remarkable results for a
class of quivers which are important in the study of (dual
semi-)canonical bases. They use an analogue \cite{GeissLeclercSchroeer07a}
of the Caldero-Chapoton map due ultimately to Lusztig \cite{Lusztig00}.
The class they consider has been further
enlarged by Buan-Iyama-Reiten-Scott \cite{BuanIyamaReitenScott09}.
Thanks to their results, an analogue
of Caldero-Chapoton's formula and a weakened version of theorem~\ref{thm:acyclic-main-thm}
was proved in \cite{FuKeller10} for an even larger class.

\subsection{Calabi-Yau reduction} \label{ss:Calabi-Yau-reduction}
Suppose that $(\cc,T)$ is a $2$-Calabi-Yau
realization of a quiver $Q$ so that for $1 \leq i\leq n$,
the vertex $i$ of $Q$ corresponds to the indecomposable summand
$T_i$ of $T$. Let $J$ be a subset of the set of
vertices of $Q$ and let $Q'$ be the quiver obtained from $Q$ by deleting
all vertices in $J$ and all arrows incident with one of these vertices.
Let $\cu$ be the full subcategory of $\cc$ formed by the objects $U$
such that $\Ext^1(T_j, U)=0$ for all $j \in J$. Note that all $T_i$,
$1\leq i\leq n$, belong to $\cu$. Let $\langle T_j | j\in J\rangle$
denote the ideal of $\cu$ generated by the identities of the objects
$T_j$, $j\in J$. By imitating the construction of the triangulated
structure on a stable category, one can endow the quotient
\[
\cc' = \cu/\langle T_j | j\in J\rangle
\]
with a canonical structure of triangulated category, \cf \cite{IyamaYoshino08}.

\begin{theorem}[Iyama-Yoshino \protect{\cite{IyamaYoshino08}}]
The pair $(\cc', T)$ is a $2$-Calabi-Yau
realization of the quiver $Q'$. Moreover, the projection $\cu\to \cc'$
induces a bijection between the cluster-tilting sets of $\cc$
containing the $T_j$, $j\in J$, and the cluster-tilting sets of
$\cc'$.
\end{theorem}

\subsection{Mutation} Let $(\cc,T)$ be a $2$-Calabi-Yau category
with cluster-tilting object. Let $T_1$ be an indecomposable direct
factor of $T$.

\begin{theorem}[Iyama-Yoshino \protect{\cite{IyamaYoshino08}}]
\label{thm:2-CY-mutation} Up to isomorphism, there is a unique
indecomposable object $T_1^*$ not isomorphic to $T_1$ such that the
object $\mu_1(T)$ obtained from $T$ by replacing the indecomposable
summand $T_1$ with $T_1^*$ is cluster-tilting.
\end{theorem}

We call $\mu_1(T)$ the {\em mutation} of $T$ at $T_1$. If $\cc$ is the
cluster category of a finite quiver without oriented cycles, this
operation specializes of course to the one defined in
section~\ref{ss:mutation-cluster-tilting-sets}.
However, in general, the quiver $Q$ of the endomorphism algebra of
$T$ may contain loops and $2$-cycles and then the quiver of
the endomorphism algebra of $\mu_1(T)$ is not determined by
$Q$. Let us illustrate this phenomenon
on the following example (taken from Proposition~2.6 of
\cite{BurbanIyamaKellerReiten08}): Let $\cc$ be the orbit category of the bounded
derived category $\cd_{\vec{D}_6}$ under the action of the autoequivalence
$\tau^2$. Then $\cc$ satisfies the assumptions 1) and 2) of section~\ref{ss:CY-definition}.
Its category of indecomposables is equivalent to the mesh category of
the quiver
\[
\begin{xy} 0;<1pt,0pt>:<0pt,-1pt>::
(0,100) *+{A} ="0",
(0,50) *+{1} ="1",
(0,25) *+{B} ="2",
(0,0) *+{C} ="3",
(25,75) *+{4} ="4",
(25,25) *+{5} ="5",
(50,100) *+{A'} ="6",
(50,50) *+{7} ="7",
(50,25) *+{B'} ="8",
(50,0) *+{C'} ="9",
(75,75) *+{10} ="10",
(75,25) *+{11} ="11",
(100,100) *+{A} ="12",
(100,50) *+{13} ="13",
(100,25) *+{B} ="14",
(100,0) *+{C} ="15",
(125,75) *+{16} ="16",
(125,25) *+{17} ="17",
"0", {\ar"4"},
"1", {\ar"4"},
"1", {\ar"5"},
"2", {\ar"5"},
"3", {\ar"5"},
"4", {\ar"6"},
"4", {\ar"7"},
"5", {\ar"7"},
"5", {\ar"8"},
"5", {\ar"9"},
"6", {\ar"10"},
"7", {\ar"10"},
"7", {\ar"11"},
"8", {\ar"11"},
"9", {\ar"11"},
"10", {\ar"12"},
"10", {\ar"13"},
"11", {\ar"13"},
"11", {\ar"14"},
"11", {\ar"15"},
"12", {\ar"16"},
"13", {\ar"16"},
"13", {\ar"17"},
"14", {\ar"17"},
"15", {\ar"17"},
\end{xy}
\]
where the vertices labeled $A$, $1$, $B$, $C$, $4$, $5$ on the left have to be
identified with the vertices labeled $A$, $13$, $B$, $C$, $16$, $17$ on the right.
In this case, there are exactly $6$ indecomposable rigid objects, namely
$A$, $B$, $C$, $A'$, $B'$ and $C'$. There are exactly $6$ cluster-tilting
sets. The following is the exchange graph: Its vertices are cluster-tilting
sets (we write $AC$ instead of $\{A,C\}$) and its edges represent mutations.
\[
\begin{xy} 0;<0.7pt,0pt>:<0pt,-0.7pt>::
(36,0) *+{AC} ="0",
(86,0) *+{AB} ="1",
(123,51) *+{C'B} ="2",
(86,100) *+{C'A'} ="3",
(36,100) *+{B'A'} ="4",
(0,50) *+{B'C} ="5",
"0", {\ar@{-}"1"},
"5", {\ar@{-}"0"},
"1", {\ar@{-}"2"},
"2", {\ar@{-}"3"},
"3", {\ar@{-}"4"},
"4", {\ar@{-}"5"},
\end{xy}
\]
The quivers of the endomorphism algebras are as follows:
\begin{align*}
AC, AB  &: \quad\quad \xymatrix@C=0.5cm{\circ \ar@<0.5ex>[r] & \circ \ar@<0.5ex>[l] \ar@(ur,dr)[]} \\
B'C, C'B &: \quad\quad \xymatrix@C=0.5cm{\circ \ar@<0.5ex>[r] & \circ \ar@<0.5ex>[l]} \\
B'A', C'A' &: \quad\quad \xymatrix@C=0.5cm{\circ \ar@(ul,dl)[]\ar@<0.5ex>[r] & \circ \ar@<0.5ex>[l]}
\end{align*}

In the setting of the above theorem, there are still exchange
triangles as in theorem~\ref{thm:mutation} but their description is
different: Let $\ct'$ be the full subcategory of $\cc$ formed by
the direct sums of indecomposables $T_i$, where $i$ is different
from $k$. A {\em left $\ct'$-approximation of $T_k$} is a morphism
$f: T_k \to T'$ with $T'$ in $\ct'$ such that any morphism
from $T_k$ to an object of $\ct'$ factors through $f$. A left
$\ct'$-approximation $f$ is {\em minimal} if for each endomorphism
$g$ of $T'$, the equality $gf=f$ implies that $g$ is invertible.
Dually, one defines (minimal) right $\ct'$-approximations. It
is not hard to show that they always exist.

\begin{theorem}[Iyama-Yoshino \protect{\cite{IyamaYoshino08}}]
If $(T_1, T_1^*)$ is an exchange
pair, there are non split triangles, unique up to isomorphism,
\[
T_1 \arr{f} E' \to T_1^* \to \Sigma T_1 \quad\mbox{ and }\quad
T_1^* \to E \arr{g} T_1 \to \Sigma T_1^*
\]
such that $f$ is a minimal left $\ct'$-approximation and
$g$ a minimal right $\ct'$-approxi\-mation.
\end{theorem}

\subsection{Simple mutations, reachable cluster-tilting objects} Let $(\cc, T)$ be a $2$-CY category
with cluster-tilting object and $T_1$ an indecomposable
direct summand of $T$. Let $(T_1, T_1^*)$ be the
corresponding exchange pair and
\[
T_1^* \to E \to T_1 \to \Sigma T_1^*
\]
the exchange triangle. The long exact sequence induced in
$\cc(T_1,?)$ by this triangle yields a short exact sequence
\[
\cc_\ct'(T_1, T_1) \to \cc(T_1, T_1) \to \Ext^1(T_1, T_1^*) \to 0\ko
\]
where the leftmost term is the space of those endomorphisms of
$T_1$ which factor through a sum of copies of $T/T_1$. Now the
algebra $\cc(T_1, T_1)$ is local and its residue field is $k$
(since $k$ is algebraically closed). We deduce the following
lemma.

\begin{lemma} \label{lemma:simple-mutation}
The quiver of the endomorphism algebra of $T$ does not have
a loop at the vertex corresponding to $T_1$ iff we have
$\dim \Ext^1(T_1, T_1^*)=1$ iff
$\Ext^1(T_1, T_1^*)$ is a simple module over $\cc(T_1, T_1)$.
In this case, in the exchange triangles
\[
T_1^* \to E \to T_1 \to \Sigma T_1^* \mbox{ and }
T_1 \to E' \to T_1^* \to \Sigma T_1 \ko
\]
we have
\[
E= \bigoplus_{\stackrel{\mbox{\tiny arrows}}{i\to 1}} T_i
\mbox{ and }
E'=\bigoplus_{\stackrel{\mbox{\tiny arrows}}{1 \to j}} T_j .
\]
\end{lemma}

We say that the mutation at $T_1$ is {\em simple} if
the conditions of the lemma hold. If $\cc$ is the cluster
category of a finite quiver without oriented cycles, all
mutations in $\cc$ are simple, by part c) of theorem~\ref{thm:acyclic-main-thm}.
By a theorem of Geiss-Leclerc-Schr\"oer, if $\cc$ is the
stable module category of a preprojective algebra of
Dynkin type, the quiver of any cluster-tilting object
in $\cc$ does not have loops nor $2$-cycles. So again,
all mutations are simple in this case.

\begin{theorem}[Buan-Iyama-Reiten-Scott \protect{\cite{BuanIyamaReitenScott09}}]
\label{thm:quiver-2-CY-mutation}
Suppose that the quivers
$Q$ and $Q'$ of the endomorphism algebras of $T$ and $T'=\mu_1(T)$
do not have loops nor $2$-cycles. Then $Q'$ is the mutation of
$Q$ at the vertex $1$.
\end{theorem}

We define a cluster-tilting object $T'$ to be {\em reachable
from $T$} if there is a sequence of mutations
\[
T=T^{(0)} \leadsto T^{(1)} \leadsto \ldots \leadsto T^{(N)}=T'
\]
such that the quiver of $\End(T^{(i)})$ does not have loops
nor $2$-cycles for all $1\leq i\leq N$. We define a rigid
indecomposable of $\cc$ to be {\em reachable from $T$} if
it is a direct summand of a reachable cluster-tilting object.

\begin{corollary} If a cluster-tilting object $T'$ is reachable from $T$, then
the quiver of the endomorphism algebra of $T'$ is mutation-equivalent
to the quiver of the endomorphism algebra of $T$.
If $\cc$ is a cluster-category or the stable module category
of the preprojective algebra of a Dynkin diagram, then
all quivers mutation-equivalent to $Q$ are obtained
in this way.
\end{corollary}

\subsection{Combinatorial invariants} Let $(\cc,T)$ be
a $2$-Calabi-Yau category with cluster-tilting object.
Let $\ct$ be the full subcategory whose objects are
all direct factors of finite direct sums of copies
of $T$. Notice that $\ct$ is equivalent to the
category of finitely generated projective modules over
the endomorphism algebra of $T$. Let $K_0(\ct)$ be the Grothendieck group
of the additive category $\ct$. Thus, the group
$K_0(\ct)$ is free abelian on the isomorphism
classes of the indecomposable summands of $T$.

\begin{lemma}[Keller-Reiten \protect{\cite{KellerReiten07}}]
\label{lemma:Keller-Reiten}
For each object $L$ of $\cc$, there is
a triangle
\[
T_1 \to T_0 \to L \to \Sigma T_1
\]
such that $T_0$ and $T_1$ belong to $\ct$. The difference
\[
[T_0] - [T_1]
\]
considered as an element of $K_0(\ct)$ does not depend
on the choice of this triangle.
\end{lemma}
In the situation of the lemma,
we define the {\em index $\ind(L)$ of $L$} as
the element $[T_0]-[T_1]$ of $K_0(\ct)$.

\begin{theorem}[Dehy-Keller \protect{\cite{DehyKeller08}}] \label{thm:Dehy-Keller}
\begin{itemize}
\item[a)] Two rigid objects are isomorphic iff their indices
are equal.
\item[b)] The indices of the indecomposable summands of a cluster-tilting
object form a basis of $K_0(\ct)$. In particular, all cluster-tilting
objects have the same number of pairwise non isomorphic indecomposable
summands.
\end{itemize}
\end{theorem}

Let $B$ be the endomorphism algebra of $T$. For two finite-dimensional
right $B$-modules $L$ and $M$ put
\[
\langle L, M \rangle_a = \dim \Hom(L,M) -\dim \Ext^1(L,M)
                        -\dim \Hom(M,L) +\dim \Ext^1(M,L).
\]
This is the antisymmetrization of a truncated Euler form. A priori
it is defined on the split Grothendieck group of the category
$\mod B$ (\ie the quotient of the free abelian group on the
isomorphism classes divided by the subgroup generated by
all relations obtained from direct sums in $\mod B$).

\begin{proposition}[Palu \protect{\cite{Palu08a}}] The form $\langle, \rangle_a$ descends
to an antisymmetric form on $K_0(\mod B)$. Its matrix in the basis
of the simples is the antisymmetric matrix associated with
the quiver of $B$ (loops and $2$-cycles do not contribute
to this matrix).
\end{proposition}

Let $T_1, \ldots, T_n$ be the pairwise non isomorphic indecomposable
direct summands of $T$. For $L\in \cc$, we define the integer $g_i(L)$ to be the
multiplicity of $[T_i]$ in the index $\ind(L)$, $1\leq i\leq n$,
and we define the element $X'_L$ of the field $\Q(x_1, \ldots, x_n)$ by
\begin{equation} \label{eq:cluster-char}
X'_L= \prod_{i=1}^n x_i^{g_i(L)} \sum_e \chi(\Gr_e(\Ext^1(T,L))) \prod_{i=1}^n x_i^{\langle S_i, e\rangle_a} \ko
\end{equation}
where $S_i$ is the simple quotient of the indecomposable
projective $B$-module $P_i=\Hom(T,T_i)$. Notice that we have
$X'_{T_i}=x_i$, $1\leq i\leq n$.
If $\cc$ is the cluster-category of a finite quiver $Q$ without
oriented cycles and $T=kQ$, then we have $X'_L=X_{\Sigma L}$
in the notations of section~\ref{s:categorification-acyclic-case} and
the formula for $X'_L$ is essentially another expression for
the Caldero-Chapoton formula.

Now let $Q$ be the quiver of the endomorphism algebra of $T$
in $\cc$ and let $\ca_Q$ be the associated cluster algebra.

\begin{theorem}[Palu \protect{\cite{Palu08a}}] \label{thm:Palu-formula}
If $L$ and $M$ are objects of $\cc$ such that
$\Ext^1(L,M)$ is one-dimensional, then we have
\[
X'_L X'_M = X'_E + X'_{E'} \ko
\]
where
\[
L \to E \to M \to \Sigma L \mbox{ and } M \to E' \to L \to \Sigma M
\]
are `the' two non split triangles. Thus, if $L$ is a rigid indecomposable
reachable from $T$, then $X'_L$ is a cluster variable of $\ca_Q$.
\end{theorem}

\begin{corollary} Suppose that $\cc$ is the cluster category $\cc_Q$ of
a finite quiver $Q$ without oriented cycles. Let
$
\ca_Q \subset \Q(x_1, \ldots, x_n)
$
be the associated cluster algebra and $x\in \ca_Q$ a cluster
variable. Let $L\in \cc_Q$ be the unique (up to isomorphism)
indecomposable rigid object such that $x=X_{\Sigma L}$.
Let $u_1, \ldots, u_n$ be an arbitrary cluster of $\ca_Q$
and $T_1$, \ldots, $T_n$ the cluster-tilting set such that
$u_i=X_{T_i}$, $1\leq i\leq n$. Then the expression of $x$
as a Laurent polynomial in $u_1, \ldots, u_n$ is given by
\[
x = X'_{L}(u_1, \ldots, u_n).
\]
\end{corollary}

The expression for $X'_L$ makes it natural to define
the polynomial $F'_L\in\Z[y_1, \ldots y_n]$ by
\begin{equation}
\label{eq:F-polynomial}
F'_L=\sum_e \chi(\Gr_e(\Ext^1(T,L))) \prod_{j=1}^n y_j^{e_j}.
\end{equation}
We then have
\[
X'_L= \prod_{i=1}^n x_i^{g_i(L)}
F'_L(\prod_{i=1}^n x_i^{b_{i1}}, \ldots, \prod_{i=1}^n x_i^{b_{in}}).
\]
The polynomial $F'_L$ is related to Fomin-Zelevinsky's
$F$-polynomials \cite{FominZelevinsky03b} \cite{FominZelevinsky07}, as we will see below.

\subsection{More mutants categorified}
Let $Q$ be a finite quiver without loops nor $2$-cycles
with vertex set $\{1, \ldots, n\}$.
Let $\mathbb{T}_n$ be the regular $n$-ary tree: Its edges
are labeled by the integers $1$, \ldots, $n$ such that
the $n$ edges emanating from each vertex carry different
labels. Let $t_0$ be a vertex of $\mathbf{T}_n$.
To each vertex $t$ of $\mathbb{T}_n$ we associate
a seed $(Q_t, \mathbf{x}_t)$ (\cf section~\ref{ss:seeds})
such that at $t=t_0$, we have $Q_t=Q$ and $\mathbf{x}_t=\{x_1, \ldots, x_n\}$
and whenever $t$ is linked to $t'$ by an edge labeled $i$, we
have $(Q_{t'}, \mathbf{x}_{t'})=\mu_i (Q_t, \mathbf{x}_t)$.

Now assume that $Q$ admits a $2$-Calabi-Yau realization $(\cc,T)$.
According to the mutation theorem~\ref{thm:quiver-2-CY-mutation}, we
can associate a cluster-tilting object $T_t$ to each vertex
$t$ of $\mathbb{T}_n$ such that at $t=t_0$, we have $T_t=T$
and that whenever $t$ is linked to $t'$ by an edge
labeled $i$, we have $T_{t'}=\mu_i(T_t)$.

Now let
\[
\xymatrix{t_0 \ar@{-}[r]^{i_1} & t_1 \ar@{-}[r]^{i_2} & \ldots & t_N \ar@{-}[l]_{i_N} }
\]
be a path in $\mathbb{T}_n$ and suppose that for each
$1\leq i\leq N$, the quiver of the endomorphism algebra of
$T_{t_i}$ does not have loops nor $2$-cycles. Then it follows
by induction from theorem~\ref{thm:quiver-2-CY-mutation} that the quiver
of the endomorphism algebra of $T_{t_i}$ is $Q_{t_i}$, $1\leq i\leq N$,
and from theorem~\ref{thm:Palu-formula} that the cluster
$\mathbf{x}_{t_i}$ equals the image under $L \mapsto X'_L$ of
the set of indecomposable direct factors of $T_{t_i}$.

Following \cite{FominZelevinsky07}, let us consider three other
pieces of data associated with each vertex $t$ of $\mathbb{T}_n$:
\begin{itemize}
\item the tropical $Y$-variables $y_{1,t}, \ldots, y_{n,t}$,
\item the $F$-polynomials $F_{1,t}, \ldots, F_{n,t}$,
\item the (non tropical) $Y$-variables $Y_{1,t}, \ldots, Y_{n,t}$.
\end{itemize}
Here the tropical $Y$-variables are monomials in the indeterminates
$y_1, \ldots, y_n$ and their inverses.
At $t=t_0$, we have $y_{i,t}=y_i$, $1\leq i\leq n$. If $t$ is
linked to $t'$ by an edge labelled $i$, then
\[
y_{j,t'}=\left\{ \begin{array}{ll} y_{i,t}^{-1} & \mbox{ if } j=i\\
y_{j,t} y_{i,t}^{[b_{ij}]_+}(y_{i,t}\oplus 1)^{-b_{ij}} & \mbox{ if } j\neq i. \end{array} \right.
\]
Here, $(b_{ij})$ is the antisymmetric matrix associated with $Q_t$,
for an integer $a$, we write $[a]_+$ for $\max(a,0)$, and, for a monomial $m$
which is the product of powers $y_j^{e_j}$, $1\leq j\leq n$, we write $m\oplus 1$ for the product
of the factors $y_j^{e_j}$ with
negative exponents $e_j$. Notice that
$[b_{ij}]_+$ is the number of arrows from $i$ to $j$ in $Q_t$.

The $F$-polynomials lie in $\Z[y_1, \ldots, y_n]$. At $t=t_0$,
they all equal $1$. If $t$ is linked to $t'$ by an edge labeled $i$,
then
\begin{align} \label{eq:rec-F-polynomials}
F_{j,t'} &= F_{j,t} \quad \mbox{if } j\neq i \ko \\
F_{i,t'} &= \frac{1}{F_{i,t}}
( \prod_{c_{li}>0} y_l^{c_{li}} \prod_{j=1}^n F_{j,t}^{[b_{ij}]_+} +
\prod_{c_{li}<0} y_l^{-c_{li}} \prod_{j=1}^n F_{j,t}^{[-b_{ij}]_+}) \ko
\end{align}
where the $c_{li}$ are the exponents in the tropical $Y$-variables
$y_{i,t}=\prod_{l=1}^n y_l^{c_{li}}$ and $B=(b_{ij})$ is the antisymmetric
matrix associated with $Q_t$.

Finally, the non tropical $Y$-variables $Y_{j,t}$ lie in the field
$\Q(y_1, \ldots, y_n)$. At $t=t_0$, we have $Y_{j,t}=y_j$, $1\leq j\leq n$, and if $t$ and
$t'$ are linked by an edge labeled $i$, then
\begin{equation} \label{eq:Y-variable-mutation}
Y_{j,t'}=\left\{
\begin{array}{cc} Y_{i,t}^{-1}              & \mbox{ if } j= i \\
Y_{j,t} Y_{i,t}^{[b_{ij}]_+} (Y_{i,t}+1)^{-b_{ij}} & \mbox{ if } j\neq i.
\end{array}
\right.
\end{equation}

Now let
\[
\xymatrix{t_0 \ar@{-}[r]^{i_1} & t_1 \ar@{-}[r]^{i_2} & \ldots & t_N \ar@{-}[l]_{i_N} }
\]
be a path in $\mathbb{T}_n$ and suppose that for each
$1\leq i\leq N$, the quiver of the endomorphism algebra of
$T_{t_i}$ does not have loops nor $2$-cycles. Let $t=t_N$
and let $T'_1, \ldots, T'_n$ be the indecomposable summands
of $T'=T_t$. Let $\ct'$ be the full subcategory of $\cc$ formed
by the direct summands of finite direct sums of copies of $T'$.
Notice that $\ct'^{op}$ is a cluster-tilting subcategory of
$\cc^{op}$ so that each object $X$ of $\cc$ also has a well-defined
index in $\cc^{op}$ with respect to $\ct'^{op}$; we denote it by $\ind^{op}_{\ct'}(X)$.
If we identify the Grothendieck groups of $\ct'$ and $\ct'^{op}$,
this index equals $-\ind_{\ct'}(\Sigma X)$.
\begin{theorem} \label{thm:categorify-mutants}
\begin{itemize}
\item[a)] The exchange matrix $B_t=(b_{ij})$ associated with $t$ is the
antisymmetric matrix associated with the quiver of the endomorphism
algebra of $T_t$.
\item[b)] We have $y_{l,t}=\prod_{j=1}^n y_j^{c_{lj}}$, $1\leq l\leq n$,
where $c_{lj}$ is defined by
\[
\ind^{op}_{\ct'}(T_j)=\sum_{l=1}^n c_{lj} [T'_l].
\]
\item[c)] We have $F_{j,t}=F'_{T'_j}$, $1\leq j\leq n$, where $F'_{T'_j}$ is
defined by equation~\ref{eq:F-polynomial}.
\item[d)] We have
\[
Y_{j,t}= y_{j,t} \prod_{i=1}^n F_{i,t}^{b_{ij}}.
\]
\end{itemize}
\end{theorem}
Here part a) follows by induction from theorem~\ref{thm:quiver-2-CY-mutation}.
Part b) is proved by induction using theorem~3.1 of \cite{DehyKeller08}.
Part c) is proved from the recursive definition~\ref{eq:rec-F-polynomials}
by a multiplication formula extracted from \cite{Palu08a}.
Part d) is Proposition~3.12 of \cite{FominZelevinsky07}.

\subsection{$2$-CY categories from algebras of global dimension $2$}
\label{ss:2-CY-from-global-dim-2}
The results of the preceding sections only become interesting if we are
able to construct $2$-Calabi-Yau realizations for large classes of
quivers. In the case of a quiver without oriented cycles, this problem
is solved by the cluster category. The results of Geiss-Leclerc-Schr\"oer
provide another large class of quivers admitting $2$-Calabi-Yau realizations.
Here, we will exhibit a construction which generalizes both cluster categories
of quivers without oriented cycles and the stable categories of preprojective
algebras of Dynkin diagrams.

Let $k$ be an algebraically closed field and $A$ a finite-dimensional
$k$-algebra of global dimension $\leq 2$. For example, $A$ can be
the path algebra of a finite quiver without oriented cycles. Let $\cd_A$
be the bounded derived category of the category $\mod A$ of
$k$-finite-dimensional right $A$-modules. It admits a Serre functor,
namely the total derived functor of the tensor product $?\ten_A DA$
with the $k$-dual bimodule of $A$ considered as a bimodule over itself.
We can form the orbit category
\[
\cd_A/(S^{-1}\Sigma^2)^\Z.
\]
This is a $k$-linear category endowed with a suspension functor
(induced by $\Sigma$) but in general it is no longer triangulated.
Nevertheless, one can show that it embeds fully faithfully into a
`smallest triangulated overcategory' \cite{Keller05}. We denote this overcategory
by $\cc_A$ and call it the {\em generalized cluster category of $A$}.

\begin{theorem}[Amiot  \protect{\cite{Amiot09}}] If the functor
\[
\Tor_2^A(?,DA): \mod A \to \mod A
\]
is nilpotent (\ie some power of it vanishes),
then $\cc_A$ is $\Hom$-finite and $2$-Calabi-Yau.
Moreover, the image $T$ of $A$ in $\cc_A$ is a cluster-tilting
object. The quiver of its endomorphism algebra is obtained
from that of $A$ by adding, for each pair of vertices
$(i,j)$, a number of arrows equal to
\[
\dim \Tor_2^A(S_j, S_i^{op})
\]
from $i$ to $j$, where $S_j$ is the simple right module
associated with $j$ and $S_i^{op}$ the simple left module
associated with $i$.
\end{theorem}

We consider two classes of examples obtained from this
theorem: First, let $\Delta$ be a simply laced Dynkin diagram
and $k\vec{\Delta}$ the path algebra of a quiver with underlying
graph $\Delta$. Let $A$ the {\em Auslander algebra of $k\vec{\Delta}$}, \ie
the endomorphism algebra of the direct sum of a system
of representatives of the indecomposable $B$-modules modulo isomorphism. Then
it is not hard to check that the assumptions of the
theorem hold. The quiver of $A$ is simply the quiver
of the category of indecomposables of $k\vec{\Delta}$
and the minimal relations correspond to its meshes. For
example, if we choose $\Delta=\vec{A}_4$ with the linear
orientation, the quiver of the endomorphism algebra of
the image $T$ of $A$ in $\cc_A$ is the quiver obtained
from Geiss-Leclerc-Schr\"oer's construction at the
end of section~\ref{ss:CY-definition}.

\begin{figure}
\[
\xymatrix@C=0.2cm@R=0.5cm{ & \bullet \ar[dd] &  &  & \circ \ar[lll] \ar[dd] \ar[rrr] & & &
\bullet \ar[dd]  & & & \circ \ar[lll] \ar[dd] \\
     \bullet  \ar[rd] & &  & \circ \ar[lll]|!{[lld];[llu]}\hole \ar[rrr]|!{[ru];[rd]}\hole \ar[rd] & & &
\bullet \ar[rd] & & & \circ \ar[rd] \ar[lll]|!{[lld];[llu]}\hole & &  \\
           & \circ \ar[rrrd] \ar[urr] \ar[rrruu]  & & & \circ \ar[lll] \ar[rrr] & & &
\circ \ar[llluu] \ar[rru]  \ar[rrruu] \ar[llld] \ar[rrrd] \ar[llllu]|!{[lll];[llluu]}\hole & & & \circ \ar[lll] & \\
           & \bullet \ar[u] \ar[d] & & & \circ \ar[u] \ar[d] \ar[lll] \ar[rrr] & & &
\bullet \ar[u] \ar[d] & & & \circ \ar[u] \ar[d] \ar[lll] & \\
           & \circ \ar[rrru]   & & & \circ  \ar[lll] \ar[rrr] &  & &
\circ \ar[lllu] \ar[rrru] & & & \circ \ar[lll] &
}
\]
\caption{The quiver $\vec{A}_4\ten\vec{D}_5$}
\label{fig:a4-ten-d5}
\end{figure}
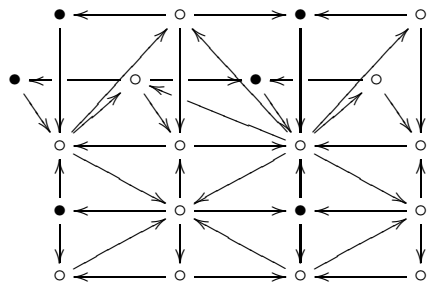

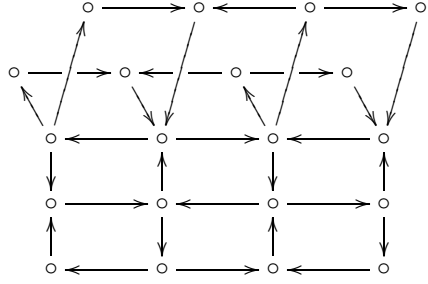
\begin{figure}
\[
\xymatrix@C=0.1cm@R=0.5cm{ & & \circ \ar[rrr] &  &  & \circ \ar[ldd] & & &
\circ \ar[lll] \ar[rrr] & & & \circ \ar[ldd] \\
     \circ \ar[rrr]|!{[dr];[urr]}\hole & &  & \circ \ar[rd] & & &
\circ \ar[lll]|!{[lu];[dll]}\hole \ar[rrr]|!{[rd];[rru]}\hole & & & \circ \ar[rd] & &  \\
           & \circ \ar[lu] \ar[ruu] \ar[d] & & & \circ \ar[lll] \ar[rrr] & & &
\circ \ar[ruu] \ar[lu] \ar[d] & & & \circ \ar[lll] & \\
           & \circ \ar[rrr] & & & \circ \ar[d] \ar[u] & & &
\circ \ar[lll] \ar[rrr] & & & \circ \ar[u] \ar[d] & \\
           & \circ \ar[u] & & & \circ \ar[lll] \ar[rrr] &  & &
\circ \ar[u] & & & \circ \ar[lll] &
}
\]
\caption{The quiver $\vec{A}_4 \square \vec{D}_5$}
\label{fig:a4-square-d5}
\end{figure}

As a second class of examples, we consider an algebra
$A$ which is the tensor product $kQ\ten_k kQ'$ of
two path algebras of quivers $Q$ and $Q'$ without
oriented cycles. Such an algebra is clearly of
global dimension $\leq 2$. Let us assume that
$Q$ and $Q'$ are moreover Dynkin quivers. Then
it is not hard to check that the functor $\Tor_2(?,DA)$
is indeed nilpotent. Thus, the theorem applies. Another
elementary exercise in homological algebra shows that the
space $\Tor_2^A(S_{j,j'}, S_{i,i'}^{op})$ is at most one-dimensional
and that it is non zero if and only if there is an arrow $i\to j$ in $Q$ and
an arrow $i'\to j'$ in $Q'$. This immediately yields the shape
of the quiver of the endomorphism algebra of the image $T$
of $A$ in $\cc_A$: It is the {\em tensor product $Q\ten Q'$}
obtained from the product $Q\times Q'$ by adding an arrow
$(j,j')\to (i,i')$ for each pair of arrows $i\to j$ of $Q$ and
$i'\to j'$ of $Q'$. For example, for suitable orientations
of $A_4$ and $D_5$, we obtain the quiver of figure~\ref{fig:a4-ten-d5}.
Notice that if we perform mutations at the six vertices of the
form $(i,i')$, where $i$ is a sink of $Q=\vec{A}_4$ and $i'$ a source
of $Q'=\vec{D}_5$ (they are marked by $\bullet$), we obtain the
quiver of figure~\ref{fig:a4-square-d5} related to the periodicity conjecture for $(A_4, D_5)$, \cf below.

\section{Application: The periodicity conjecture}
\label{s:periodicity}
Let $\Delta$ and $\Delta'$ be two Dynkin
diagrams with vertex sets $I$ and $I'$. Let $A$ and $A'$ be the
incidence matrices of $\Delta$ and $\Delta'$, \ie if $C$ is the
Cartan matrix of $\Delta$ and $J$ the identity matrix of the
same format, then $A=2 J -C$.
Let $h$ and $h'$ be the Coxeter numbers of $\Delta$ and $\Delta'$.

The {\em $Y$-system of algebraic equations} associated with the
pair of Dynkin diagrams $(\Delta,\Delta')$ is a system of countably
many recurrence relations in the variables $Y_{i,i',t}$, where
$(i,i')$ is a vertex of $\Delta\times\Delta'$ and $t$ an integer.
The system reads as follows:
\begin{equation} \label{eq:Y-system}
Y_{i,i',t-1} Y_{i,i',t+1} =
\frac{\prod_{j\in I} (1+Y_{j,i',t})^{a_{ij}}}{\prod_{j'\in I'} (1+Y_{i,j',t}^{-1})^{a'_{i'j'}}}.
\end{equation}

\begin{periodicity-conjecture} \label{conj:periodicity} All solutions to this
system are periodic in $t$ of period dividing $2(h+h')$.
\end{periodicity-conjecture}

Here is an algebraic reformulation:
Let $K$ be the fraction field of the ring of integer polynomials
in the variables $Y_{i i'}$, where $i$ runs through the
set of vertices $I$ of $\Delta$ and $i'$ through the set of
vertices $I'$ of $\Delta'$. Since $\Delta$ is a tree, the
set $I$ is the disjoint union of two subsets $I_+$ and $I_-$
such that there are no edges between any two vertices of $I_+$
and no edges between any two vertices of
$I_-$. Analogously, $I'$ is the disjoint union of two sets
of vertices $I'_+$ and $I'_-$. For a vertex $(i,i')$ of the
product $I\times I'$, define $\eps(i,i')$ to be $1$ if
$(i,i')$ lies in $I_+\times I'_+ \cup I_-\times I'_-$ and
$-1$ otherwise. For $\eps=\pm 1$, define an automorphism
$\tau_\eps$ of $K$ by $\tau_\eps(Y_{i,i'})=Y_{i,i'}$ if $\eps(i,i')\neq \eps$
and
\begin{equation} \label{eq:automorphism-tau-eps}
\tau_\eps(Y_{i i'}) =
Y_{i i'}^{-1}\prod_{j} (1 + Y_{ji'})^{a_{ij}}  \prod_{j'} (1+ Y_{ij'}^{-1})^{-a'_{i'j'}}
\end{equation}
if $\eps(i,i')=\eps$. Finally, define an automorphism $\phi$ of $K$ by
\begin{equation} \label{eq:automorphism-phi}
\phi = \tau_+ \tau_-.
\end{equation}
Then, as in \cite{FominZelevinsky03b}, it is not hard to check that the periodicity conjecture
holds iff the order of the automorphism $\phi$ is finite and divides
$h+h'$.

The conjecture was formulated
\begin{itemize}

\item by Al.~B.~Zamolodchikov for $(\Delta, A_1)$, where $\Delta$ is
simply laced \cite[(12)]{Zamolodchikov91};

\item by Ravanini-Valleriani-Tateo for $(\Delta, \Delta')$, where $\Delta$
and $\Delta'$ are simply laced or `tadpoles' \cite[(6.2)]{RavaniniTateoValleriani93};

\item by Kuniba-Nakanishi for $(\Delta, A_n)$, where $\Delta$ is not
necessarily simply laced \cite[(2a)]{KunibaNakanishi92}, see also
Kuniba-Nakanishi-Suzuki \cite[B.6]{KunibaNakanishiSuzuki94}; notice
however that in the non simply laced case, the $Y$-system
given by Kuniba-Nakanishi is different from the one above and that they
conjecture periodicity with period dividing twice the sum of the
{\em dual} Coxeter numbers;

\end{itemize}

\noindent It was proved

\begin{itemize}
\item for $(A_n, A_1)$  by Frenkel-Szenes \cite{FrenkelSzenes95}
(who produced explicit solutions) and by Gliozzi-Tateo
\cite{GliozziTateo96} (via volumes of threefolds computed using
triangulations),

\item by Fomin-Zelevinsky \cite{FominZelevinsky03b} for $(\Delta,
A_1)$, where $\Delta$ is not necessarily simply laced (via the
cluster approach and a computer check for the exceptional types; a
uniform proof can be given using \cite{YangZelevinsky08}),

\item for $(A_n, A_m)$ by Volkov \cite{Volkov07}, who exhibited explicit
solutions using cross
ratios, and by Szenes \cite{Szenes09}, who interpreted the system as a
system of flat connections on a graph; an equivalent statement
was proved by Henriques \cite{Henriques07};

\item by Hernandez-Leclerc for $(A_n, A_1)$ using
representations of quantum affine algebras (which yield formulas for
solutions in terms of $q$-characters). They expect to treat $(A_n,
\Delta)$ similarly, \cf \cite{HernandezLeclerc09}.
\end{itemize}

\begin{theorem}[\cite{Keller10a}]
The periodicity conjecture~\ref{conj:periodicity} is true.
\end{theorem}

Let us sketch a proof based on $2$-Calabi-Yau categories (details
will appear elsewhere \cite{Keller10a}): First, using
the folding technique of Fomin-Zelevinsky's \cite{FominZelevinsky03b}
one reduces the conjecture to the case where both $\Delta$ and $\Delta'$
are simply laced. Thus, from now on, we assume that
$\Delta$ and $\Delta'$ are simply laced.

{\em First step:} We choose quivers $Q$
and $Q'$ whose underlying graphs are
$\Delta$ and $\Delta'$. We assume that $Q$ and $Q'$ are {\em alternating},
\ie that each vertex is either a source or a sink. If $i$ is a vertex
of $Q$ or $Q'$, we put $\eps(i)=1$ if $i$ is a source and $\eps(i)=-1$
if $i$ is a sink. For example, we can
consider the following quivers
\begin{align*}
\vec{A}_4 &: \xymatrix{ 1 & 2 \ar[l] \ar[r] & 3 & 4 \ar[l] } \ko\\
\vec{D}_5 &: \raisebox{0.75cm}{\xymatrix@R=0.2cm{   &  &  & 4 \ar[dl] \\ 1 & 2 \ar[l] \ar[r] & 3 & \\
& & & 5. \ar[ul] }}
\end{align*}
We define the {\em square product $Q \square Q'$} to be the quiver
obtained from $Q\times Q'$ by reversing all
the arrows in the full subquivers of the form $\{i\}\times Q'$
and $Q\times \{i'\}$, where $i$ is a sink of $Q$ and $i'$
a source of $Q'$.
The square product of the above quivers $\vec{A}_4$ and $\vec{D}_5$ is
depicted in figure~\ref{fig:a4-square-d5}. The {\em initial $Y$-seed}
associated with $Q\square Q'$ is the pair $\mathbf{y}_0$
formed by the quiver $Q\square Q'$ and the family of variables $Y_{i,i'}$,
$(i,i')\in I\times I'$. We can apply mutations to it using the quiver
mutation rule and the mutation rule for (non tropical) $Y$-variables
given in equation~\ref{eq:Y-variable-mutation}.

{\em A general construction.}
Let $R$ be a quiver and $v$ a sequence of vertices $v_1, \ldots, v_N$
of $R$. We assume that the composed mutation
\[
\mu_v = \mu_{v_N} \ldots \mu_{v_2} \mu_{v_1}
\]
transforms $R$ into itself. Then clearly the same holds for the
inverse sequence
\[
\mu_{v}^{-1} = \mu_{v_1} \mu_{v_2} \ldots \mu_{v_N}.
\]
Now the {\em restricted $Y$-pattern} associated
with $R$ and $\mu_v$ is the sequence of $Y$-seeds obtained from the
initial $Y$-seed $\mathbf{y}_0$ associated with $R$ by applying
all integral powers of $\mu_v$. Thus this pattern is given by
a sequence of seeds $\mathbf{y}_t$, $t\in \Z$, such that $\mathbf{y}_0$
is the initial $Y$-seed associated with $R$ and, for all $t\in \Z$,
$\mathbf{y}_{t+1}$ is obtained from $\mathbf{y}_t$ by the
sequence of mutations $\mu_v$.

{\em Second step.}
For two elements $\sigma$, $\sigma'$ of $\{+, -\}$
define the following composed mutation of $Q\square Q'$:
\[
\mu_{\sigma, \sigma'} = \prod_{ \eps(i)=\sigma, \eps(i')=\sigma'} \mu_{(i,i')}.
\]
Notice that there are no arrows between any two vertices of the
index set so that the order in the product does not matter.
Then it is easy to check that $\mu_{+,+} \mu_{-,-}$ and $\mu_{-,+} \mu_{+,-}$
both transform $Q\square Q'$ into $(Q\square Q')^{op}$ and vice versa. Thus
the composed sequence of mutations
\[
\mu_{\square}=\mu_{-,-} \mu_{+,+} \mu_{-,+} \mu_{+,-}
\]
transforms $Q\square Q'$ into itself.
We define the {\em $Y$-system $\mathbf{y}_\square$} associated with $Q\square Q'$
to be the restricted $Y$-pattern associated with $Q\square Q'$
and $\mu_\square$. As in section~8 of \cite{FominZelevinsky07}, one checks
that the identity $\phi^{h+h'}=\id$ follows if one shows
that the system $\mathbf{y}_\square$ is periodic of period dividing $h+h'$.

{\em Third step.}
One checks easily that we have
\[
\mu_{+-}(Q\square Q') = Q\ten Q' \ko
\]
where the tensor product $Q\ten Q'$ is defined at the end of
section~\ref{ss:2-CY-from-global-dim-2}. For the above quivers
$\vec{A}_4$ and $\vec{D}_5$, the tensor product is depicted
in figure~\ref{fig:a4-ten-d5}. Therefore, the periodicity
of the restricted $Y$-system associated with $Q\square Q'$ and $\mu_\square$
is equivalent to that of the restricted $Y$-system associated
with $Q\ten Q'$ and
\[
\mu_\ten=\mu_{+,-} \mu_{-,-} \mu_{+,+} \mu_{-,+}.
\]

{\em Fourth step.} As we have seen in section~\ref{ss:2-CY-from-global-dim-2},
the quiver $Q\ten Q'$ admits a $2$-Calabi-Yau realization given by the
cluster category $\cc_{kQ\ten kQ'}$ associated with the tensor product
of the path algebras $kQ\ten kQ'$. Let $T$ be the initial
cluster tilting object. By theorem~\ref{thm:2-CY-mutation}, we can define its iterated
mutations
\[
T_t = \mu_{\ten}^t (T)
\]
for all $t\in \Z$. Now we use the

\begin{proposition} None of the endomorphism quivers occurring in the sequence
of mutated cluster-tilting objects joining $T$ to $T_t$ contains loops nor $2$-cycles.
\end{proposition}

Now it follows from theorem~\ref{thm:categorify-mutants} that the quiver
of the endomorphism algebra of $T_t$ is $Q\ten Q'$ and that the
$Y$-variables in the $Y$-seed $\mathbf{y}_{\square, t}$ can be
expressed in terms of the triangulated category $\cc_{kQ\ten kQ'}$ and the
objects $T$ and $T_t$. Thus, it suffices to show that
$T_t$ is isomorphic to $T$ whenever $t$ is an integer multiple
of $h+h'$.

{\em Fifth step.} Let $\tau\ten\id$ denote the auto-equivalence
of the bounded derived category of $kQ\ten kQ'$ given by the total
left derived functor of the tensor product with the bimodule complex
$(\Sigma^{-1} D(kQ))\ten kQ'$, where $D$ is duality over $k$. It
induces an autoequivalence of $\cc_{kQ\ten kQ'}$ which we still
denote by $\tau\ten\id$. Define $\Phi$ to be the autoequivalence
$\tau^{-1}\ten \id$ of $\cc_{kQ\ten kQ'}$.

\begin{proposition} For each integer $t$, the image $\Phi^t(T)$
is isomorphic to $T_t$.
\end{proposition}

The proposition is proved by showing that the indices of the
two objects are equal. This suffices by theorem~\ref{thm:Dehy-Keller}.
Now we conclude thanks to the following categorical periodicity
result:

\begin{proposition} The power $\Phi^{h+h'}$ is isomorphic to the
identity functor.
\end{proposition}

Let us sketch the proof of this proposition: With the natural
abuse of notation, the Serre functor $S$ of the bounded derived
category of $kQ\ten kQ'$ is given by the `tensor product' $S\ten S$
of the Serre functors for $kQ$ and $kQ'$. In the generalized
cluster category, the Serre functor becomes isomorphic to the
square of the suspension functor. So we have
\[
S = S\ten S = \Sigma^2 = \Sigma \ten \Sigma
\]
as autoequivalences of $\cc_{kQ\ten kQ'}$. Now
recall from equation~\ref{eq:tau-Sigma-equals-S} that
$\tau = \Sigma^{-1} S$. Thus we get the isomorphism
\[
\tau \ten \tau = \id
\]
of autoequivalences of $\cc_{kQ\ten kQ'}$. So we have
\[
\tau^{-1} \ten \id = \id \ten \tau.
\]
Now we compute $\Phi^{h+h'}$ by using the left hand side for the
first $h$ factors and the right hand side for the last $h'$ factors:
\[
\Phi^{h+h'} = (\tau^{-1}\ten \id)^h (\id \ten \tau)^{h'}.
\]
But we know from equation~\ref{eq:tau-h-equals-Sigma2}
that $\tau^{-h}=\Sigma^2$. So we find
\[
\Phi^{h+h'} = (\Sigma^2 \ten \id) (\id \ten \Sigma^{-2}) = \Sigma^2 \Sigma^{-2} = \id
\]
as required.

\section{Quiver mutation and derived equivalence}
\label{s:quiver-mutation-and-derived-equivalence}

Let $Q$ be a finite quiver
and $i$ a {\em source} of $Q$, \ie no arrows have target $i$. Then the mutation $Q'$
of $Q$ at $i$ is simply obtained by reversing all the arrows starting
at $i$. In this case, the categories of representations of $Q$ and
$Q'$ are related by the Bernstein-Gelfand-Ponomarev
reflection functors \cite{BernsteinGelfandPonomarev73}
and these induce equivalences in the derived categories \cite{Happel87}. We would
like to present a similar categorical interpretation for mutation
at arbitrary vertices.
For this, we use recent work by Derksen-Weyman-Zelevinsky
\cite{DerksenWeymanZelevinsky08} and a construction due to Ginzburg
\cite{Ginzburg06}. We first recall the classical reflection functors
in a form which generalizes well.

\subsection{A reminder on reflection functors} We keep the above notations.
In the sequel, $k$ is a field and $kQ$ is the path algebra
of $Q$ over $k$. For each vertex $j$ of $Q$, we write
$e_j$ for the lazy path at $j$, the idempotent of $kQ$
associated with $j$. We write $P_j=e_j kQ$ for the
corresponding indecomposable right $kQ$-module,
$\Mod kQ$ for the category of (all)
right $kQ$-modules and $\cd(kQ)$ for its derived
category.

The module categories over $kQ$ and $kQ'$ are linked by
a pair of adjoint functors \cite{BernsteinGelfandPonomarev73}
\[
\xymatrix{
\Mod kQ' \ar[d]<-1ex>_{F_0} \\ \Mod kQ. \ar[u]<-1ex>_{G_0} }
\]
The right adjoint $G_0$ takes a representation $V$ of $Q^{op}$
to the representation $V'$ with $V'_j=V_j$ for $j\neq i$
and where $V'_i$ is the kernel of the map
\[
\bigoplus_{\stackrel{\mbox{\tiny arrows of $Q$}}{i\to j}} V_j  \longrightarrow  V_i
\]
whose components are the images under $V$ of the arrows $i \to j$.

Then the left derived functor of $F_0$ and the right derived
functor of $G_0$ are quasi-inverse equivalences \cite{Happel87}
\[
\xymatrix{
\cd kQ' \ar[d]<-1ex>_{F} \\ \cd kQ. \ar[u]<-1ex>_{G} }
\]
The functor $F$ sends the indecomposable projective $P'_j$, $j\neq i$,
to $P_j$ and the indecomposable projective $P'_i$ to the cone over
the morphism
\[
P_i \to \bigoplus_{\stackrel{\mbox{\tiny arrows of $Q$}}{i\to j}} P_j.
\]
In fact, the sum of the images of all the $P'_j$ has a natural structure
of complex of $kQ'$-$kQ$-bimodules and $F$ can also be described
as the derived tensor product over $kQ'$ with this complex of bimodules.

The example of the following quiver $Q$
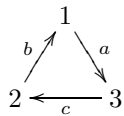
\[
\begin{xy} 0;<0.2pt,0pt>:<0pt,-0.2pt>::
(94,0) *+{1} ="0",
(0,156) *+{2} ="1",
(188,156) *+{3} ="2",
"1", {\ar"0"^b},
"0", {\ar"2"^a},
"2", {\ar"1"^c},
\end{xy}
\]
and its mutation $Q'$ at $1$
\[
Q': 2 \la 1 \la 3
\]
shows that if we mutate at a vertex which is neither a sink nor
a source, then the derived categories of representations of $Q$ and
$Q'$ are not equivalent in general. The cluster-tilted algebra
associated with the mutation of $Q'$ at $1$
(\cf section~\ref{ss:from-cluster-categories-to-cluster-algebras})
is the quotient of $kQ$ by the ideal generated by $ab$, $bc$ and $ca$.
It is of infinite global dimension and therefore not derived equivalent
to $kQ$, either.

Clearly, in order to understand mutation from a
representation-theoretic point of view, more structure is needed. Now
quiver mutation has been independently invented and investigated in
the physics literature, \cf for example equation (12.2) on page 70 in
\cite{CachazoFiolIntriligatorKatzVafa01} (I thank to S.~Fomin for this
reference).  In the physics context, a crucial role is played by the
so-called superpotentials, \cf for example \cite{BerensteinDouglas02}.
This lead Derksen-Weyman-Zelevinsky to their systematic study of
quivers with potentials and their mutations in
\cite{DerksenWeymanZelevinsky08}.  We now sketch their main result.

\subsection{Mutation of quivers with potentials} Let $Q$ be a finite quiver.
Let $\hat{kQ}$ be the {\em completed path algebra},
\ie the completion of the path algebra at the ideal generated by the arrows
of $Q$. Thus, $\hat{kQ}$ is a topological algebra and the paths of $Q$ form
a topologial basis so that the underlying vector space of $\hat{kQ}$ is
\[
\prod_{p \mbox{ \tiny path}} kp.
\]
The {\em continuous Hochschild homology} of $\hat{kQ}$ is the vector
space $\hh_0$ obtained as the quotient of $\hat{kQ}$ by the closure of
the subspace generated by all commutators. It admits a topological
basis formed by the {\em cycles} of $Q$, \ie the orbits
of paths $p=(i|\alpha_n| \ldots |\alpha_1| i)$ of length $n\geq 0$ with
identical source and target under the action of the cyclic group of order $n$.
In particular, the space $\hh_0$ is a product of copies of $k$ indexed
by the vertices if $Q$ does not have oriented cycles.
For each arrow $a$ of $Q$, the {\em cyclic derivative with respect to $a$}
is the unique linear map
\[
\del_a : \hh_0 \to \hat{kQ}
\]
which takes the class of a path $p$ to the sum
\[
\sum_{p=uav} vu
\]
taken over all decompositions of $p$ as a concatenation
of paths $u$, $a$, $v$, where $u$ and $v$ are of length $\geq 0$.
A {\em potential} on $Q$ is an element $W$ of $\hh_0$ whose expansion
in the basis of cycles does not involve cycles of length $0$.

Now assume that $Q$ does not have loops or $2$-cycles.

\begin{theorem}[Derksen-Weyman-Zelevinsky
\protect{\cite{DerksenWeymanZelevinsky08}}] The mutation operation
\[
Q \mapsto \mu_i(Q)
\]
admits a good extension to quivers with potentials
\[
(Q,W) \mapsto \mu_i(Q,W)= (Q', W') \ko
\]
\ie the quiver $Q'$ is isomorphic to $\mu_i(Q)$ if
$W$ is generic.
\end{theorem}

Here, `generic' means that $W$ avoids a certain countable union of
hypersurfaces in the space of potentials. The main ingredient of the
proof of this theorem is the construction of a `minimal model', which,
in a different language, was also obtained in \cite{Kajiura07} and
\cite{KontsevichSoibelman06}.

If $W$ is not generic,
then $\mu_i(Q)$ is not necessarily isomorphic to $Q'$ but
still isomorphic to its $2$-reduction, \ie the quiver obtained from
$Q'$ by removing the arrows of a maximal set of pairwise
disjoint $2$-cycles.  For example, the mutation of the quiver
\begin{equation} \label{eq:3-cycle}
Q:
\raisebox{0.5cm}{\begin{xy} 0;<0.2pt,0pt>:<0pt,-0.2pt>::
(94,0) *+{1} ="0",
(0,156) *+{2} ="1",
(188,156) *+{3} ="2",
"1", {\ar"0"^b},
"0", {\ar"2"^a},
"2", {\ar"1"^c},
\end{xy}
}
\end{equation}
with the potential $W=abc$ at the vertex $1$ is the quiver with potential
\[
Q': 2 \la 1 \la 3 \ko \quad W=0 .
\]
On the other hand, the mutation of the above cyclic quiver $Q$
with the potential $W=abcabc$ at the vertex $1$ is
the quiver
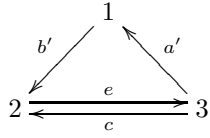
\[
\xymatrix{  & 1 \ar[dl]_{b'} & \\
2 \ar[rr]<0.5ex>^e & & 3 \ar[ll]<0.5ex>^{c} \ar[lu]_{a'}
}
\]
with the potential $ecec+eb'a'$.

Two quivers with potentials $(Q,W)$ and $(Q',W')$ are {\em right
  equivalent} if there is an isomorphism $\phi: \hat{kQ}\to \hat{kQ'}$
taking $W$ to $W'$. It is shown in \cite{DerksenWeymanZelevinsky08}
(\cf also \cite{Zelevinsky07}) that the mutation $\mu_i$ induces an
involution on the set of right equivalence classes of quivers with
potentials, where the quiver does not have loops and does not have a
$2$-cycle passing through $i$.

\subsection{Derived equivalence of Ginzburg dg algebras} We will
associate a derived equivalence of differential graded (=dg) algebras with
each mutation of a quiver with potential. The dg algebra in question is the
algebra $\Gamma=\Gamma(Q,W)$ associated by Ginzburg
with an arbitrary quiver $Q$ with potential $W$, \cf section~4.3 of
\cite{Ginzburg06}. The underlying
graded algebra of $\Gamma(Q,W)$ is the (completed) graded
path algebra of a {\em graded quiver}, \ie a quiver where
each arrow has an associated integer degree. We first
describe this graded quiver $\tilde{Q}$: It has the
same vertices as $Q$. Its arrows are
\begin{itemize}
\item the arrows of $Q$ (they all have degree $0$),
\item an arrow $a^*: j \to i$ of degree $-1$ for each arrow $a:i\to j$ of $Q$,
\item loops $t_i : i \to i$ of degree $-2$ associated with each vertex $i$
of $Q$.
\end{itemize}
We now consider the completion $\Gamma=\hat{k\tilde{Q}}$ (formed in the
category of graded algebras) and endow it with the unique continuous
differential of degree $1$ such that on the generators we have:
\begin{itemize}
\item $da=0$ for each arrow $a$ of $Q$,
\item $d(a^*) = \del_a W$ for each arrow $a$ of $Q$,
\item $d(t_i) = e_i (\sum_{a} [a,a^*]) e_i$ for each vertex $i$ of $Q$, where
$e_i$ is the idempotent associated with $i$ and the sum runs over the
set of arrows of $Q$.
\end{itemize}
One checks that we do have $d^2=0$ and that the homology in degree $0$
of $\Gamma$ is the {\em Jacobi algebra} as defined in
\cite{DerksenWeymanZelevinsky08}:
\[
\cp(Q,W)= \hat{kQ}/(\del_\alpha W \; | \; \alpha\in Q_1) \ko
\]
where $(\del_\alpha W \; | \; \alpha\in Q_1)$ denotes the
closure of the ideal generated by the cyclic derivatives
of the potential with respect to all arrows of $Q$.

Let us consider two typical examples: Let $Q$ be the quiver
with one vertex and three loops labeled $X$, $Y$ and $Z$.
Let $W=XYZ-XZY$. Then the cyclic derivatives of $W$ yield
the commutativity relations between $X$, $Y$ and $Z$ and the
Jacobi algebra is canonically isomorphic to the power series
algebra $k[[X,Y,Z]]$. It is not hard to check that in this
example, the homology of $\Gamma$ is concentrated in degree $0$
so that we have a quasi-isomorphism $\Gamma \to \cp(Q,W)$.
Using theorem~5.3.1 of Ginzburg's \cite{Ginzburg06}, one can show
that this is the case if and only if the full subcategory
of the derived category of the Jacobi algebra formed by the
complexes whose homology is of finite total dimension is
$3$-Calabi-Yau as a triangulated category (\cf also below).

As a second example, we consider the Ginzburg dg algebra associated with
the cyclic quiver~\ref{eq:3-cycle} with the potential $W=abc$.
Here is the graded quiver $\tilde{Q}$:
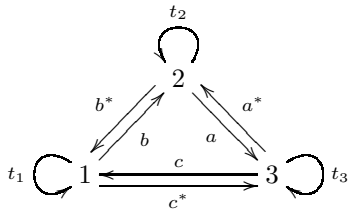
\[
\xymatrix{ & 2 \ar[dr]_a  \ar@<-1ex>[ld]_{b^*} \ar@(ur,ul)[]_{t_2} & \\
1 \ar@(ul,dl)[]_{t_1} \ar[ur]_b \ar@<-1ex>[rr]_{c^*}& &
3 \ar[ll]_c \ar@<-1ex>[ul]_{a^*} \ar@(ur,dr)[]^{t_3} }
\]
The differential is given by
\[
d(a^*)=bc\ko d(b^*)=ca \ko d(c^*)=ab \ko d(t_1)= cc^* - b^* b \ko \ldots .
\]
It is not hard to show that for each $i\leq 0$, we have a canonical isomorphism
\[
H^i(\Gamma) = \cc_{\vec{A}_3}(T,\Sigma^i T)
\]
where $\vec{A}_3$ is the quiver $1 \to 2 \to 3$, $\cc_{\vec{A}_3}$ its
cluster category and $T$ the sum of the images of the modules
$P_1$, $P_3$ and $P_3/P_2$. Since $\Sigma$ is an autoequivalence
of finite order of the cluster category $\cc_{\vec{A}_3}$,
we see that $\Gamma$ has non vanishing homology in infinitely
many degrees $<0$. In particular, $\Gamma$ is not quasi-isomorphic to
the Jacobi algebra.

Let us denote by $\cd \Gamma$ the derived category of $\Gamma$.
Its objects are the differential $\Z$-graded right $\Gamma$-modules and its
morphisms obtained from morphisms of dg $\Gamma$-modules (homogeneous
of degree $0$ and which commute with the differential) by formally
inverting all quasi-isomorphisms, \cf \cite{Keller94}. Let us
denote by $\per \Gamma$ the {\em perfect derived category}, \ie
the full subcategory of $\cd \Gamma$ which is the closure
of the free right $\Gamma$-module $\Gamma_\Gamma$.
under shifts, extensions and passage to direct factors.
Finally, we denote
by $\cd_{fd}(\Gamma)$ the {\em finite-dimensional derived category}, \ie
the full subcategory of $\cd\Gamma$ formed by the dg modules
whose homology is of finite total dimension (!). We recall that
the objects of $\per(\Gamma)$ can be intrinsically characterized
as the {\em compact} ones, \ie those whose associated covariant
$\Hom$-functor commutes with arbitrary coproducts. The objects $M$
of the bounded derived category are characterized by the
fact that $\Hom(P,M)$ is finite-dimensional for each object
$P$ of $\per(\Gamma)$.

The following facts will be shown in \cite{Keller09a}.
Notice that they hold for arbitrary $Q$ and $W$.

1) The dg algebra $\Gamma$ is {\em homologically smooth} (as
a topological dg algebra), \ie it is
perfect as an object of the derived category of
$\Gamma^e=\Gamma\hat{\ten}\Gamma^{op}$.
This implies that the bounded derived category is contained
in the perfect derived category, \cf \eg \cite{Keller08d}.

2) The dg algebra $\Gamma$ is {\em $3$-Calabi-Yau as a bimodule}, \ie
there is an isomorphism in the derived category of $\Gamma^e$
\[
\RHom_{\Gamma^e}(\Gamma, \Gamma^e) \iso \Gamma[-3].
\]
As a consequence, the finite-dimensional derived category $\cd_{fd}(\Gamma)$
is $3$-Calabi-Yau as a triangulated category, \cf \eg
lemma~4.1 of \cite{Keller08d}.

3) The triangulated category $\cd(\Gamma)$ admits a $t$-structure
whose left aisle $\cd_{\leq 0}$ consists of the dg modules
$M$ such that $H^i(M)=0$ for all $i>0$. This $t$-structure
induces a $t$-structure on $\cd_{fd}(\Gamma)$ whose heart is
equivalent to the category of finite-dimensional (and hence nilpotent) modules
over the Jacobi algebra.

The {\em generalized cluster category} \cite{Amiot08} is the
idempotent completion $\cc_{(Q,W)}$ of the quotient
\[
\per(\Gamma)/\cd_{fd}(\Gamma).
\]
The name is justified by the fact that if $Q$ is a quiver
without oriented cycles (and so $W=0$), then $\cc_{(Q,W)}$
is triangle equivalent to $\cc_Q$. In general, the category
$\cc_{(Q,W)}$ has infinite-dimensional $\Hom$-spaces. However,
if $H^0(\Gamma)$ is finite-dimensional, then $\cc_{(Q,W)}$ is
a $2$-Calabi-Yau category with cluster-tilting object $T=\Gamma$,
in the sense of section~\ref{ss:CY-definition}, \cf \cite{Amiot08}
(the present version of \cite{Amiot08} uses non completed
Ginzburg algebras; the complete case will be included in a
future version).

From now on, we suppose that $W$ does not involve cycles of length $\leq 1$.
Then the simple $\tilde{Q}$-modules $S_i$ associated with the
vertices of $Q$ yield a basis of the Grothendieck group
$K_0(\cd_{fd}(\Gamma))$. Thanks to the Euler form
\[
\langle P, M\rangle = \chi(\RHom(P,M)) \ko P\in \per(\Gamma) \ko M\in \cd_{fd}(\Gamma)\ko
\]
this Grothendieck group is dual to $K_0(\per(\Gamma))$ and the
dg modules $P_i = e_i \Gamma$ form the basis dual to the
$S_i$. The Euler form also induces an antisymmetric bilinear
form on $K_0(\cd_{fd}(\Gamma))$ (`Poisson form').
If $W$ does not involve cycles
of length $\leq 2$, then the quiver of the Jacobi algebra is $Q$ and
the matrix of the form on $\cd_{fd}(\Gamma)$ in the basis of the $S_i$ is
given by
\[
\langle S_i, S_j \rangle = |\{\mbox{arrows $j\to i$ of $Q$}\}| - |\{\mbox{arrows $i\to j$ of $Q$}\}.
\]
The dual of this form is given by the map (`symplectic form')
\[
K_0(k) \to K_0(\per(\Gamma\hat{\ten}\Gamma^{op}) = K_0(\per(\Gamma))\ten_\Z K_0(\per(\Gamma))
\]
associated with the $k$-$\Gamma^e$-bimodule $\Gamma$.

Now suppose that $W$ does not involve cycles of length $\leq 2$ and
that $Q$ does not have loops nor $2$-cycles. Then each simple $S_i$
is a ($3$-)spherical object in the $3$-Calabi-Yau category $\cd_{fd}(\Gamma)$,
\ie we have an isomorphism of graded algebras
\[
\Ext^*(S_i, S_i) \iso H^*(S^3, k).
\]
Moreover, the quiver $Q$ encodes the dimensions of the $\Ext^1$-groups
between the $S_i$: We have
\[
\dim \Ext^1(S_i, S_j) = |\{\mbox{arrows $j\to i$ in $Q$}\}|.
\]
In particular, for $i\neq j$, we have either $\Ext^1(S_i, S_j)=0$ or
$\Ext^1(S_j, S_i)=0$ because $Q$ does not contain $2$-cycles.
It follows that the $S_i$, $i\in Q_0$, form a {\em spherical collection}
in $\cd_{fd}(\Gamma)$ in the sense of \cite{KontsevichSoibelman08}.
Conversely, each spherical collection in an ($A_\infty$-) $3$-Calabi-Yau category
in the sense of \cite{KontsevichSoibelman08} can be obtained in this
way from the Ginzburg algebra associated to a quiver with potential.

The categories we have considered so far are summed up in the
sequence of triangulated categories
\begin{equation} \label{eq:ex-seq-tria-cat}
\xymatrix{ 0 \ar[r] & \cd_{fd}(\Gamma) \ar[r] & \per(\Gamma) \ar[r] & \cc_{(Q,W)} \ar[r] & 0}.
\end{equation}
This sequence is `exact up to factors', \ie the left hand term is a thick
subcategory of the middle term and the right hand term is the idempotent
closure of the quotient.
Notice that the left hand category is $3$-Calabi-Yau, the middle one does
not have a Serre functor and the right hand one is $2$-Calabi-Yau if the
Jacobi algebra is finite-dimensional.

We keep the last hypotheses: $Q$ does not have loops or $2$-cycles
and $W$ does not involve cycles of length $\leq 2$. Let $i$ be a vertex
of $Q$ and $\Gamma'$ the Ginzburg algebra of the mutated quiver
with potential $\mu_i(Q,W)$. The following theorem
improves on Vit\'oria's \cite{Vitoria09}, \cf also
\cite{MukhopadhyayRay04} \cite{BerensteinDouglas02}.

\begin{theorem}[\cite{KellerYang09}] \label{thm:KellerYang}
There is an equivalence of
derived categories
\[
F : \cd(\Gamma') \to \cd(\Gamma)
\]
which takes the dg modules $P'_j$ to $P_j$ for $j\neq i$ and
$P'_i$ to the cone on the morphism
\[
P_i \longrightarrow \bigoplus_{\stackrel{\mbox{\tiny arrows}}{i \to j}} P_j.
\]
It induces triangle equivalences $\per(\Gamma')\iso \per(\Gamma)$ and
$\cd_{fd}(\Gamma')\to\cd_{fd}(\Gamma)$.
\end{theorem}

In fact, the functor $F$ is given by the left derived functor
of the tensor product by a suitable $\Gamma'$-$\Gamma$-bimodule.

By transport of the canonical $t$-structure on $\cd(\Gamma')$,
we obtain new $t$-structures on $\cd(\Gamma)$ and $\cd_{fd}(\Gamma)$.
They are related to the canonical ones by a tilt
(in the sense of \cite{HappelReitenSmaloe96}) at the simple
object corresponding to the vertex $i$. By iterating mutation
(as far as possible), one obtains many new $t$-structures on
$\cd_{fd}(\Gamma)$.

In the context of Iyama-Reiten's \cite{IyamaReiten06}, the dg algebras
$\Gamma$ and $\Gamma'$ have their homologies
concentrated in degree $0$ and the equivalence
of the theorem is given by the tilting module constructed in [loc. cit.].

\subsection{A geometric illustration} \label{ss:geometric-illustration}
To illustrate the sequence~\ref{eq:ex-seq-tria-cat}, let us consider the
example of the quiver
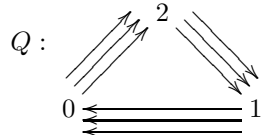
\[
Q: \raisebox{0.5cm}{\xymatrix{ & 2  \ar[rd]<2ex> \ar[rd]<1ex> \ar[rd] & \\
0 \ar[ru] \ar[ru]<1ex> \ar[ru]<2ex> & & 1\ar[ll] \ar[ll]<1ex> \ar[ll]<2ex>
}
}
\]
where the arrows going out from $i$ are labeled $x_i$, $y_i$, $z_i$, $0\leq i\leq 2$,
endowed with the potential
\[
W=\sum_{i=0}^2 (x_i y_i z_i - x_i z_i y_i).
\]
Let $\Gamma'$ be the non completed Ginzburg algebra associated
with $(Q,W)$.
Let $p: \omega \to \P^2$ be the canonical bundle on $\P^2$. Then we
have a triangle equivalence
\[
\cd^b(\coh(\omega)) \iso \per(\Gamma')
\]
which sends $p^*(O(-i))$ to the dg module $P_i$, $0\leq i\leq 2$.
Under this equivalence, the subcategory $\cd_{fd}\Gamma'$ corresponds to
the subcategory $\cd^b_Z(\coh(\omega))$ of complexes of coherent
sheaves whose homology is supported on the zero section $Z$ of the
bundle $\omega$. This subcategory is indeed Calabi-Yau of dimension
$3$. Bridgeland has studied its $t$-structures \cite{Bridgeland05}
by linking them to $t$-structures and mutations (in the sense of
Rudakov's school) in the derived category of coherent sheaves
on the projective plane. In this example, the Ginzburg
algebra has its homology concentrated in degree $0$. So it is
quasi-isomorphic to the Jacobi algebra. The Jacobi algebra is
infinite-dimensional so that the generalized cluster category
$\cc_{(Q,W)}$ is not $\Hom$-finite.  The category $\cc_{(Q,W)}$
identifies with the quotient
\[
\cd^b(\coh(\omega))/\cd^b_Z(\coh(\omega))
\]
and thus with the bounded derived category $\cd^b(\coh(\omega\setminus
Z))$ of coherent sheaves on the complement of the zero section of
$\omega$.  In order to understand why this category is `close to
being' Calabi-Yau of dimension $2$, we consider the scheme $C$
obtained from $\omega$ by contracting the zero section $Z$ to a point
$P_0$. The projection $\omega \to C$ induces an isomorphism from
$\omega \setminus Z$ onto $C\setminus\{P_0\}$ and so we have an
equivalence
\[
\cc_{(Q,W)} \iso \cd^b(\coh(\omega\setminus Z)) \iso
\cd^b(\coh(C\setminus\{P_0\})).
\]
Let $\hat{C}$ denote the completion of $C$ at the singular point
$P_0$. Then $\hat{C}$ is of dimension $3$ and has $P_0$ as its unique
closed point. Thus the subscheme $\hat{C}\setminus \{P_0\}$ is of
dimension $2$. The induced functor
\[
\cc_{(Q,W)} \iso \cd^b(\coh(C\setminus\{P_0\})) \to
\cd^b(\coh(\hat{C}\setminus\{P_0\}))
\]
yields a `completion' of $\cc_{(Q,W)}$ which is of `dimension $2$' and
Calabi-Yau in a generalized sense, \cf \cite{ChuangRouquier07}.

\subsection{Ginzburg algebras from algebras of global dimension $2$}
Let us link the cluster categories obtained from algebras of
global dimension~$2$ in section~\ref{ss:2-CY-from-global-dim-2}
to the generalized cluster categories $\cc_{(Q,W)}$
constructed above.

Let $A$ be an algebra given as the quotient $kQ'/I$ of the
path algebra of a finite quiver $Q'$ by an ideal $I$ contained
in the square of the ideal $J$ generated by the arrows of $Q'$.
Assume that $A$ is of global dimension $2$ (but not necessarily
of finite dimension over $k$). We construct a quiver with
potential $(Q,W)$ as follows: Let $R$
be the union over all pairs of vertices $(i,j)$ of a set
of representatives of the vectors belonging to a basis of
\[
\Tor_2^A(S_j, S_i^{op}) = e_j(I/(IJ+JI)) e_i .
\]
We think of these representatives as `minimal relations' from $i$ to $j$,
\cf \cite{Bongartz83}.
For each such representative $r$, let $\rho_r$ be a new
arrow from $j$ to $i$. We define $Q$ to be obtained from $Q'$ by
adding all the arrows $\rho_r$ and the potential is given by
\[
W=\sum_{r\in R} r \rho_r .
\]
Recall that a {\em tilting module} over an algebra $B$
is a $B$-module $T$ such that the total derived functor
of the tensor product by $T$ over the endomorphism
algebra $\End_B(T)$ is an equivalence
\[
\cd(\End_B(T)) \iso \cd(B).
\]
The second assertion of part a) of the following theorem
generalizes a result by Assem-Br\"ustle-Schiffler
\cite{AssemBruestleSchiffler08}.

\begin{theorem}[\protect{\cite{Keller09a}}]
\begin{itemize}
\item[a)]
The category $\cc_{(Q,W)}$ is triangle equivalent
to the cluster category $\cc_A$. This equivalence takes
$\Gamma$ to the image $\pi(A)$ of $A$ in $\cc_A$ and thus
induces an isomorphism from the Jacobi algebra $\cp(Q,W)$ onto
the endomorphism algebra of the image of $A$ in $\cc_A$.
\item[b)]
If $T$ is a tilting module over $kQ''$ for a quiver
without oriented cycles $Q''$ and $A$ is the
endomorphism algebra of $T$, then
$\cc_{(Q,W)}$ is triangle equivalent to $\cc_{Q''}$.
\end{itemize}
\end{theorem}

\subsection{Cluster-tilting objects, spherical collections, decorated representations}
We consider the setup of theorem~\ref{thm:KellerYang}.
The equivalence $F$ induces equivalences $\per(\Gamma')\to \per(\Gamma)$
and $\cd_{fd}(\Gamma')\to \cd_{fd}(\Gamma)$. The last equivalence takes
$S'_i$ to $\Sigma S_i$ and, for $j\neq i$, the module $S'_j$ to
$S_j$ if $\Ext^1(S_j,S_i)=0$ and, more generally, to the middle term $FS'_j$ of
the universal extension
\[
S_i^d \to FS'_j \to S_j \to \Sigma S_i^d \ko
\]
where the components of the third morphism form a basis of $\Ext^1(S_j, S_i)$.
This means that the images $FS'_j$, $j\in Q_0$, form the
{\em left mutated spherical collection} in $\cd_{fd}(\Gamma)$ in the
sense of \cite{KontsevichSoibelman08}.

If $H^0(\Gamma)$ is finite-dimensional so that $\cc_{(Q,W)}$ is
a $2$-Calabi-Yau category with cluster-tilting object, then the
images of the $FP_j$, $j\in Q_0$, form the mutated cluster-tilting
object, as it follows from the description in lemma~\ref{lemma:simple-mutation}.

Now suppose that $H^0(\Gamma)$ is finite-dimensional.
Then we can establish a connection with decorated representations
and their mutations in the sense of \cite{DerksenWeymanZelevinsky08}:
Recall from [loc. cit.] that a decorated representation of $(Q,W)$
is given by a finite-dimensional (hence nilpotent) module $M$
over the Jacobi algebra and a collection of vector spaces $V_j$
indexed by the vertices $j$ of $Q$. Given an object $L$ of $\cc_{(Q,W)}$,
we put $M=\cc(\Gamma, L)$ and, for each vertex $j$, we choose a
vector space $V_j$ of maximal dimension such that the triangle
\begin{equation}
T_1 \to T_0 \to L \to \Sigma T_1
\end{equation}
of lemma~\ref{lemma:Keller-Reiten} admits a direct factor
\[
V_j\ten P_j \to 0 \to \Sigma V_j\ten P_j \to \Sigma V_j \ten P_j.
\]
Let us write $GL$ for the decorated representation thus constructed.
The assignment $L\mapsto GL$ defines a bijection between isomorphism
classes of objects $L$ in $\cc_{(Q,W)}$ and right equivalence
classes of decorated representations. It is compatible with
mutations: If $(M',V')$ is the image $GL'$ of an object $L'$
of $\cc_{(Q',W')}$, then the mutation $(Q,W,M,V)$ of
$(Q',W',M',V')$ at $i$ in the sense of \cite{DerksenWeymanZelevinsky08}
is right equivalent to $(Q,W,M'',V'')$ where $(M'',V'')=GFL'$
and $F$ is the equivalence of theorem~\ref{thm:KellerYang}.



\def\cprime{$'$} \def\cprime{$'$}
\providecommand{\bysame}{\leavevmode\hbox to3em{\hrulefill}\thinspace}
\providecommand{\MR}{\relax\ifhmode\unskip\space\fi MR }
\providecommand{\MRhref}[2]{%
  \href{http://www.ams.org/mathscinet-getitem?mr=#1}{#2}
}
\providecommand{\href}[2]{#2}

\end{document}